\documentclass[12pt]{article}
\usepackage{graphicx}
\usepackage{subfig}
\usepackage{extsizes}
\usepackage{epsfig}
\usepackage[numbers,sort&compress]{natbib}
\usepackage{citehack}
\usepackage{float}
\usepackage{amsmath}
\usepackage{amssymb}
\usepackage{multicol}
\usepackage{graphicx}
\usepackage{subfig}
\usepackage{extsizes}
\usepackage{epsfig}

\newcommand{\eps}{\varepsilon}
\newcommand{\vte}{\vartheta}
\newcommand{\vfi}{\varphi}

\newcommand{\R}{\mathbb {R}}

 \newcommand{\sx}{separatrix}   

\newtheorem{theorem}{ТTheorem}
\newtheorem{lemma}{ЛLemma}[section]
\newtheorem{corollary}{СCorollary}[section]
\newtheorem{proposition}{Proposition}[section]

\makeatletter \@addtoreset{equation}{section} \makeatother
\date{}
\title{Averaging method for systems  with separatrix crossing}\author{Anatoly NeishtadtАН$^{1,2}$  \\
$^1$ Loughborough University, Loughborough, LE11 3TU, UK\\
 $^2$ Space Research Institute, Moscow, 117997, Russia}

\begin{document}
\maketitle

\begin{abstract}
The averaging method  provides a powerful tool for studying evolution in near-integrable systems. Existence of separatrices  in the phase space of the underlying integrable system is an obstacle for application of standard results that justify using of  averaging. We establish estimates that allow to use averaging method when the underlying integrable system is a system with one rotating phase, and the evolution leads to separatrix crossings.     
\end{abstract}


\section{Introduction}
\label{intro} 
An averaging method (see, e.g., \cite{bm}) is a powerful tool for study a long-term evolution  in  systems which are small perturbations of integrable systems. Many applications of  this method are for one-frequency systems (also called systems with one fast rotating phase).  In these cases in the phase space of the corresponding unperturbed system there is a  domain  foliated by closed trajectories - invariant circles. Averaging of perturbations over these circles provides a closed system for an approximate description of perturbed dynamics in this domain. Typically, in systems under consideration there are several  such domains. These domains are bounded by surfaces on which this foliation has singularities. Classical results justifying averaging method \cite{bm} guarantee its applicability for description of evolution  not too close to these separating surfaces.  However, it is rather typical that evolution leads to crossing of these  surfaces. The goal of this paper is to provide justification of a modified   version of the averaging method for description of such an evolution.

\begin{figure}
 \begin{center}
            \includegraphics[scale=0.4, angle=0.2]{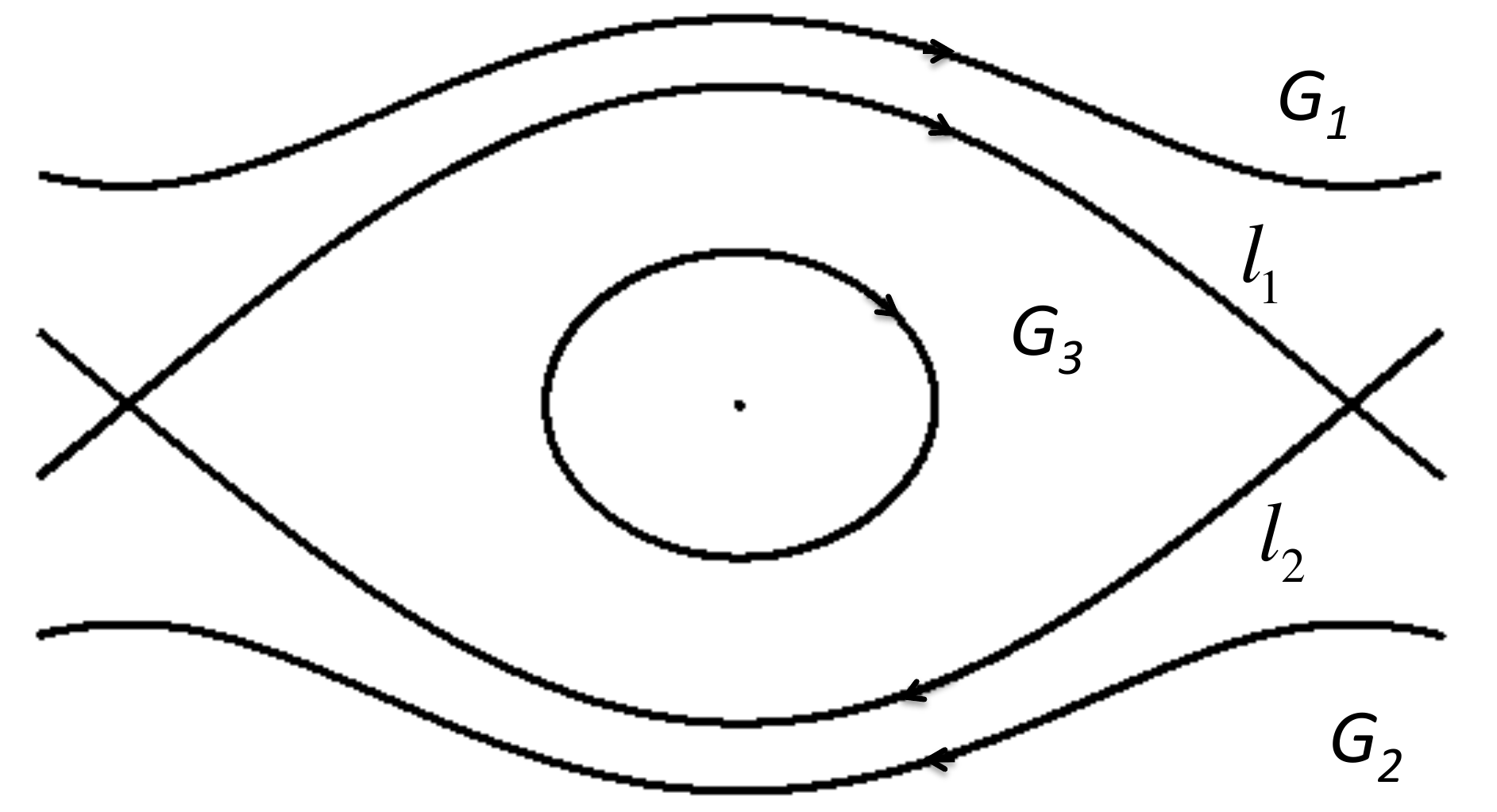}
            \end{center}
           \caption{Phase portrait of a pendulum.}
            \label{pendulum} 
\end{figure}

A paradigmatic example of problems considered in this paper is a pendulum under the action of perturbations, e.g., of  a small friction, a small constant torque, and a slow change of its length. An unperturbed pendulum could be in one of three regimes of motion: it could rotate in one or other direction, or oscillate. In the phase portrait of the pendulum these three regimes are demarcated by separatrices (Fig. \ref{pendulum}). Motion of the pendulum evolves slowly under the action of perturbations. In the process  of this evolution the pendulum can change the regime of its motion. In the phase plane the phase point crosses an instant separatrix of the unperturbed  pendulum. Evolution of energy far from the separatrix can be described by the averaging method. Classical results justifying  this method \cite{bm} are not applicable  in the case of crossings of a separatrix. Moreover, this crossing leads to a remarkable probabilistic scattering.  Initial data for different outcomes of separatrix crossings are mixed, and it is reasonable to consider each outcome as a random event with a definite  probability.  This probabilistic approach was first described in a similar problem   in \cite{lif} and then independently in \cite{goldr}. 

A natural way to describe evolution in the considered system is to use averaging method up to arrival to the separatrix, to calculate probabilities  of capture into different domains at the separatrix, and to use the averaging method starting from the separatrix in the domain in which the system continues its  evolution.   In the current paper we justify such an  approach for a rather wide class of one-frequency systems that change a qualitative character  of their motion in the process of evolution. The obtained estimates of the accuracy of the averaging method are sharp. 

Part of results of this paper was announced (without proofs) in \cite{nei_chaos} (see also \cite{akn}, Subsection 6.1.10) on the basis of estimates in \cite{nei_dd}.



\section{Averaging method  and  averaging theorem  for the  separatrix crossing} \label{S2}
In this section an averaging theorem is formulated that  justifies  the averaging method  for description of the  separatrix crossing. The proof is based on  propositions given in Sections  3, 4, 5. There are probability phenomena due to  separatrix crossing. Hence the recipe of the averaging method here includes calculations of the corresponding probabilities, and the  averaging theorem justifies these calculations. All considerations are  for systems of the form (\ref{perturbed}) below. We explain in the Appendix  relation of this form to the general form of one-frequency systems with separatrix crossings. 

\subsection{Outline of the problem}\label{ss2.1}
\begin{figure}
 \begin{center}
            \includegraphics[scale=0.65, angle=0.0]{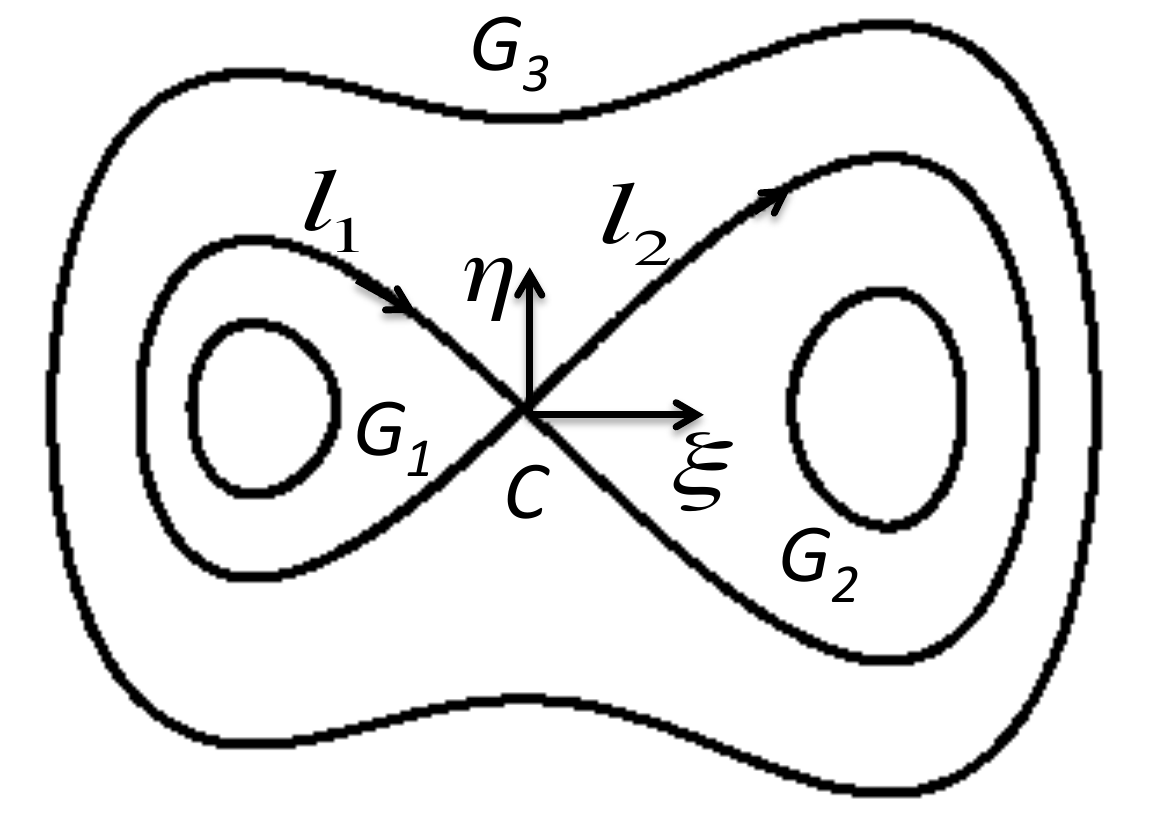}
            \end{center}
           \caption{Phase portrait of the unperturbed system.}
            \label{unperturbed_plane} 
\end{figure}

We  consider   systems described by  differential equations of the form
 \begin{eqnarray} \label{perturbed}
 \dot q&=&\frac{\partial E}{\partial p}+\eps f_1, \, \dot p=-\frac{\partial E}{\partial q}+\eps f_2, \, \dot z=\eps f_3\, ,
 \\
        E&=&E(p,q,z),\, f_i=f_i(p,q,z,\varepsilon),i=1,2,3,\, (p,q)\in \R^2,z\in \R^{l-2}\,. \ \nonumber
       \end{eqnarray}
       Here $\varepsilon >0$  is a small parameter  characterising the rate of evolution. For $\varepsilon=0,\,  z={\rm const}$ we have {\it an unperturbed system} for $p,q$, which is a Hamiltonian system with one degree of freedom. The function $E$ is an unperturbed  Hamiltonian, and the functions $\varepsilon f_i$ are the perturbations.
       It is supposed, that there are separatrices in the phase portrait of the unperturbed system Fig. \ref{unperturbed_plane}. In the course of evolution the projection of the phase point onto the plane $(p,q)$ crosses a separatrix. 

       Far from the separatrices instead of  $(p,q)$ it is possible to use the variables $h=E$ and $\varphi$, where  $\varphi$ is  ``the angle'' (from the pair ``action-angle''  variables \cite{arn_1} of the unperturbed system). Then for $h,z,\varphi$  we get the perturbed system having the standard form of  system with one rotating phase \cite{bm}: in this system $h,z$ are called {\it slow variables}, $\vfi$ is {\it the rotating phase}. It is a classical result  that the averaged with respect to $\varphi$ system describes the evolution of $h,z$ far from separatrices with accuracy $O(\varepsilon)$ during the time interval of order $1/\varepsilon$. 
        At the {separatrices} the frequency of the      
 unperturbed motion vanishes, and also the equations in variables               
 $h$,$z$,$\varphi $ have singularities.    In a region that             
 includes a separatrix, the conditions of  the classical theorem about accuracy of the averaging  \cite{bm}  fail and        
 the applicability of the averaging method for the description of the           
 evolution near the {separatrices} requires a justification. 
 
\begin{figure}
 \begin{center}
            \includegraphics[scale=0.5]{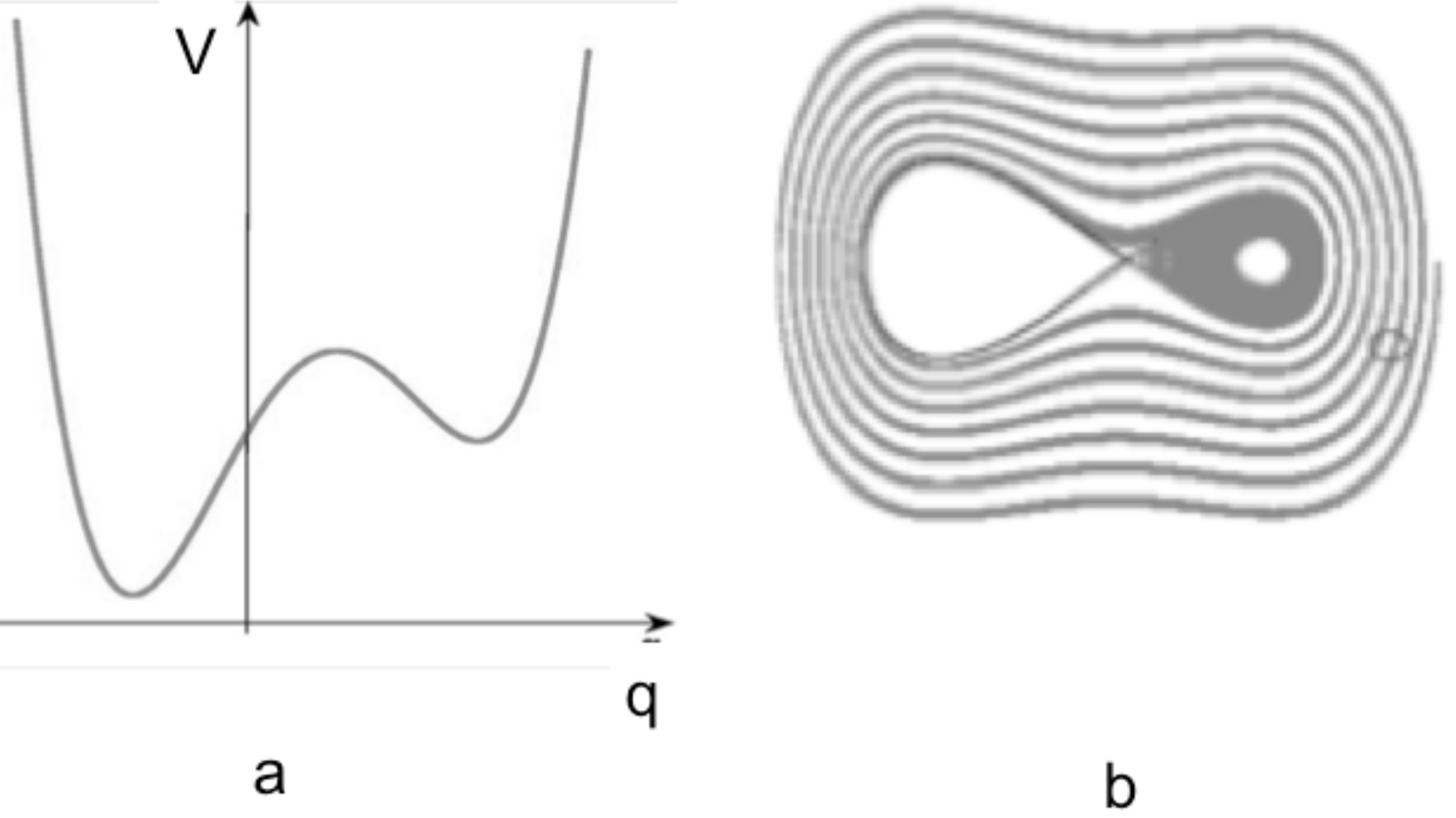}
            \end{center}
           \caption{a) Double-well potential. b) Effect of a small friction. }
            \label{double_well} 
\end{figure}

       Separatrix crossing leads to  probability phenomena \cite {lif,goldr, arn_2,gurevich}.  As a simple example let as consider the motion of a particle in one dimension in double-well potential, Fig. \ref{double_well}a, perturbed by a small, of order  $\eps$, dissipation \cite {arn_2}.  Phase portrait of the perturbed system is shown in Fig.  \ref{double_well}b, where the initial conditions for the capture into  the region, surrounded by the  right separatrix loop, are shaded.
The shaded strips far from the saddle have  width of order $\eps$ and  form a spiral with a step of order $\eps$. Therefore small, of order $\eps$, change of the initial conditions can change the result of evolution. As the initial conditions are always known only with some finite accuracy, the deterministic approach to the problem fails when  $\eps \to 0$. But it is possible to define in some natural way  and to calculate the probabilities of  capture into different regions after the separatrix crossing  \cite{arn_2}.

For systems of such types, a procedure of  an approximate description of the evolution consists of using  the averaged system up to the separatrix and calculation of  the probability of capture into one or another region on the separatrix. It will be seen, that for  majority of  initial conditions this procedure describes the behaviour of the slow variables with accuracy $O(\eps\ln\eps)$ during  time of order $1/\eps$. The measure of the ``bad'' set  of initial conditions, for which this description is not valid, tends to $0$ faster, than any given power of $\eps$ as $\eps\to 0$.  The general formula for the probability of the capture into one or another region (in the sense of definition in \cite {arn_2}) also will be proved.
       
       
       \subsection{Formulation of  Hypotheses}\label{ss2.2}

 System (\ref{perturbed}) is considered for $(p,q,z)\in D, \mid\eps\mid<\eps_1$, where $ D$ is a bounded domain in $\R^l, \ \eps_1= {\rm const}$.  We denote $B$ the projection of  $D$ onto $z$-space, and $G=G(z)$  the section of $ D$ by the two-dimensional plane $z= {\rm const}$. It is supposed that each domain  $G(z) \subset \R^2$ is composed of the whole trajectories of the unperturbed system. It is supposed that the following assumptions  are satisfied.
 
 \medskip
 $\bf{A}$. The function $E$ is of smoothness $C^3$, and the functions $f_i$ are of smoothness $C^2$ with respect to $p, q, z$. Functions $f_i$ have one continuous derivative with respect to $\eps$.
 
 \medskip
 $\bf{B}$. For $z\in B$ the phase portrait of the unperturbed system in the domain $G(z)$ has the form shown in  Fig. \ref {unperturbed_plane}. The unstable stationary point $C$ is a non-degenerate saddle point. The separatrices $l_1=l_1(z)$ and $l_2=l_2(z)$ divide the unperturbed phase portrait into three regions   $G_{\nu}=G_{\nu}(z), \, \nu=1, 2, 3$.
 In what follows we assume that  the Hamiltonian $E$ is normalised in such a way, that $E=0$ at the saddle point $C$, and therefore, on the separatrices. Then $E>0$ in the region $G_3$,  $E<0$  in the regions $G_1$ and $G_2$. We denote $l_3=l_1 \cup l_2$.
 
 \medskip
 $\bf {C}$. Introduce  quantities
 \begin{eqnarray} \label{2.2}
       \Theta_\nu(z)&=&-\oint_{l_\nu}\left(\frac{\partial E}{\partial q}f_1^0+\frac{\partial E}{\partial p}f_2^0+\frac{\partial E}{\partial z}f_3^0\right)dt\,, \  \nu=1, 2,\\
       \Theta_3(z)&=&\Theta_1(z) +\Theta_2(z)\, .\nonumber
       \end{eqnarray}
        Integrals in (\ref{2.2}) are calculated along the unperturbed separatrices parametrized by the time $t$ of the unperturbed motion; $f_{i}^0=f_i(p,q,z,0),\, i=1,2,3$.   Integrals  (\ref {2.2}) are improper, because the motion along a separatrix takes infinite time. Our normalisation of $E$ guarantees  the convergence of the integrals as it is proved at the end of this section. We assume that the values $\Theta_\nu, \,\nu=1, 2, 3$ are different from zero.  In what follows, for certainty, the values  $\Theta_\nu,\, \nu=1, 2, 3$ are supposed to be   positive.
       
       \medskip       
       Let us explain the meaning of the condition $\bf {C}$. In the region $G_\nu$ for small $|E|>\eps$, in the perturbed motion,  a phase point 
       makes  rounds that  are close to the unperturbed separatrix $l_\nu$. The change of the value of $E$ during one  such round is close to the value $-\eps\Theta_\nu$ 
       . Therefore, for  phase points  with small $|E|$,  condition $\bf C$ ensures an approach the separatrix in the region $G_3$ and a departure from the separatrix in the regions $G_1$ and $G_2$.
       The convergence of  integrals (\ref{2.2}) is a corollary of the following assertion.
         \begin{lemma} 
        \label{L2.1}
        The first derivatives of the function $E$ with respect to $p, q, z$ vanish at the point $C$.
                   \end{lemma} 

{\it Proof}. The derivatives with respect to $p, q $ vanish at the point $C$ because the point $C$ is an equilibrium position of the unperturbed system. Let $p_C(z), q_C(z)$ be  coordinates of the point $C$. The Hamiltonian is normalised by the condition $E(p_C(z), q_C(z),z)\equiv0$. Calculating  the derivative of this equality with respect to $z$ and taking into account that  $\frac{\partial E}{\partial p}=\frac{\partial E}{\partial q} =0$ at the point $C$, we get that $\frac{\partial E}{\partial z}=0$ at the point $C.\hskip 7cm  \square$

\bigskip
Lemma \ref{L2.1} implies that  
$$
 (\partial E/\partial\alpha)/\sqrt{(\partial E/\partial p)^2+(\partial E/\partial q)^2}, \alpha=p, q, z,
 $$
 tend to  finite limits as a point $(p, q)$ tends to the point $C$ along a separatrix. Let us use in the integrals (\ref {2.2}) the arc length along the separatrix as a new independent variable. Then  integrands  do not have  singularities, and therefore integrals (\ref {2.2}) converge. Moreover,  $\Theta_\nu$ are smooth functions of $z$.

  \subsection{Averaged system}\label{ss2.3}
  Let us define the averaged system separately for each region $G_\nu$ first. 
   Let
  $$
  \varSigma_\nu=\{(h, z)\colon z\in B, h=E(p, q, z), (p,q)\in G_\nu(z)\}, \ \nu=1, 2, 3 \,.
  $$
  The averaged in the  region $G_{\nu}$ system is, by definition, the following system of differential equations in $\Sigma_\nu$: 
  \begin{eqnarray} \label{2.3}                                                             
   \dot h&=&{\eps \over T}\oint_{E=h}( \frac{\partial E} {\partial            
     q}f_1^0+\frac{\partial E} {\partial p}f_2^0+ \frac{\partial E}                
     {\partial z}f_3^0)dt \label{avs}\, ,\\ 
     \dot z&=&{\eps \over                 
     T}\oint_{E=h}f_3^0 dt\, . \nonumber                                            
 \end {eqnarray}   
       Here integrals are calculated along the level line  $E=h$ of the Hamiltonian situated in the domain  $G_\nu(z)$. This level line  is parametrized by the time $t$ of the unperturbed motion along it, and  $T=T(h, z)$ is the period of this motion. To write down this averaged system we  calculate the rate of changing of $h=E(p, q, z)$ and $z$ in the perturbed system, and then  average the obtained expressions over $t$ along the level line $E=h$ (for $\eps=0$ in the arguments of $f_i$). This averaging is equivalent to the averaging over the angular variable $\varphi$ discussed in Subsection \ref{ss2.1}.
       
\begin{figure}
 \begin{center}
            \includegraphics[scale=0.5]{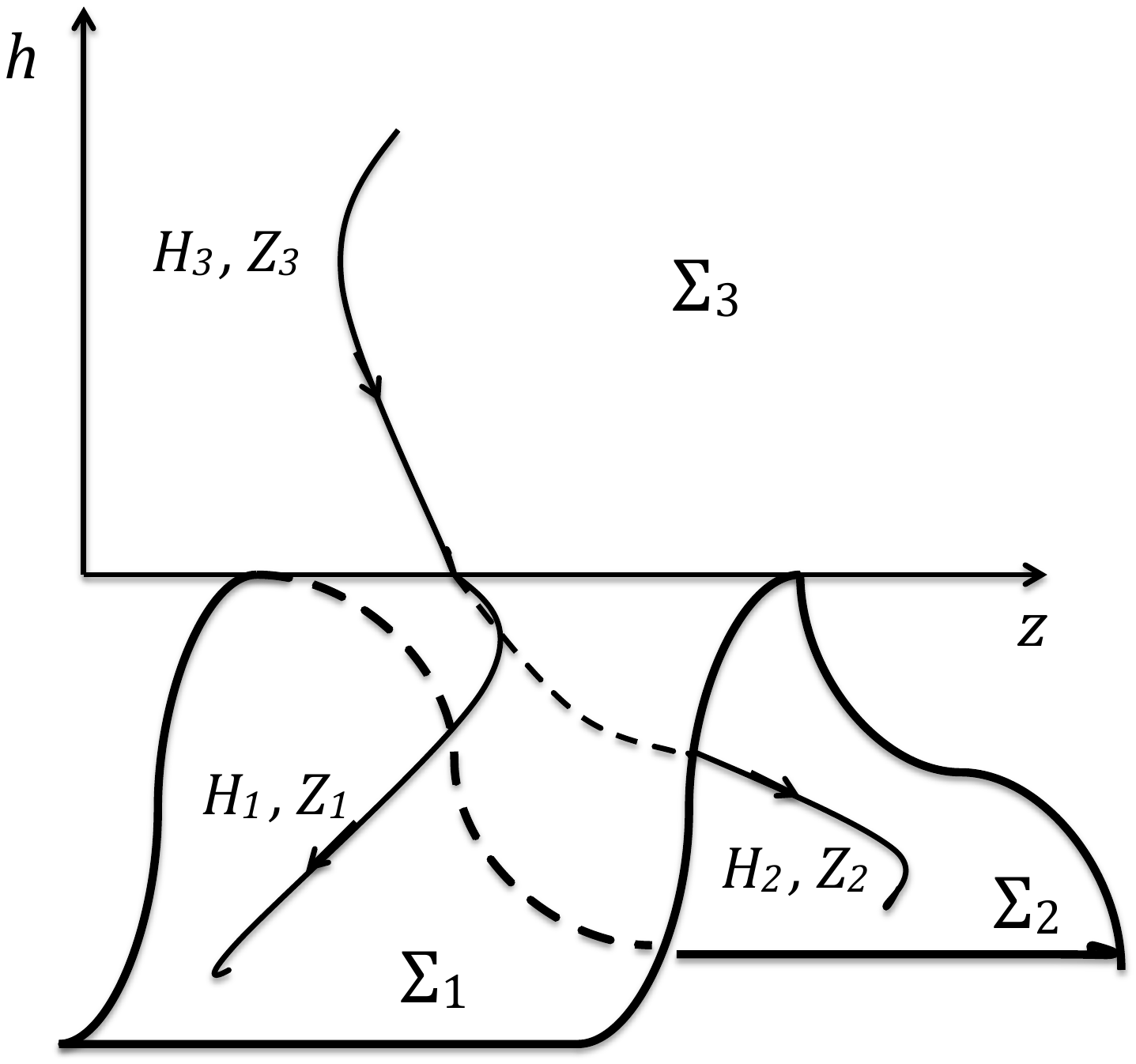}
            \end{center}
            
           \caption{The phase space of the averaged system.}
            \label{glued_space} 
\end{figure}

  The period $T$ grows proportionally  to  $-\ln|h|$  in the principal  approximation  as $h\to 0$ (see Lemma \ref{3.3}). When $h=0$ it is reasonable to  extent  the definition of the  right hand sides of (\ref{2.3}) by continuity, putting
       $$
       \dot h\mid_{h=0}\,=0, \  \dot z\mid_{h=0}\,=\eps f_{3C}^0\, ,
       $$
 where  $f_{3C}^0$ is the value of the function $f_3^0$ at the point $C$. Now we can combine three averaged systems in different regions into one ``whole'' averaged system. The phase space of the ``whole'' averaged system is a singular  manifold, glued of  three parts $\Sigma_1, \Sigma_2$ and $\Sigma_3$ along the set $\{h=0\}$, Fig. \ref{glued_space}.  We will call  the set  $\{h=0\}$ the separatrix for the averaged system.

 According  to condition $\bf {C}$,   \ $\dot h<0$ in the averaged system for small $h\ne 0$ in all regions.
 
  {\bf Definition.} A solution $H(\tau ),\; Z(\tau ),\;                
   \tau =\eps t$, of the averaged system in $G_{\nu }$ such that            
   $H(\tau ) \to 0$, $Z(\tau )\to z_*$ as $ \tau \to \tau{}_*-0$ (for           
   the region $G_3$) or $\tau \to \tau{}_*+0$ ( for the regions $G_{1,2}$) is           
   called a {\em solution crossing the separatrix at the point $z_*$ at the         
   moment $\tau{}_*$.}
                                                             
 \smallskip 
  The following lemma is  proved  in Subsection \ref{ss3.6}.                          
 \begin{lemma} \label{L2.2}                                                     
   a) For any $z_*\in B$, $\tau{}_*$, and $\nu =1,2,3$ there exists              
  a unique solution of the averaged system in $G_{\nu }$ crossing the           
   separatrix at the point $z_*$ at the moment $\tau{}_*$.                      
                                                                                
   \noindent b). If $z_0 \in B$, then for small enough $|h_0|$ the              
   solution $H(\tau ),\; Z(\tau )$ of the averaged system with the              
   initial conditions $H(\tau{}_0)=h_0, Z(\tau{}_0)=z_0$ crosses the            
   separatrix at some moment $\tau =\tau{}_*$ ($\tau{}_*>\tau_0$,  if             
   $h_0>0$ and $\tau{}_*<\tau{}_0$ if $h_0<0$).                                 
  \end{lemma}                                                                   
  Take any point $(p_0,q_0,z_0)\in D$ such that  $(p_0,q_0)\in               
  G_3(z_0)$.  Denote  $h_0=E(p_0,q_0,z_0)$. Consider the solution $               
  (H_3(\tau ),\; Z_3(\tau ))$ of the averaged system in $G_3$ with              
  the initial condition $(h_0,z_0)$ at $\tau =0$ (Fig. \ref{glued_space}). Suppose that             
  this solution crosses the \sx\ at some $\tau= \tau{}_*$, i.e.                 
  $H_3(\tau{}_*)=0,\; Z_3(\tau{}_*)=z_*$. According to   Lemma \ref{L2.2},
  we can consider the solution  $(H_{\nu}(\tau ), Z_{\nu}(\tau))$ of the system averaged in $G_{\nu}$,  $\nu =1,2$ with the initial              
  condition $(0,z_*)$ at $\tau =\tau{}_*$. This solution is well defined 
  for  $\tau$ close enough to $\tau{}_*$, $\tau >\tau{}_*$.
    For an approximate  description of the behaviour of the values $E,\;z$ in the perturbed            
  system (\ref{perturbed}) we use the solution $(H_3,\; Z_3)$ for $0\le        
  \eps t\le \tau{}_*$, and the solution $(H_{\nu}, Z_{\nu})$ for $        
  \eps t > \tau{}_*$ with $ \nu = 1$ or $2$,  if the phase point has been captured  into  the region          
  $G_{\nu}$  after the \sx \ crossing. We define                            
  for $\tau \le \tau{}_* $ the  functions $H_\nu ,\; Z_\nu ,\ \nu =1,2$, by the       
  relations $H_\nu (\tau )=H_3(\tau ),\; Z_\nu (\tau )=Z_3(\tau )$.            
  These functions $(H_\nu , Z_\nu )$ are called the solutions of the               
  averaged system with the initial condition $(h_0,z_0)$ at $t=0$.              
                                                                                
  We attribute the  probability $\Theta_{\nu}(z_*)/\Theta{}_3(z_*)$           
  to the capture into $G_\nu$ of the initial point  $(p_0,q_0,z_0)$. The meaning of this definition of the          
  probability will be clear from the results of Subsection~\ref{s2.est}. The             
  function $P_\nu$ defined by the formula                                          
 \begin{equation}                                                               
 P_{\nu}(z)=                                                                    
  \Theta_{\nu}(z)/\Theta_3(z) \label{prob.def}                                 
 \end{equation}                                                                 
  will be called the probability of the                                         
  capture into $G_\nu ,\, \nu =1,2$, at the \sx{}.                                 
                                                                                
  We will need a lemma which allows to estimate a distance        
  between two solutions of the averaged system with  initial                 
  conditions near the \sx .                                                     
                                                                                
 \begin{lemma}                                                                  
 \label{L2.3}                                                                   
  Let two solutions of the averaged system,                                     
  $(H_\nu(\tau),Z_\nu(\tau))$ and                                               
  $({{H}'_\nu(\tau)},{{Z}'_\nu(\tau))},\ 0\le\tau\le K$, be               
  given. Suppose that for some $\tau_0\in [0,K]$ and some $\delta>0$ these      
  solutions satisfy the following condition:                                    
 $$|H_\nu(\tau_0)|+|{H}'_\nu(\tau_0)|+|Z_\nu(\tau_0)-{Z}'_\nu(\tau_0)|
<                                                                               
  \delta . $$                                                                   
 If $\delta$ is small enough, then for $0\le\tau\le K$ the following            
 estimate is valid:                                                             
 $$|H_\nu(\tau)-{H}'_\nu(\tau)|+|Z_\nu(\tau)-{Z}'_\nu(\tau)|<O\left(\delta+
{                                                                               
  \delta|\ln\delta|\over                                                          
   1+|\ln|H_\nu(\tau)||}\right)\,.$$                                                     
 \end{lemma}                                                                    
 The proof is given in Subsection  \ref{ss3.6}.                                               
 \smallskip                                                                     
                                                                                
 In what follows,  the action $I$ of the unperturbed              
 system  will be important. The action $I=I(h,z)$ of the                       
 unperturbed trajectory $E=h$ in the region $G_\nu$ is the area                 
 enclosed by this trajectory, divided by $2\pi$.
 We have \cite{LL}
 \begin{equation}                                                               
 \frac{\partial I}{\partial h} = \frac{1}{2\pi} T\left(h,z\right), \label{didh} 
 \end{equation} 
 \begin{equation}                                                               
 \frac{\partial I}{\partial z} = -\frac{1}{2\pi} \oint _{E=h} \frac{\partial    
 E}{\partial z} dt\, . \label{didz}                                                
 \end{equation} 
 
 With the aid of the formulas (\ref{didh}) and (\ref{didz})
 the rate of change of $I$ along a trajectory                                   
 of the averaged system is found to be:                                         
 \begin{eqnarray}                                                               
   2\pi {dI\over dt}&=&\eps \oint_{E=h}( {\partial E\over \partial          
     q}f_1^0+{\partial E\over \partial p}f_2^0+ {\partial E\over                
     \partial z}f_3^0)dt \nonumber \\                                           
 &-&{\eps \over T}\oint_{E=h} {\partial E\over                              
     \partial z}dt \ \oint_{E=h}f_3^0 dt\,. 
      \label{2.7}                         
  \end{eqnarray}                                                                
  A corollary of the formula (\ref{didz}) when $h\to 0$ is the                  
  following useful formula for the areas $S_{1,2}=S_{1,2}(z)$ of             
  the regions $G_{1,2}$ and the area $S_3=S_3(z)$ of  the region $G_1\cup            
  G_2$:                                                                         
  \begin{equation}                                                              
    {\partial S_\nu \over \partial z} =- \oint_{l_{\nu}} {\partial              
      E\over \partial z}dt, \nu =1,2,3  \,.
 \end{equation}

 From (\ref{2.2})  for the averaged system in  the region         
 $G_\nu$ we get:                                                             
 \begin{equation}                                                               
   \lim_{h\to 0}{d \over d t}(2\pi I-S_\nu)=-\eps \Theta{}_\nu (z),         
   \ \nu =1,2,3  \,.
  \end{equation}                                                                
  A consequence of this equation is the above-mentioned property               
  (see Lemma ~\ref{L2.2}) that for solutions of the averaged system with            
  small $|h|$ the arrival at the separatrix takes a finite time (in            
  the region $G_3$ this time is positive, in  the regions $G_{1,2}$ it is negative). 
                                                                                
 \subsection{Estimates in the averaging method\label{s2.est}}                      
 \markboth{Separatrix-Crossing Orbits}{Estimates}                               
                                                                                
 Let a point $\hat M_0 = (\hat p_0, \hat q_0,\hat z_0)$ belong to the         
 region $D$, and let $\hat I_0,\hat \varphi_0$ be the  values of the            
 action-angle variables $I,\varphi$ at                          
 this point.  The following sets are well defined and lie in $D$ for small                                                          
 enough $\delta $:                
 \begin{eqnarray}                                                               
   U^\delta&=&\{ p,q,z \ : \ |z-\hat z_0|<\delta ,|I-\hat I_0 |<\delta          
   , |\varphi -\hat \varphi_0|<\delta \}\,, \\ \nonumber W^\delta&=&\{             
   p,q,z \ : \ |z-\hat z_0|<\delta ,|I-\hat I_0 |<\delta \}\,.\nonumber            
  \end{eqnarray}                                                                
  Denote $\hat h_0 = E\left(\hat p_0, \hat q_0,\hat z_0 \right)$. We assume that   
  solutions of the averaged system with                                         
  initial data $(\hat h _0,\hat z_0)$ are                                       
  well defined for $0 \le \tau \le K$ and cross the separatrix at some $ \tau         
  =\hat\tau_*$,  $z=\hat z_*$. Denote:                                
 \begin{itemize}                                                                
  \item $(p(t),q(t),z(t))$ the solution of the perturbed system                 
  (\ref{perturbed}) with initial data $(p_0,q_0,z_0) \in U^\delta $ at $t=0$,   
                                                                                
  \item $h(t)=E(p(t),q(t),z(t))$ the value of the Hamiltonian along this solution, 
    
     \item      $ h_0=h\left( 0 \right) $,                
                                                                                
 \item $(H_\nu ,Z_\nu ),\ \nu =1,2$, the solutions of the averaged               
  system with the initial condition $ ( h_0 ,z_0)$ at $\tau =0$,                
                                                                                
  \item $\tau_*$ the moment of \sx \ crossing for the                           
  solutions $(H_\nu ,Z_\nu )$.                                                  
 \end{itemize}                                                                  
                                                                                
  \noindent The value $\delta $ is supposed to be  small enough so that the               
  solutions $(H_\nu ,Z_\nu ),\ \nu =1,2$, are well defined for $0\le\tau          
  \le K$, and $\tau_* <K$.  Fix any natural number $r\ge 2$.  In what follows $K_i$ (and afterwards $k_i,             
  c_i, d_i, \nu_i$) are positive constants, i.e. values independent of         
  $\eps ,\delta$ and initial conditions $(p_0,q_0,z_0 )\in W^\delta$.       
  The appearance  of $K_i$ in some relation 
  is equivalent to the assertion that there exists $K_i$ satisfying this relation for small      
  enough $\eps >0,\delta>0, \eps<\delta^2$ (and similarly for           
  other constants). The following theorem summarises the principal              
  features of the averaging method for separatrix-crossing orbits:              
 \begin{theorem} \label{T2.1}                                                      
   There exists a representation $U^\delta=U_1^\delta\cup U_2^\delta \cup     
   v $ with the following properties.                                           
                                                                                
   \noindent I.\ If $(p_0,q_0,z_0 ) \in U_\nu^\delta,\ \nu=1,2$, then the behaviour of     
   $E,z$ in the perturbed system is described approximately by the              
   solution $(H_\nu ,Z_\nu)$ of the averaged system, and the following          
   estimates hold                                                         
 \begin{eqnarray}                                                               
   |h(t)-H_\nu (\eps t)|+|z(t)-Z_\nu (\eps t)|&<&K_1\eps \quad    \mbox{for} \quad     
   0\le \eps t\le\tau_*,\\
   \nonumber |h(t)-H_\nu (\eps                   
   t)|+|z(t)-Z_\nu (\eps t)|&<&K_1\eps +{K_2 \eps |\ln              
     \eps |\over 1+|\ln |H_\nu (\eps t)||} \quad    \mbox{for} \quad     \tau_*\le \eps           
   t\le K \,. \nonumber                                                             
 \end{eqnarray}                                                                 
 For $0\le \eps t\le \tau_*-K_3\eps |\ln\eps| $ the point               
 $(p(t),q(t))$ moves in the region $G_3(z(t))$, while for $\tau_*               
 +K_3\eps |\ln\eps |\le \eps t\le K$ it moves in the region         
 $G_\nu (z(t))$.                                                                
    
    \medskip                                                                            
 \noindent II.                                                                  
 $$\left | {\hbox{\em mes} \, U_{\nu}^{\delta}\over \mbox{\em mes} \, U^{\delta} }-         
 {\Theta_\nu(\hat z_*)\over\Theta_3 (\hat z_*)}\right | <K_4\left(\delta +{\eps         
   |\ln \eps|\over \delta }\right).$$                                              
   
   \medskip                                                                               
 \noindent III.\ $\hbox{\em mes}\,  v < k_5\eps^r\delta^{-1}\hbox{\em          
   mes} \, U^\delta$.                                                        
  
  \medskip                                                                              
 \noindent Here $\hbox{\em mes}\, (\cdot )$ is the standard phase volume in            
 $\R^l$ 
.                                                      
 \end{theorem}                                                                  
 This theorem will be proved by means of a series of propositions               
 established in the following three Sections (Propositions \ref{Pr2.1}, \ref{Pr2.2},  \ref{Pr2.3}).                
 \medskip                                                                       
                                                                                
 It is natural to consider the relative measure of the set of                   
 points from a small neighbourhood of $\hat M_0$ that will be                    
 captured into the region $G_\nu , \nu=1,2$ for small $\eps $ as the           
 value at the point $\hat M_0$ of the probability density of                    
 capture into $G_\nu$. This approach is formalised as follows (cf. \cite{arn_2}).                                                              
  
  \medskip                                                                              
 {\bf Definition } The value at the point $\hat M_0$ of {\em the                
   probability density of capture into $G_\nu\, ,\nu=1,2\,,$} (or, for               
   brevity, {\em the probability of capture of $\hat M_0$ into $G_\nu$})        
   is                                                                           
 \begin{equation}                                                               
   Q_\nu(\hat M_0)=\lim_{\delta \to 0}\lim_{\eps\to                         
     0}{\hbox{mes}\, U_\nu^\delta\over \hbox{mes}\,U^\delta}\label{eq:def}\,  .          
 \end{equation}                                                                 
  
  \begin{corollary}  
  \label{col_propability}                                                                            
  The probability of capture of the point $\hat   
   M_0$ into the region $G_\nu$ is given by the formula                               
 \begin{equation}                                                               
   Q_\nu(\hat M_0)={\Theta_\nu(\hat z_*)\over\Theta_3(\hat z_*)},\              
   \nu=1,2\label{eq:cor}\,.                                                        
 \end{equation}   
 \end{corollary}                                                              
 By means of the function $P_\nu$ defined by equation (\ref{prob.def}) --         
  the                                                                    
 probability of the capture into $G_\nu$ at the separatrix -- the last           
 formula can be rewritten in the form $$Q_\nu (\hat M_0)=P_\nu(\hat             
 z_*),\ \nu =1,2\,.$$  
 
 \bigskip
 {\bf Remarks} 
 
 \medskip
{\bf 1.}   The formulation of the problem of separatrix-crossing has been described       
 in the context of an eight-figure separatrix that is being approached by          
 orbits of the perturbed and averaged systems. It will be clear from           
 the nature of the estimates in the succeeding sections,  that the basic results embodied in           
 Theorem \ref{T2.1} can be carried over to the other geometric pictures           
 of separatrix-crossing that are possible in $\R^2$. In fact, the phase space    
need                                                                            
  not be                                                                        
 $\R^2$. It could, for example, be a cylinder or a sphere. The important          
 hypotheses are those demanding the non-degeneracy of the saddle                 
 equilibrium and the non-vanishing of the numbers $\Theta _\nu,\;\nu =           
 1,2,3$ (see Appendix). These conditions may be relaxed: they are needed only            
 when the orbit approaches the separatrix. During the evolution prior           
 to that time, these conditions are  not needed. Metamorphoses of the          
 phase portrait may take place as long as the phase point is                    
 far from separatrices when this happens. Situations wherein the                 
 crossing of a separatrix occurs for values of $z$ for which $\Theta _          
 \nu = 0$ or at which the non-degeneracy condition for the saddle   fails are viewed as         
 degenerate, and the estimates of the averaging method will in general          
 be poorer in these cases. Examples of separatrix-crossing near a             
 ``newborn'' saddle are  considered in \cite {L_Pesci, Kiselev_Nonlinearity, Haberman_Chaos, Haberman_PhysD}.                                      
                                                                                
    
   \medskip                                                                              
 {\bf 2.} The conclusion that the error in the use of the averaging method is            
 $O(\eps \ln \eps)$, which follows from Theorem \ref{T2.1},               
 cannot be improved. This follows from asymptotic formulas  for  change of            
 the adiabatic invariant at  a separatrix  in Hamiltonian systems (\cite{Timofeev}, \cite{Cary}, \cite{nei_fplazmy},              
 \cite{nei_pmm}) and also for systems perturbed by a weak dissipation \cite {Bourland}.

 \medskip                                                                               
{\bf 3.}  The assertion of the Theorem~\ref{T2.1} remains valid if the
standard volume     
 $\hbox{mes}\,(\cdot)$ is replaced by any measure in $\R^l$ that has a             
 smooth density, independent of $\eps$, with respect to the 
 standard volume.   
 The formula for the probability given by (\ref{eq:def}) and                    
 (\ref{eq:cor}) remains valid if $U^\delta$ is any domain with                  
 piecewise-smooth boundary having diameter $\delta$ (but in this case           
 the estimate in the right hand side in part II of the Theorem ~\ref{T2.1} may              
 not be satisfied).  

\medskip  
{\bf 4.} For the example of   motion of a particle in one dimension in double-well potential perturbed by a small dissipation the formula for probability is given in  \cite {arn_2}. The  proof is contained  in \cite{brin_freidlin}.

   \medskip                                                                             
 {\bf 5.} A different approach for introducing probability was suggested in \cite{wolansky,  Freidlin}. White noise of order $\eps \delta$ was added to the                  
 right-hand side of equations (\ref{perturbed}), in the case when the         
 parameter $z$ is absent from the problem. In this problem,                     
 capture into one or another region becomes a genuinely random                  
 event. Again taking the limit of the probability of capture as                 
 $\eps \to 0$ (first) and $\delta \to 0$, one recovers the same             
 formula for the probability as that found in Subsection \ref{s2.est} above.   
 
  \medskip                                                                             
 {\bf 6.} D.V.Anosov has suggested yet another   approach for introducing probability in the considered problem\footnote {This was a comment in a meeting of the Moscow   Mathematical Society.}. Denote $\varkappa_\nu(\eps_0)$ the measure of the set of values $\eps\in (0, \eps_0]$ such that $\hat M_0\in U_\nu^\delta$ for these values of $\eps$, $\nu=1,2$. Then we define the probability  of capture of  $\hat M_0$ into  $G_\nu$ as $\lim_{\eps_0\to 0} \varkappa_\nu(\eps_0)/ \eps_0$.  This definition was not discussed in a literature. One can show that, e.g. for the case when $f_1=f_2=0$, $f_3={\rm const}$ in equations  (\ref{perturbed})  this   probability is again given by  formula (\ref{eq:cor}).   It looks plausible that this is the case for the general form of  system  (\ref{perturbed}) as well. 
 
 {\bf 7.} An important open question is when   results of consecutive crossings of separatrices can be considered as statistically independent. This is not  the case, for example, if  (\ref {perturbed}) is a Hamiltonian system with slowly varying parameter, and $\Theta_1\equiv\Theta_2$ \cite{cary_scodje}.  However, it looks as a reasonable hypothesis that typically there is  such an independence.

  \subsection{Derivation of the estimates in the averaging  method}
  \label{ss2.5} 
  
  In this subsection Theorem \ref{T2.1} is derived from several principal Propositions describing approach the separatrix, passage through its narrow neighbourhood, and departure from the separatrix. These Propositions are proved in Sections  \ref{S3}, \ref{S4}, \ref{S5}. Approach the separatrix is described by the following assertion.
  
 \begin{proposition}
 \label{Pr2.1}
  
 For all initial conditions $(p_0, q_0, z_0)\in W^\delta$ and for $t>0$ while $H_3(\eps t)\ge k_1\eps$ the following estimates are valid:
 $$
|h(t)-H_3(\eps t)|+|z(t)-Z_3(\eps t)|=O(\eps), \  \frac{1}{2}H_3(\eps t)\le h(t)\le 2H_3(\eps t) \, .
$$
 \end{proposition}
Departure from the separatrix in the regions $G_{1,2}$ \, is described by an  analogous assertion. 
                                                                       
 \begin{proposition}
 \label{Pr2.2}
 Let  a  moment of time $t'\in[0, K/\eps]$   exist such that  $h(t')=-k_2\eps, (p(t'), q(t'))\in G_\nu(z(t')),\  \nu=1, 2, \  |\eps t'-\tau_*|+|Z(t')-Z_3(\tau_*)|<k_3^{-1}$.
 Then for $t'\le t\le K/\eps$ the behaviour of $E, z$ is approximately described by the solution $(H'_\nu, Z'_\nu)$ of the averaged in $G_\nu$ system with the initial condition $h(t'), z(t')$ at $\tau=\eps t'$ as follows:
  $$
|h(t)-H'_\nu(\eps t)|+|z(t)-Z'_\nu(\eps t)|=O\left(\frac{\eps\ln\eps}{1+|\ln|H'_\nu(\eps t)||}\right) \!, \, 2H'_\nu(\eps t)\le h(t)\le  \frac{1}{2} H'_\nu(\eps t)
$$
 \end{proposition}
 It follows from Proposition \ref {Pr2.1} that there exists a moment of time $t_-=t_-(p_0, q_0, z_0, \eps)$ such that  at this moment $h(t)=2k_1\eps$ for the first time. Let  $t_+=t_+(p_0, q_0, z_0, \eps)$ be the moment of time such that  at this moment $h(t)=-k_2\eps$ for the first time at the segment $[0, K/\eps]$. The moment of time $t_+$ is defined, in general, not for all initial conditions.
 The following assertion describes  passage through a narrow neighbourhood of the separatrices during the interval of time $t_-\le t\le\ t_+$.     
                                               
  \begin{proposition}
 \label{Pr2.3}
 There exists  a representation $W^\delta=W^\delta_3\cup w$ such that  ${\rm mes}\,  w=O(\eps^r)\,  {\rm mes}\, W^\delta$ and for $(p_0, q_0, z_0) \in W_3^\delta$ the moment  of time  $t_+$ is well defined and            $t_+=t_-+O(\ln\eps)$. For
  $t_-\le t\le\ t_+$  and for $\nu=1, 2$ the following  estimate holds: 
  $$
|h(t)|+|H_\nu(\eps t)|+|z(t)-Z_\nu(\eps t)|=O(\eps)\,.
$$
\end{proposition}

\medskip
Let $W_\nu^\delta,\,  \nu=1, 2$, be the sets of points $(p_0, q_0, z_0) \in W_3^\delta$ such that $(p(t_+), q(t_+)) \in G_\nu(z(t_+))$.

 \begin{proposition}
 \label{Pr2.4}
 The measure of the set $W_\nu^\delta$ is  estimated as 
  $$
       {\rm mes}\, W_\nu^\delta=\int_{W^\delta}P_\nu(Z_3(\tau_*))dp_0dq_0dz_0+O(\eps\ln\eps\, \delta^{-1}\, {\rm mes} \, W^\delta), \ \nu=1, 2\,.
$$
Here $(H_3(\tau), Z_3(\tau))$ is the solution of the averaged system with the initial condition $H_3(0)=E(p_0, q_0, z_0), Z_3(0)=z_0$, and $\tau_*$ is the moment of separatrix crossing for this solution, $dz_0$ is the standart volume element in $\R^{l-2}$.
 \end{proposition}

Denote $U_\nu^\delta=U^\delta\cap W_\nu^\delta, \ v =U^\delta\cap w$.
    
  \begin{proposition}
 \label{Pr2.5}
 The measure of the set $U_\nu^\delta$ meets the estimate of part  $I\!I$ in  Theorem \ref {T2.1}. 
  \end{proposition}
 For $0\le t \le t_+$  the estimates in  Theorem  \ref{T2.1}  hold due to the estimates of  Propositions \ref{Pr2.2} and \ref{Pr2.3}. From these estimates we get also that it is possible to choose $t_+$ as $t'$ in  Propositions \ref {Pr2.2}. Then at $t=t_+$ the distance between  solutions $(H'_\nu, Z'_\nu)$ and  $(H_\nu, Z_\nu)$ is $O(\eps)$. Therefore from Lemma \ref{L2.3} and Proposition \ref{Pr2.2} we get the estimates in  Theorem \ref{T2.1}  for $t_+\le t\le K/\eps$. This completes the proof of  Theorem \ref{T2.1}.
 
 \medskip
  To proof  Corollary \ref{col_propability}   it is enough to consider a cover of $U^\delta$ by the union  of sets
 \begin{eqnarray}                                                               
  U_{j_1, j_2, j_3}^{\kappa} &=& \{p, q, z\colon |z-\hat z_{j_1}|\le\kappa, |I-\hat I_{j_2}|\le\kappa, |\varphi-\hat\varphi_{j_3}|\le\kappa\}  ,\\       
 {\hat z_{j_1}} &=& \hat z_0+2\kappa j_1, \  \hat I_{j_2}=\hat I_0+2\kappa j_2, \  \hat \varphi_{j_3}=\hat \varphi_0+2\kappa j_3  \, .     
    \nonumber                                                             
 \end{eqnarray} 
 Here $j_2, j_3$ are integer numbers, and $j_1$ is an integer $(l-2)$-dimensional vector. Then one  should apply   Theorem \ref{T2.1} to those of sets   $U_{j_1, j_2, j_3}^\kappa$, which intersect $U^\delta$. Now proceed to the limit first as $\eps\to 0$ and then as $\kappa\to 0$. We get the formula
 $$
 \lim_{\eps\to0}\hbox{mes}\, U_\nu^\delta=\int_{U^\delta}P_\nu(Z_3(\tau_*))dp_0dq_0dz_0
 $$                                                       
which implies  Corollary \ref{col_propability}. 

\section{ Estimates of the accuracy of the averaging method up to  separatrix} \label{S3}   
 
 In this section the proof of   Proposition \ref{Pr2.1}   is given. The proof of   Proposition \ref{Pr2.2} is entirely analogous,  and only a sketch of it is given. We use several lemmas, which are formulated in Subsections \ref{3.1}, \ref{3.2}. This lemmas are proven in Subsections \ref{ss3.4}, \ref{ss3.5}. 
 
 Only  motion in the region $G_3$ is considered in the principal part of this section. Thus we will omit index $\nu=3$ at solution $(H_\nu, Z_\nu)$. The following notation is used: $I(h, z)$ - the value of the ``action" variable for the trajectory $E=h$ of the unperturbed system, $j(t)=I(h(t), z(t)), J(\eps t)=I(H(\eps t), Z(\eps t))$.
 
 Below  $z\in B-c_1^{-1}, \, (p,q,z)\in D-c_1^{-1}$ in all estimates, and the constant $c_1$ is chosen in such a way that $3c_1^{-1}$-neighbourhood of the set $\{z\colon z=Z_\nu(\tau), \nu=1, 2; \tau\in [0, K]\}$ belongs to $B$, and  $3c_1^{-1}$-neighbourhood of the set 
 $$
 \{(p, q, z)\colon z=Z_\nu(\tau), E(p, q, z)=H_\nu(\tau), \  \nu=1, 2;\, \tau\in [0, K]\}
 $$
  belongs to $D$.
  
    \subsection{Lemmas on unperturbed motion}
  \label{ss3.1}
  Let $\varphi=\varphi(p, q, z, \eps)$ and $ \psi=\psi(p, q, z, \eps)$ be  smooth functions, and $\varphi$ vanish at the saddle point $C$ identically with respect to $z,\eps$. Let $0<h<1/2, z\in B-c_1^{-1}$.
   
   \begin{lemma} 
        \label{L3.1}
          \begin{equation}  \label{3.1}     
        \oint_{E=h}\varphi dt= \oint_{l_3}\varphi dt+O(h\ln h)\,.
         \end{equation} 
           \end{lemma} 
           \begin{corollary}
         For $0<h<c_2^{-1}$ the following estimate is valid:
          \begin{equation}   \label{3.2}                                                              
    \oint_{E=h}\left(\frac{\partial E}{\partial q}f_1^0+\frac{\partial E}{\partial p}f_2^0+\frac{\partial E}{\partial z}f_3^0\right)dt=-\Theta_3+O(h\ln h)<-c_3^{-1} \,.      
 \end{equation} 
 \end{corollary}

 \begin{lemma} 
        \label{L3.2}
           \begin{eqnarray}      \label{3.3}     
      (1) \ \oint_{E=h}|\varphi|dt=O(1);\phantom{*****} \ (2) \ \oint_{E=h}|\psi|dt=O(\ln h) ;\\
    (3)\  \frac{\partial}{\partial h} \oint_{E=h}\psi dt=O(h^{-1}); \  \ \ 
     (4) \ \frac{\partial}{\partial z} \oint_{E=h}\psi dt=O(\ln h) ; \   \nonumber \\
      (5) \ \frac{\partial}{\partial h} \oint_{E=h}\varphi dt=O(\ln h); \phantom{**}\ (6) \ \frac{\partial}{\partial z} \oint_{E=h}\varphi dt=O(1). \   \nonumber  
         \end{eqnarray}    
           \end{lemma} 
         \begin {corollary} 
         \begin {eqnarray}   \label{3.4}                                                              
    \frac{\partial I}{\partial h}=\frac{1}{2\pi}T=O(\ln h), \quad \frac{\partial I}{\partial z}=-\frac{1}{2\pi}\oint_{E=h}\frac{\partial E}{\partial z}=O(1), \\
     \frac{\partial^2 I}{\partial h^2}=O(h^{-1}), \quad  \frac{\partial^2I}{\partial h\partial z}=O(\ln h), \quad  \frac{\partial^2I}{\partial z^2}=O(1).  \   \nonumber       
  \end{eqnarray}     
\end{corollary}
  \begin{lemma} 
        \label{L3.3}
        $T=-2a\ln h+b_3+O(h\ln h)$, where $a=a(z), b_3=b_3(z), a=1/\omega_0,$ and $\omega_0>0$ is the eigenvalue of the saddle point $C$.
                   \end{lemma} 
          \begin{corollary}  
         $T>c_4^{-1}|\ln h|$. 
         \end {corollary}  
                                                                    
   {\bf Remark.} For the period $T_i$ of the unperturbed trajectory $E=h<0$ in the region $G_i, i=1, 2$, the following expansion is valid:
   $$
   T_i=-a\ln |h|+b_i+O(h\ln|h|),
   $$
   where $b_3=b_1+b_2$. 
    \begin{lemma} 
        \label{L3.4}
          \begin{equation}   \label{3.5} 
   \frac{\partial}{\partial h}\frac{1}{T} \oint_{E=h}\psi dt=O(h^{-1}\ln^{-2} h) .
    \end{equation}                                                  
     \end{lemma}

      \subsection{Lemmas on perturbed motion}
  \label{ss3.2}
  Let $C\xi\eta=C\xi\eta(z)$ be the system of principal axes for the saddle point $C$, oriented as it is shown in  Fig. \ref {unperturbed_plane}.
   \begin{lemma} 
        \label{L3.5}
        Let at a moment of time $t'$ the  point $(p(t'), q(t'))$ lie on the axis  $C\eta (z(t'))$ in $c_5^{-1}$-neighbourhood of the point $C$, and $c_6\eps\le h(t')\le c_7^{-1}\le c_2^{-1}, \ z(t') \in B-2c_1^{-1}$. Then there exists  a moment of time $t''=t'+O(\ln h(t'))$ such that
        a) for $t'\le t\le t''$ the solution $(p(t), q(t), z(t))$ is well defined, and the point  $(p(t''), q(t''))$ lies on the axis $C\eta(z(t''))$ in $c_5^{-1}$-neighbourhood of the point $C$; 
       b)  for $t'\le t\le t''$ the following estimates are satisfied: 
          \begin{eqnarray}      \label{3.6}     
       h(t)=h(t')+O(\eps), \ \frac1{2}h(t')\le h(t)\le 2h(t') ,\ \nonumber   \\
         |z(t)-z(t')|+ |j(t)-j(t')|=O(\eps\ln h(t')); \
              \end{eqnarray} 
              c) integrals  of  functions $\psi, \varphi$  (see Section \ref{3.1})  along the trajectory are estimated as follows:
                 \begin{eqnarray}      \label{3.7}     
        \int_{t'}^{t''}\psi dt =\oint_{E=h(t')}\psi dt+\eps O(h^{-1}(t')), \ \nonumber   \\
          \int_{t'}^{t''}\varphi dt= \oint_{E=h(t')}\varphi dt+\eps O(h^{-1/2}(t')). \
              \end{eqnarray} 
              Here integrals in the right hand side are calculated along the unperturbed trajectory for $z=z(t'), \,\eps=0$.
  \end{lemma} 
         \begin {corollary} 
         \begin {eqnarray}   \label{3.8}                                                              
   t''-t'=  \int_{t'}^{t''}dt=T(h(t'), z(t'))+\eps O(h^{-1}(t')),\phantom{*******}  \ \\
    h(t')-h(t'')=-\eps \int_{t'}^{t''}\left(\frac{\partial E}{\partial q}f_1+\frac{\partial E}{\partial p}f_2+\frac{\partial E}{\partial z}f_3\right)dt \phantom{*******}  \ \nonumber \\
    =-\eps\oint_{E=h(t')}\left(\frac{\partial E}{\partial q}f_1^0+\frac{\partial E}{\partial p}f_2^0+\frac{\partial E}{\partial z}f_3^0\right)dt +\eps^2 O(h^{-1/2}(t')) >\frac1{2}c_3^{-1}\eps. \nonumber  \          
  \end{eqnarray}     
\end{corollary}
  \begin{lemma} 
        \label{L3.6}
       Let at a moment of time $t'$ the conditions $c_6\eps\le h(t')\le c_7^{-1}, \ z(t')\in B-2c_1^{-1}$ be satisfied. Then there exists a moment of time $t''=t'+O(\ln h(t'))$  such that for $t'\le t\le t''$ the solution $(p(t), q(t), z(t))$ is well defined, it meets  estimates (\ref{3.6}), and the point  $(p(t'', q(t''))$ lies on the axis $C\eta(z(t''))$ in  $c_5^{-1}$-neighbourhood of the point $C$.
  \end{lemma} 
         
           \subsection{Proof of Proposition \ref{Pr2.1}}
  \label{ss3.3}
 I. Let $t_*$ be the maximal  moment of time  at the interval $[0, K/\eps]$ such that for $0\le t\le t_*$ the following estimates hold:
  $$
  (p(t), q(t),z(t))\in D-c_1^{-1}, \ h(t)\ge \frac1{2}c_7^{-1}.
  $$
 For $0\le t\le t_*$ the frequency of the unperturbed motion is separated from $0$ by a positive constant. Therefore, usual estimates of the averaging method for one-frequency systems \cite{bm} are valid: 
      \begin{equation}   \label{3.9} 
 |h(t)-H(\eps t)|+  |z(t)-Z(\eps t)|+ |j(t)-J(\eps t)|=O(\eps).
    \end{equation} 
    Thus, for  $0\le t\le t^*$ estimates of  Proposition \ref{Pr2.1} hold, and $h(t_*)=\frac1{2}c_7^{-1}$.       
                                             
  The moment $t_*$ meets  conditions of  Lemma \ref{L3.6}. Due to this lemma a moment of time $t_1=t_*+O(1)$ is defined  such that for  $t_*\le t\le t_1$ the estimate (\ref{3.9}) holds, and  the point  $ (p(t_1), q(t_1))$   lies on the axis $C\eta(z(t_1))$ in  $c_5^{-1}$-neighbourhood of the saddle point $C$.
  
  \medskip
II.  Let $t_{**}$ be  the maximal moment of time in the interval  $[t_1, K/\eps]$ such that for $t_1\le t\le t_{**}$  the following estimates hold:
     \begin{equation}   \label{3.10} 
  (p(t), q(t),z(t))\in D-c_1^{-1}, \ h(t)\ge c_6 \eps,  \
  \frac1{2}H(\eps t)\le h(t)\le 2H(\eps t).
   \end{equation} 
  Lemma \ref{3.5} allows to define  moments of time  $t_2, ..., t_{n^*}$ of consecutive arrivals of the point    $ (p(t), q(t))$  at the ray $C\eta$, where $n^*$ is the maximal  number $n$ such that $t_n<t_{**}$. Denote  $h_n=h(t_n)$ and,  analogously, $z_n, j_n, H_n, Z_n, J_n$. From (\ref{3.6}) we have that for   $t_n\le t\le t_{n+1}$ the following estimates are satisfied:
        \begin{eqnarray}      \label{3.11}     
       |h(t)-h_n|=O(\eps), \  \frac1{2} h_n\le h(t)\le 2h_n ,\ \nonumber   \\
         |z(t)-z_n|+ |j(t)-j_n|=O(\eps\ln h_n). \
              \end{eqnarray} 
 Lemma \ref{L3.2} allows to estimate the right hand sides of the averaged system:
 $$
 \dot H=O(\eps/\ln H), \ \dot Z=O(\eps), \ \dot J=O(\eps).
 $$        
Considering motion in the averaged system for   $t_n\le t\le t_{n+1}$ and making use of estimates (\ref{3.6}), (\ref{3.10}),  we get for   $t_n\le t\le t_{n+1}$    
    \begin{eqnarray}      \label{3.12}     
       |H(\eps t)&-&H_n|=O(\eps), \  \nonumber   \\
         |Z(\eps t)&-&Z_n|+ |J(\eps t)-J_n|=O(\eps\ln h_n). \
              \end{eqnarray} 
If $n=n^*$, then (\ref{3.11}), (\ref{3.12}) hold  for  $t_{n^*}\le t\le t_{**}$.

From formulas
   \begin {eqnarray}                                                            
     j_{n+1}-j_n&=&\eps \int_{t_n}^{t_{n+1}}\left(\frac{\partial I}{\partial h}\left(\frac{\partial E}{\partial q}f_1+\frac{\partial E}{\partial p}f_2+\frac{\partial E}{\partial z}f_3\right)+\frac{\partial I}{\partial z}f_3\right)dt , \  \nonumber  \\ 
      z_{n+1}-z_n&=& \eps \int_{t_n}^{t_{n+1}}f_3\,dt \, , \nonumber 
    \end {eqnarray}     
making use of (\ref{3.3}), (\ref{3.4}), (\ref{3.7}), (\ref{3.11}), we get
 \begin{eqnarray}        
      j_{n+1}=j_n+\eps T(h_n, z_n)F(h_n, z_n)+\eps^2O(h_n^{-1} ),\  \nonumber \\
      z_{n+1}=z_n+\eps T(h_n, z_n)\Phi(h_n, z_n)+\eps^2O(h_n^{-1}), \  \nonumber  \
                     \end{eqnarray} 
where $\eps F(h, z)$ and $\eps \Phi(h, z)$ are right hand sides of the averaged equations for $ I$ and  $ z$ respectively:
 \begin {eqnarray}                                                            
     2\pi F(h, z)&=&\oint_{E=h}\left(\frac{\partial E}{\partial q}f_1^0+\frac{\partial E}{\partial p}f_2^0+\frac{\partial E}{\partial z}f_3^0\right)\!dt-\frac{1}{T}\oint_{E=h}\frac{\partial E}{\partial z}\,dt\cdot\oint_{E=h}f_3^0\, dt, \  \nonumber  \\ 
      \Phi(h, z)&=&\frac{1}{T}\oint_{E=h}f_3^0\,dt. \ \nonumber 
    \end {eqnarray}     
In the averaged system
$$
\dot J=\eps F(H, Z), \ \dot Z=\eps \Phi(H, Z).
$$
Therefore, making use of (\ref{3.3}), (\ref{3.4}), (\ref{3.6}), (\ref{3.12}),  we get
 \begin{eqnarray}        
      J_{n+1}-J_n&=&\eps \int_{t_{n}}^{t_{n+1}}F(H, Z)dt=\eps \int_{t_{n}}^{t_{n+1}}F(H_n, Z_n)dt+\eps^2O\left(\frac{h_n^{-1}}{\ln h_n}\right) \  \nonumber \\
      &=&\eps T(h_n, z_n)F(H_n, Z_n)+\eps^2O(h_n^{-1} ),\  \nonumber  \\
       Z_{n+1}-Z_n&=&\eps \int_{t_{n}}^{t_{n+1}}\Phi(H, Z)dt=\eps \int_{t_{n}}^{t_{n+1}}\Phi(H_n, Z_n)dt+\eps^2O\left(\frac{h_n^{-1}}{\ln h_n}\right)\  \nonumber \\
      &=&\eps T(h_n, z_n)\Phi(H_n, Z_n)+\eps^2O(h_n^{-1}). \  \nonumber  \
  \end{eqnarray} 
  So
 \begin{eqnarray}        
      J_{n+1}=J_n+\eps T(h_n, z_n)F(H_n, Z_n)+\eps^2O(h_n^{-1}), \  \nonumber \\
     Z_{n+1}=Z_n+\eps T(h_n, z_n)\Phi(H_n, Z_n)+\eps^2O(h_n^{-1}). \  \nonumber 
       \end{eqnarray} 
       Denote
       $$
       u_n=|j_n-J_n|+|z_n-Z_n|.
$$
From the previous estimates by means of  (\ref{3.3}), (\ref{3.5}), we get
  \begin{eqnarray}      \label{3.13}     
      u_{n+1}\le u_n+ |h_n-H_n|\eps O(h_n^{-1}\ln^{-1} h_n)\\
      +|z_n-Z_n|\eps O(\ln h_n) +\eps^2O(h_n^{-1}). \nonumber\  
      \end{eqnarray} 
       We have  
      \begin{eqnarray}
           2\pi\partial I/\partial h=T>c_4^{-1}|\ln h|, \ \partial I/\partial z=O(1). \ \nonumber 
 \end{eqnarray}
 Therefore, solving relation $I=I(h,z)$ for $h$ we get $h$ as a function of $I,z$ for which
 $$
\frac{\partial h}{\partial I}=\frac{2\pi}{T}=O(\ln^{-1} h), \ \frac{\partial h}{\partial z}=-\frac{\partial I/\partial z}{\partial I/\partial h} =O(\ln^{-1} h).
$$
Therefore, by Lagrange's formula, 
 \begin{eqnarray}        
      h_n-H_n=\left(\frac{\partial h}{\partial I}\right)_*(j_n-J_n)+ \left(\frac{\partial h}{\partial z}\right)_*(z_n-Z_n)=\  \nonumber \\
    = \left( |j_n-J_n|+|z_n-Z_n|\right) O(\ln^{-1}h_n). \  \nonumber 
       \end{eqnarray} 
Symbol  `` * '' here means that derivatives are calculated at some point in the straight line interval    with  endpoints  $(j_n, z_n)$ and $(J_n, Z_n)$. Using  this estimate in   (\ref{3.13}), we get
\begin{equation}   \label{3.14} 
u_{n+1}\le \left[1+\eps O(h_n^{-1}\ln^{-2} h_n)\right]u_n+\eps^2O(h_n^{-1}).
\end{equation} 
Consecutive use of this relation gives
  \begin{equation*}    
  u_{n+1}\le \left[\prod\limits^n_{s=1}(1+\eps |O(h_s^{-1}\ln^{-2} h_s)|\right]\left(u_1+\eps^2 \sum\limits_{s=1}^n|O( h_s^{-1})|\right). 
       \end{equation*} 
In accordance with  (\ref{3.8}), $h_s-h_{s+1}>\frac1{2}c_3^{-1}\eps$.
Therefore
 \begin{eqnarray} 
 \sum\limits_{s=1}^n\eps^2|O( h_s^{-1})|\le\eps^2 c_8\sum\limits_{s=1}^n h_s^{-1}\le2c_8c_3\eps \sum\limits_{s=1}^n h_s^{-1}(h_s-h_{s+1})\ \nonumber \\
 <\eps c_9 \int\limits_{h_{n+1}}^{h_1}h^{-1}dh=\eps c_9\ln(h_1/h_{n+1})=\eps O(\ln h_{n+1}).  \nonumber \
   \end{eqnarray} 
In an analogous way, from the convergence of $ \int\limits_0^{1/2}h^{-1}\ln^{-2}h\,dh$,  we get
 \begin{eqnarray} 
 \prod\limits_{s=1}^n\left(1+\eps|O( h_s^{-1}\ln^{-2}h_s)|\right)=\exp\left[\sum\limits_{s=1}^n\ln(1+\eps| O( h_s^{-1}\ln^{-2}h_s)|\right]   \ \nonumber \\
\le\exp(\eps\sum\limits_{s=1}^nO( h_s^{-1}\ln^{-2}h_s)) =O(1) .\nonumber \
   \end{eqnarray} 
Using  these estimates in   (\ref{3.14}) and taking into account that $u_1=O(\eps)$, we get for $n=1, ...,n^*-1$:
  \begin{equation}  
  u_{n+1}=|j_{n+1}-J_{n+1}|+|z_{n+1}-Z_{n+1}|=O( \eps\ln h_{n+1} ).
   \ \nonumber \
\end{equation}
From here we get with the help of   (\ref{3.10}),  (\ref{3.11}),  (\ref{3.12}) for $t_1\le t\le t_{**}$:
 \begin{eqnarray}      \label{3.15} 
|j(t)-J(\eps t)|&=&O(\eps\ln H(\eps t)), \\
|z(t)-Z(\eps t)|&=&O(\eps\ln H(\eps t)) , \ |h(t)-H(\eps t)|<c_{10}\eps. \ \nonumber \
         \end{eqnarray}
         Let us choose $k_1=2c_{10}+c_6$. While $H(\eps t)\ge k_1\eps$, the condition  (\ref{3.10}) can not be violated. Hence, the estimates in  Proposition \ref{Pr2.1} for $h(t)$ are proved.
         
         \medskip
          III. The estimate for $z(t)$ in  (\ref{3.15}) is less accurate, than in formulation of  Proposition \ref{Pr2.1}. Now we will   improve this estimate.
         Denote
 \begin {eqnarray} 
 (\Delta z)_n&=&z_n-Z_n, \   (\Delta h)_n=h_n-H_n, \    \nonumber  \\                                                       
     a_n&=&-\oint_{E=h_n}\left(\frac{\partial E}{\partial q}f_1^0+\frac{\partial E}{\partial p}f_2^0+\frac{\partial E}{\partial z}f_3^0\right)dt, \ b_n=\oint_{E=h_n}(f_3^0-f_{3C}^0)dt. \ \nonumber 
    \end {eqnarray}     
    Here $f_{3C}^0=f_{3C}(z)$ is the value of the function $f_3^0$ at the saddle point $C$. The integrals  are calculated for $z=z_n$.
    We have an identity
 \begin {eqnarray}   \label{3.16} 
  \dot z&-&\dot Z=\eps[(f_{3C}^0(z)-f_{3C}^0(Z))+ (f_{3}^0(p, q, z)-f_{3C}^0(z)) \\
  &-&\frac{1}{T(H, Z)}\oint_{E=H}(f_3^0(p, q, Z)-f_{3C}^0(Z))dt+(f_3(p, q, z, \eps)-(f_3^0(p, q, z))] . \   \nonumber   
    \end {eqnarray}  
    Let us integrate both sides of this identity with respect to $t$ from $t_n$ to $t_{n+1}$ and estimate the right hand side making use of already established estimates   (\ref{3.15}) and Lemmas \ref{L3.2}, \ref{L3.3}, \ref{L3.5}. In the left hand side we have  $ (\Delta z)_{n+1}-(\Delta z)_n$. Terms in the right hand side are
 \begin {eqnarray}  
  \int_{t_n}^{t_{n+1}}(f_{3C}^0(z)-f_{3C}^0(Z))dt&=&O(\eps\ln^2h_n),  \   \nonumber \\
   \int_{t_n}^{t_{n+1}}(f_3^0(p, q, z)-f_{3C}^0(z))dt&=&b_n+ O(\eps h_n^{-1/2}), \   \nonumber \\  
     \int_{t_n}^{t_{n+1}} \frac{dt}{T(H, Z)}\left(\oint_{E=H}(f_3^0(p, q, Z)-f_{3C}^0(Z))dt\right)&=&b_n \int_{t_n}^{t_{n+1}} \frac{dt}{T(H, Z)} + O(\eps\ln^2h_n),  \nonumber \\
     \int_{t_n}^{t_{n+1}}(f_3^0(p, q, z, \eps)-f_{3}^0(p, q, z))dt&=&O(\eps\ln h_n). \  \nonumber 
    \end {eqnarray} 
    Finally, we have 
 \begin {equation}   \label{3.17} 
 (\Delta z)_{n+1}-(\Delta z)_n=\eps b_n\left(1- \int_{t_n}^{t_{n+1}} \frac{dt}{T(H, Z)}\right) +\eps^2O(h_n^{-1/2}). \
 \end {equation}
 In an analogous way we get
 \begin {equation}   \label{3.18} 
 (\Delta h)_{n+1}-(\Delta h)_n=-\eps a_n\left(1- \int_{t_n}^{t_{n+1}} \frac{dt}{T(H, Z)}\right) +\eps^2O(h_n^{-1/2}) .\
 \end {equation}
In accordance with  (\ref{3.2}), $a_n>c_3^{-1}$. Let us denote $\mu_n=a_n^{-1}b_n$. It follows from  (\ref{3.17}),  (\ref{3.18}) that 
$$
(\Delta z)_{n+1}-(\Delta z)_n=-[(\Delta h)_{n+1}-(\Delta h_n)]\mu_n+\eps^2O(h_n^{-1/2}).
$$
Consecutive  use of this relation gives
 \begin {eqnarray} 
(\Delta z)_{n+1}=(\Delta z)_1-\sum\limits_{s=1}^n[(\Delta h)_{s+1}-(\Delta h)_s)]\mu_s+\eps^2\sum\limits_{s=1}^nO(h_s^{-1/2}). \   \nonumber 
\end {eqnarray}
Moreover, 
\begin {eqnarray} 
\sum\limits_{s=1}^n [(\Delta h)_{s+1}-(\Delta h)_s]\mu_s=(\Delta h)_{n+1}\mu_n-(\Delta h)_1\mu_1-\sum\limits_{s=2}^n(\Delta h)_s(\mu_s-\mu_{s-1}). \nonumber \
   \end {eqnarray} 
   
The definition of $\mu_s$, Lemma \ref{L3.2} and estimates  (\ref{3.11}) imply that  \\ $\mu_s-\mu_{s-1}=O(\eps\ln h_s)$. As 
$(\Delta h)_{n+1}=O(\eps), \  (\Delta z)_1=O(\eps), \ \mu_n=O(1)$, so
$$
 z_{n+1}-Z_{n+1}\equiv (\Delta z)_{n+1}=O(\eps)+\eps^2\sum\limits_{s=1}^nO(h_s^{-1/2})=O(\eps).
 $$
Integrating  both sides of   (\ref{3.16}) with respect to time from $t_n$ to some $t \in (t_n, t_{n+1})$ and estimating the right hand side,  we get for $ t_n \le t\le t_{n+1}$
$$
z(t)-Z(\eps t)=z_n-Z_n+O(\eps)
$$
(in the case when  $n=n^*$ this estimate is valid for  $ t_{n^*} \le t< t_{**}$). Thus, for $0\le t\le t_{**}$ we have $z(t)-Z(\eps t)=O(\eps)$. Hence, this estimate is valid while $H(\eps t)\ge k_1\eps$ (as the last inequality certainly holds for  $0\le t\le t_{**}$). The Proposition \ref{Pr2.1} is proved.

  \medskip
  {\bf{Remark}}. In the proof of the last estimate  representation  (\ref{3.16}) was used.  We can not use this representation in  problems  where separatrices connect different saddle points, and at these points the function $f_3^0(p, q, z)$ has different values. The accuracy of the description of behaviour of $z$ in these cases  is worse  than $O(\eps)$ in general. This accuracy is given by the obtained above estimate   (\ref{3.15}): $z(t)-Z(\eps t)=O(\eps\ln H(\eps t))$.
 \medskip

IV. The proof of  Proposition \ref{Pr2.2} is completely analogous to the given above proof of the Proposition \ref{Pr2.1}. But final estimates are different. The reason is as  follows.
 As above,  in nn. I, II, we prove that  
 $$
 |j(t)-J(\eps t)|+|z(t)-Z(\eps t)|=O(\eps)+O\left(\eps\ln\frac{H(\eps t)}{H(\eps t')}\right)
$$
under  assumptions  of Proposition \ref{Pr2.2}
(we will omit subscript $\nu$ and superscript  $``  \ '    "$ in the notation $(H_\nu', Z_\nu'))$.
Here $t'$ is the initial moment of time for the motion in    Proposition \ref{Pr2.2}. 
Because  $H(\eps t')\sim\eps$,  we have
 $$
 |j(t)-J(\eps t)|+|z(t)-Z(\eps t)|=O(\eps\ln\eps)
$$
(for  Proposition \ref {Pr2.1} there  was $H(\eps t')\sim 1$, and there  was used the estimate $O(\eps\ln H(\eps t))$ in the right hand side of this equality).
As in n. II, to estimate $|h(t)-H(\eps t)|$  from here, we use Lagrange's formula
$$
h(t)-H(\eps t)= \left(\frac{\partial h}{\partial I}\right)_*(j(t)-J(\eps t))+\left(\frac{\partial h}{\partial z}\right)_*(z(t)-Z(\eps t)).
$$
The symbol  ``\, *\, " here indicates  that derivatives are calculated at some point in the straight line interval     with  endpoints  $(j(t), z(t))$ and $ (J(\eps t), Z(\eps t))$. As by assumption   $\frac{1}{2}|H(\eps t)|\le|h(t)|\le2|H(\eps t)|$,  we have
\begin {eqnarray} 
 \left(\frac{\partial h}{\partial I}\right)_*&=& \frac{2\pi}{T_*}= O\left(\frac{1}{1+|\ln|H(\eps t)||}\right),  \nonumber \\
   \left(\frac{\partial h}{\partial z}\right)_*&=&-\left( \frac{\partial I/\partial z}{\partial I/\partial h}\right)_*= O\left(\frac{1}{1+|\ln|H(\eps t)||}\right). \ \nonumber \
\end {eqnarray}

Therefore
$$
h(t)-H(\eps t)= O\left(\frac{\eps\ln\eps}{1+|\ln|H(\eps t)||}\right),
$$
and this is the assertion  of   Proposition \ref{Pr2.2} concerning  an accuracy of description of $h$.

Now, the estimate of $z$ can be  improved as in n. III, and we get
$$
z(t)-Z(\eps t)= O\left(\frac{\eps\ln\eps}{1+|\ln|H(\eps t)||}\right).
$$
This is the assertion of   Proposition \ref{Pr2.2} concerning  the accuracy of description of $z$.
 \medskip
 
   \subsection{Proofs of  Lemmas on unperturbed motion}
  \label{ss3.4}
    \subsubsection{Preliminary estimates}
  \label{sss3.4.1}

  Let  $\eta, \xi$ be  such variables  that the quadratic part of the Hamiltonian near  the saddle  $C$ is $\frac12\omega_0(z)(\eta^2-\xi^2), \ \omega_0>0$, and the transformation $(p,q)\mapsto  (\eta, \xi)$ is a canonical transformation containing $z$ as a parameter.   Let us denote ${\cal {H}}(\eta, \xi, z)$ the Hamiltonian $E$ expressed via $\eta, \xi, z$.
 \begin{lemma} 
        \label{L3.7}
        Let $|\xi|<d_1^{-1}, \  |\eta|<d_1^{-1}$.
        
       \noindent  If ${\cal {H}}(\eta, \xi, z)>0$, then
       \begin {eqnarray}   \label{3.19}                                                              
 \eta&=&\pm\sqrt{2\omega_0^{-1}{\cal {H}}+\xi^2} \ +O({\cal {H}}+\xi^2), \ \nonumber \\
   \partial {\cal {H}}/\partial \eta&=&\omega_0\eta+O(\eta^2), \ \nonumber \\
   d_2^{-1}\sqrt{{\cal {H}}+\xi^2}&<&| \partial {\cal {H}} / \partial \eta |< d_2\sqrt{{\cal {H}}+\xi^2}. \
     \end{eqnarray} 
      If ${\cal {H}}(\eta, \xi, z)<0$, then
       \begin {eqnarray}                                                             
 \xi&=&\pm\sqrt{2\omega_0^{-1}|{\cal {H}}|+\eta^2} \ +O(|{\cal {H}}|+\eta^2), \ \nonumber \\
   \partial {\cal {H}}/\partial \xi&=&-\omega_0\xi+O(\xi^2), \ \nonumber \\
   d_2^{-1}\sqrt{|{\cal {H}}|+\eta^2}&<&| \partial {\cal {H}} / \partial \xi|< d_2\sqrt{|{\cal {H}}|+\eta^2}.  \ \nonumber \
     \end{eqnarray} 
 \end{lemma} 
 {\it {Proof.}}
 
 Consider the case ${\cal {H}}>0$. The  case ${\cal {H}}<0$ is analogous. We have
  \begin {eqnarray} 
     {\cal {H}}&=&\frac12\omega_0(\eta^2-\xi^2)+O(\eta^3)+O(\eta^2\xi)+O(\eta\xi^2)+O(\xi^3),   \ \nonumber \\
     \eta^2(1+O(\eta)+O(\xi))&=&2\omega_0^{-1}{\cal {H}}+\xi^2+\xi^2(O(\eta)+O(\xi)), \ \nonumber \\
     \eta^2&=&\frac{2\omega_0^{-1} {\cal {H}}+\xi^2}{1+O(\eta)+O(\xi)}+ \xi^2(O(\eta)+O(\xi)). \ \nonumber \
      \end{eqnarray} 
      From this equality we derive that
       \begin {eqnarray} 
 \nu_1^{-1}|\xi|< \nu_1^{-1}\sqrt{{\cal {H}}+\xi^2}<|\eta|< \nu_1\sqrt{{\cal {H}}+\xi^2},  \ \nonumber \\
     \eta^2=2\omega_0^{-1}{\cal {H}}+\xi^2+O((\sqrt{{\cal {H}}+\xi^2})^3), \ \nonumber \\
     \eta=\pm\sqrt{2\omega_0^{-1}{\cal {H}}+\xi^2}+O({\cal {H}}+\xi^2). \ \nonumber \\
\partial{\cal {H}}/\eta=\omega_0\eta+O(|\eta|^2+|\xi|^2)=\omega_0\eta+O(\eta^2), \ \nonumber \\
   d_2^{-1}\sqrt{{\cal {H}}+\xi^2}<\left| \frac{\partial {\cal {H}}}{\partial \eta}\right|< d_2\sqrt{{\cal {H}}+\xi^2}.   \ \nonumber \
      \end{eqnarray} 
Lemma \ref{L3.7} is proved.
    \begin {corollary}
    \label{C3.7}
  a)  For   $|\xi|<d_3^{-1}<d_1^{-1}, \ 0\le h\le d_4^{-1}$ the equation ${\cal {H}}(\eta, \xi, z)=h$ defines a unique
  $\eta=\tilde\eta(h, \xi, z)$ such that $0\le\tilde\eta(h, \xi, z)<\frac12 d_2^{-1}$. For $h+\xi^2>0$ the function $\tilde\eta$ is smooth and
    \begin {equation}   \label{3.20}                                                              
     \frac{\partial \tilde\eta}{\partial h}=\frac{1}{\partial {\cal {H}} /\partial \eta}, \   \frac{\partial \tilde\eta}{\partial z}=-\frac{\partial {\cal {H}}/ \partial z}{\partial {\cal {H}} / \partial \eta} \,.
          \end{equation} 
        The same  equation ${\cal {H}}(\eta, \xi, z)=h$ defines also a unique 
  $\eta=\tilde{\tilde\eta}(h, \xi, z)$ such that $-\frac12d_2^{-1}<\tilde{\tilde\eta}\le0$ with analogous properties.
  
  b)  For $|\eta|<d_3^{-1}<d_1^{-1}, \ -d_4^{-1}\le h\le0$ the equation ${\cal {H}}(\eta, \xi, z)=h$ defines a unique
  $\xi=\tilde\xi(h, \eta, z)$ such that $0\le\tilde\xi(h, \eta, z)<1/2 d_2$. For $|h|+\eta^2>0$ the function $\tilde\xi$ is smooth and
    \begin {equation}                                                                
     \frac{\partial \tilde\xi}{\partial h}=\frac{1}{\partial {\cal {H}} / \partial \xi}, \    \frac{\partial \tilde\xi}{\partial z}=-\frac{\partial {\cal {H}} / \partial z}{\partial {\cal {H}} / \partial \xi} \, .   \nonumber \
          \end{equation} 
   The same  equation ${\cal {H}}(\eta, \xi, z)=h$ defines also a unique 
  $\xi=\tilde{\tilde\xi}(h, \xi, z)$ such that $-\frac12 d_2^{-1}<\tilde{\tilde\xi}\le0$ with analogous properties.
\end{corollary}
 \medskip
 
 \subsubsection{Proof of  Lemma \ref{L3.1}}
  \label{sss3.4.2}
  To save notation we consider the case when  $\vfi$ does not depend on  $\eps$: $\vfi=\vfi(p,q,z)$. 
   Let $0<h<\frac12d_4^{-1}$. Denote
  $$
  R(h)=\oint\limits_{E=h}\varphi dt, \ R(0)=\oint_{l_3}\varphi dt\, .
  $$
  
  $1^0$. Split  $R(h)$ into integrals $R'(h)$ and $R''(h)$  over the segments of the curve $E=h$, situated  inside  and outside of the  rectangle $|\xi|<d_3^{-1}, \  |\eta|<d_1^{-1}$ respectively. Split  $R'(h)$ into integrals $R'_i(h), i=1, 2, 3, 4$, where $R'_i(h)$ is calculated  over the segment of the curve $E=h$, situated into the $i$th quadrant of the coordinate system $C\xi\eta$. Split $R(0)$ in the same way.
  
  $2^0$.  The integral $R''(h)$ is calculated over the arcs separated from the singularity (i.e. from the point $C$). It is easy to check that
  $$
  R''(h)=R''(0)+O(h).
  $$
  
$3^0$. Estimate  $R'_1(h)$. Let us use  $\xi$ as an  independent variable in this integral. According to  Corollary \ref {C3.7} we get

 $$
  R_1'(h)=\int_0^{d_3^{-1}}\left(\frac{\varphi}{\dot\xi}\right)_hd\xi=\int_0^{d_3^{-1}}\left(\frac{\varphi}{\partial {\cal {H}} / \partial \eta}\right)_hd\xi,
  $$
where the subscript  `` $h$ '' indicates that we should plug $\eta=\tilde\eta(h, \xi, z)$ into the integrand.
The function $\varphi$ has the form
$$
\varphi=\alpha\xi+\beta\eta+\varphi_2(\eta, \xi, z), \ \varphi_2=O(\xi^2+\eta^2), \ \alpha=\alpha(z), \ \beta=\beta(z).
$$
Represent
  \begin {eqnarray}  \label{3.21} 
  R_1'(h)&-&R'_1(0)
  =\alpha\int\limits_0^{d_3^{-1}}\left[\frac{\xi}{(\partial {\cal {H}} / \partial \eta)_h}-\frac{\xi}{(\partial {\cal {H}} / \partial \eta)_0}\, \right]d\xi \\
  &+&\beta\int\limits_{\sqrt h}^{d_3^{-1}}\left[\left(\frac{\eta}{\partial {\cal {H}} / \partial \eta}\right)_h-\left(\frac{\eta}{\partial {\cal {H}} / \partial \eta}\right)_0\, \right]d\xi  \  
  + \beta\int\limits_0^{\sqrt h}\left[\left(\frac{\eta}{\partial {\cal {H}} / \partial \eta}\right)_h-\left(\frac{\eta}{\partial {\cal {H}} / \partial \eta)}\right)_0\, \right]d\xi  \   \nonumber \\ 
 &+&\int\limits_{\sqrt h}^{d_3^{-1}}\left[\left(\frac{\varphi_2}{(\partial {\cal {H}} / \partial \eta}\right)_h-\left(\frac{\varphi_2}{(\partial {\cal {H}} / \partial \eta}\right)_0\, \right]d\xi+ \int\limits_0^{\sqrt h}\left[\left(\frac{\varphi_2}{(\partial {\cal {H}} / \partial \eta}\right)_h-\left(\frac{\varphi_2}{\partial {\cal {H}} / \partial \eta}\right)_0\, \right]d\xi . \   \nonumber
  \end{eqnarray} 
  
  Estimate  the last integral in (\ref{3.21}). According to  Lemma \ref{L3.7}, the integrand in this integral is $O(\sqrt{h})$. Therefore,   this integral is $O(h)$.  
  
  Estimate the next to last  integral in (\ref{3.21}). Using Lagrange's formula and Corollary  \ref{C3.7}, we get, that this integral is equal to 
   \begin {equation}
    \int\limits_{\sqrt h}^{d_3^{-1}}h\left(\frac{\partial}{\partial h}\left(\frac{\varphi_2}{\partial {\cal {H}} / \partial \eta}\right)_h\right)_{h_*}d\xi=h\int\limits_{\sqrt h}^{d_3^{-1}}\left(\frac{\partial}{\partial \eta}\left(\frac{\varphi_2}{\partial {\cal {H}} / \partial \eta}\right)\frac{1}{\partial{\cal {H}}/ \partial \eta}\right)_{h_*}d\xi,  \nonumber \
  \end{equation} 
  where $0<h_*<h$.   By Lemma \ref{L3.7}, the integrand  is $O((h_*+\xi^2)^{-1/2})$. Therefore,   this integral  is
  \begin {equation}
   O(h) \int\limits_{\sqrt h}^{d_3^{-1}}\frac{d\xi}{\sqrt{h_*+\xi^2}}d\xi=O(h)\int\limits_{\sqrt h}^{d_3^{-1}}\frac{d\xi}{\xi }=O(h\ln h).  \nonumber \
  \end{equation} 
  Consider the third integral in (\ref{3.21}). It is equal to
 \begin {equation}
   \frac{1}{\omega_0} \int\limits_{0}^{\sqrt h}\left[\left(\frac{\eta}{\eta+\chi_2}\right)_h-\left(\frac{\eta}{\eta+\chi_2}\right)_0\right]d\xi= \frac{1}{\omega_0} \int\limits_{0}^{\sqrt h}\frac{(\eta)_h(\chi_2)_0-(\eta)_0(\chi_2)_h}{(\eta+\chi_2)_h(\eta+\chi_2)_0}d\xi . \   \nonumber
    \end{equation} 
Here $\chi_2$ is a function of $\eta, \xi, z$; expansion of $\chi_2$ with respect to  $\eta, \xi$ starts with  second order terms.  By  Lemma \ref{L3.7}, the integrand   is $O(\sqrt{h})$.  The integral is, therefore, $O(h)$.

 Estimate the second integral in (\ref{3.21}). Using Lagrange's formula and Corollary  \ref{C3.7}, we get that this integral is equal to 
 \begin {equation}
   \frac{h}{\omega_0^2} \int\limits_{\sqrt h}^{d_3^{-1}}\left(\frac{\partial}{\partial\eta}\left(\frac{\eta}{\eta+\chi_2}\right)\cdot\frac{1}{\eta+\chi_2}\right)_{h_*}d\xi= \frac{h}{\omega_0^2} \int\limits_{\sqrt h}^{d_3^{-1}}\left(\frac{\chi_2-\eta\,\partial{\chi_2/\partial\eta}}{(\eta+\chi_2)^3}\right)_{h_*}d\xi .  \   \nonumber
    \end{equation} 
The integrand  here is $O(1/\xi)$.  Therefore the second integral  in  (\ref{3.21}) is $O(h\ln h)$.

 Estimate the first integral in (\ref{3.21}). It is equal to 
 \begin{eqnarray} 
   \frac{1}{\omega_0} \int\limits_{0}^{d_3^{-1}}\left[\frac{\xi}{\sqrt{2\omega_0^{-1}{\cal {H}}+\xi^2}}-\frac{\xi}{|\xi|}\right]d\xi \phantom{******************} \   \nonumber \\   
+ \frac{1}{\omega_0} \int\limits_{0}^{d_3^{-1}}\xi\left[\left(\frac{\omega_0\sqrt{2\omega_0^{-1}{\cal {H}}+\xi^2}-\partial{\cal {H}}/\partial\eta}{\sqrt{2\omega_0^{-1}{\cal {H}}+\xi^2}\,\partial{\cal {H}}/\partial\eta}\right)_h  -\left(\frac{\omega_0|\xi|-\partial{\cal {H}}/\partial\eta}{|\xi|\partial{\cal {H}}/\partial\eta}\right)_0\right]d\xi . \   \nonumber
 \end{eqnarray} 
We can estimate the second integral in this expression by splitting  it into two integrals, from $0$ to $\sqrt{h}$ and from $\sqrt{h}$  to $d_3^{-1}$. This gives that this integral is $O(h\ln h)$.
So we get
 \begin {equation}
R'_1(h)-R'_1(0)=\frac{\alpha}{\omega_0} \int\limits_{0}^{d_3^{-1}}\left(\frac{\xi}{\sqrt{2\omega_0^{-1}h+\xi^2}}-\frac{\xi}{|\xi|}\right)d\xi+O(h\ln h). \   \nonumber
 \end{equation}
 Analogous estimates are valid for $R'_2,..., R'_4$. In particular, 
 \begin {equation}
R'_4(h)-R'_4(0)=\frac{\alpha}{\omega_0} \int\limits_{-d_3^{-1}}^{0}\left(\frac{\xi}{\sqrt{2\omega_0^{-1}h+\xi^2}}-\frac{\xi}{|\xi|}\right)d\xi+O(h\ln h). \   \nonumber
 \end{equation}
 Therefore
 $$
R'_1(h)+R'_4(h)-R'_1(0)-R'_4(0)=O(h\ln h) 
$$
(we use here that the integral of an odd function over a symmetric  with respect to  $0$ interval is equal to $0$).
The same estimate holds for $R'_2+R'_3$. Taking into account the estimate for $R''$, we  finally get
 $$
R(h)=R(0)+O(h\ln h) .
$$
Lemma \ref{L3.1} is proved.
 \medskip
 
 \subsubsection{Proof of Lemma \ref{L3.2}}
  \label{sss3.4.3}
To save notation we consider the case when  the functions $\vfi, \psi$ do not depend on  $\eps$: $\vfi=\vfi(p,q,z), \psi=\psi(p,q,z)$. The proof uses the same scheme, as the proof of  the Lemma \ref{L3.1} above. The part of the curve $E=h$, situated outside of the neighbourhood of the point $C$, gives a contribution $O(1)$ to any integral in  (\ref{3.3}). In a neighbourhood of the saddle point $C$ we split  the curve  $E=h$ into four segments situated in four  quadrants of the coordinate system $C\xi\eta$. The contribution to the integral of each of these segments is estimates in the same way. For certainty,  let us consider the segment situated in the first quadrant. The corresponding  integral we denote  $R^{(i)}$, where $i$ is the number of the integral in the Lemma \ref{L3.2}, $i\ne 5$. To estimate integrals  $R^{(i)}$ we introduce $\xi$ as an independent variable in the integral and use estimates in Lemma \ref{L3.7} and  Corollary \ref{C3.7}. We have
 \begin{eqnarray*} 
  R^{(1)}&= &\int\limits_{0}^{d_3^{-1}}\left(\frac{|\vfi|}{\partial{\cal {H}}/\partial\eta}\right)_hd\xi= \int\limits_{0}^{d_3^{-1}}\left(\frac{O(\eta)}{|\eta|}\right)_hd\xi= \int\limits_{0}^{d_3^{-1}}O(1)d\xi=O(1). \   \nonumber \\
  R^{(2)}&=& \int\limits_{0}^{d_3^{-1}}\left(\frac{|\psi|}{\partial{\cal {H}}/\partial\eta}\right)_hd\xi= \int\limits_{0}^{d_3^{-1}}\frac{O(1)}{\sqrt{h+\xi^2}}d\xi=O(\ln h). \   \nonumber \\
   R^{(3)}&=&\frac{\partial}{\partial h} \int\limits_{0}^{d_3^{-1}}\left(\frac{\psi}{\partial{\cal {H}}/\partial\eta}\right)_hd\xi
   = \int\limits_{0}^{d_3^{-1}}\left(\frac{\partial}{\partial\eta}\left(\frac{\psi}{\partial{\cal {H}}/\partial\eta}\right)\frac{1}{\partial{\cal {H}}/\partial\eta}\right)_hd\xi  \\
   &=& \int\limits_{0}^{d_3^{-1}}\frac{O(1)d\xi}{(h+\xi^2)^{3/2}}=O\left(\frac{1}{h}\right) .\   \nonumber \\
      R^{(4)}&=&\frac{\partial}{\partial z} \int\limits_{0}^{d_3^{-1}}\left(\frac{\psi}{\partial{\cal {H}}/\partial\eta}\right)_h
      d\xi\\
      &=& \int\limits_{0}^{d_3^{-1}}\left(\frac{\partial}{\partial z}\left(\frac{\psi}{\partial{\cal {H}}/\partial\eta}\right)-\frac{\partial}{\partial \eta}\left(\frac{\psi}{\partial{\cal {H}}/\partial\eta}\right)\frac{\partial{\cal {H}}/\partial z}{\partial{\cal {H}}/\partial\eta}\right)_hd\xi=  \nonumber \\
  &=& \int\limits_{0}^{d_3^{-1}}\frac{O(1)}{\sqrt{h+\xi^2}}d\xi=O(\ln h). \   \nonumber \\
   R^{(6)}&=&\frac{\partial}{\partial z} \int\limits_{0}^{d_3^{-1}}\left(\frac{\varphi}{\partial{\cal {H}}/\partial\eta}\right)_hd\xi\\
   &=& \int\limits_{0}^{d_3^{-1}}\left(\frac{\partial}{\partial z}\frac{\varphi}{\partial{\cal {H}}/\partial\eta}-\frac{\partial}{\partial \eta}\left(\frac{\varphi}{\partial{\cal {H}}/\partial\eta}\right)\frac{\partial{\cal {H}}/\partial z}{\partial{\cal {H}}/\partial\eta}\right)_hd\xi= 
  \int\limits_{0}^{d_3^{-1}}O(1)d\xi=O(1) .\   \nonumber \
 \end{eqnarray*} 
 
 In order to estimate the integral $(5)$ in  Lemma \ref{L3.2} it is useful to combine the integrals over the segments of the curve $E=h$ situated  in a neighbourhood of the point $C$ in the $1$st and $4$th quadrants.  Denote  the corresponding  integral    $R^{(5)}$. Then 
 \begin{eqnarray*} 
{\cal {H}}=\frac12\omega_0(\eta^2-\xi^2)+\chi_3(\eta, \xi, z) , \ \chi_3=O(|\xi|^3+|\eta|^3),  \   \nonumber \\
\varphi=\alpha\xi+\beta\eta+\varphi_2(\eta, \xi, z) , \ \varphi_2=O(\xi^2+\eta^2),  \   \nonumber \
  \end{eqnarray*} 
   \begin{eqnarray*} 
R^{(5)}&=&\frac{\partial}{\partial h} \int\limits_{-d_3^{-1}}^{d_3^{-1}}\left(\frac{\varphi}{\partial{\cal {H}}/\partial\eta}\right)_hd\xi=
 \int\limits_{-d_3^{-1}}^{d_3^{-1}}\left(\frac{\partial}{\partial\eta}\left(\frac{\varphi}{\partial{\cal {H}}/\partial\eta}\right)\frac{1}{\partial{\cal {H}}/\partial\eta}\right)_hd\xi  \   \nonumber \\ 
 &=& \int\limits_{-d_3^{-1}}^{d_3^{-1}}\left(\frac{(\beta+\partial\varphi_2/\partial\eta)(\omega_0\eta+\partial\chi_3/\partial\eta)-(\alpha\xi+\beta\eta+\varphi_2)(\omega_0+\partial^2\chi_3/\partial\eta^2)}{(\partial{\cal {H}}/\partial\eta)^3}\right)_hd\xi  \   \nonumber \\ 
&=&-\int\limits_{-d_3^{-1}}^{d_3^{-1}}\frac{\omega_0\alpha\xi d\xi}{(\partial{\cal {H}}/\partial\eta)^3_h}+ \int\limits_{-d_3^{-1}}^{d_3^{-1}}\frac{O(1)}{\sqrt{h+\xi^2}}d\xi=-\omega_0\alpha\int\limits_{-d_3^{-1}}^{d_3^{-1}}\frac{\xi d\xi}{(\partial{\cal {H}}/\partial\eta)^3_h}+O(\ln h) \   \nonumber \\
&=&-\omega_0\alpha\int\limits_{-d_3^{-1}}^{d_3^{-1}}\frac{\xi d\xi}{(\omega_0\sqrt{2\omega_0^{-1}h+
\xi^2})^3}\\
&-&\omega_0\alpha\int\limits_{-d_3^{-1}}^{d_3^{-1}}\xi\left(\frac{1}{(\omega_0\eta+\partial\chi_3/\partial\eta)^3_h}-\frac{1}{(\omega_0\sqrt{2\omega_0^{-1}h+\xi^2})^3}\right)d\xi+ 
O(\ln h)= \   \nonumber \\
&=& \int\limits_{-d_3^{-1}}^{d_3^{-1}}\xi\frac{O(\eta^4)d\xi}{(h+\xi^2)^3}+O(\ln h)= \int\limits_{-d_3^{-1}}^{d_3^{-1}}\frac{O(1)d\xi}{\sqrt{h+\xi^2}}+O(\ln h)=O(\ln h) .\   \nonumber \
 \end{eqnarray*} 
Combining  estimates for different segments, we get  estimates of  Lemma \ref{L3.2}.
 \medskip
 
 \subsubsection{Proof of  Lemma \ref{L3.3}}
  \label{sss3.4.4}
  Let $0<h<\frac12d_4^{-1}$. We have 
  $$
  T=T(h)=\oint\limits_{E=h}dt\,.
  $$
  Let us split $T(h)$ into two  integrals,  $T'(h)$ and $T''(h)$, calculated  over the segments of the curve $E=h$ situated  respectively inside the domain $|\xi|<d_3^{-1}, |\eta|<d_3^{-1}$  and  outside this domain. Let us, in addition, split  $T'(h)$ into  integrals  $T'_{1,4}(h)$ and $T'_{2, 3}(h)$, calculated  over the segments  situated  in the $4$th and $1$st and, respectively, in the $2$nd and $3$rd  quadrants of the system $C\xi\eta$. 
Integral  $T''(h)$ is calculated over the segments of the curve separated from the point $C$. Therefore,  $T''(h)= T''(0)+O(h)$.
To estimate  $T'_{4, 1}(h)$, let us introduce $\xi$ as an independent variable:
 \begin{eqnarray*} 
T'_{4, 1}(h)&=& \int\limits_{-d_3^{-1}}^{d_3^{-1}}\frac{d\xi}{(\partial{\cal {H}}/\partial\eta)_h}\\
&=&\int\limits_{-d_3^{-1}}^{d_3^{-1}}\frac{d\xi}{\omega_0\sqrt{2\omega_0^{-1}h+\xi^2}}+ \int\limits_{-d_3^{-1}}^{d_3^{-1}}\frac {\omega_0\sqrt{2\omega_0^{-1}h+\xi^2}\,-(\partial{\cal {H}}/\partial\eta)_h}{\omega_0\sqrt{2\omega_0^{-1}h+\xi^2}\ (\partial{\cal {H}}/\partial\eta)_h}d\xi. \   \nonumber \
  \end{eqnarray*}
 The integrand in the last integral remains  bounded as $h\to 0$. Acting as in the proof of   Lemma \ref{L3.1} we get that this integral can be calculated at $h=0$ (i.e. over the separatrix)  with an accuracy $O(h\ln h)$. Another integral, forming a part of $T'_{4, 1}(h)$, has an explicit form:
 \begin {eqnarray}  \label{3.22}  
\int\limits_{-d_3^{-1}}^{d_3^{-1}}\frac{d\xi}{\omega_0\sqrt{2\omega_0^{-1}h+\xi^2}}=\left. \frac{1}{\omega_0}\ln(\xi+\sqrt{2\omega_0^{-1}h+\xi^2}\ )\right|_{-d_3^{-1}}^{d_3^{-1}}\\
=- \frac{1}{\omega_0}\ln h+ \frac{1}{\omega_0}\ln\frac{\omega_0}{2}+\frac{2}{\omega_0}\ln(d_3^{-1}+\sqrt{d_3^{-2}+2\omega_0^{-1}h} \ ) \   \nonumber \\
=- \frac{1}{\omega_0}\ln h+ \frac{1}{\omega_0}\ln\frac{\omega_0}{2}+\frac{2}{\omega_0}\ln(d_3^{-1})+O(h)\, . \   \nonumber \
 \end{eqnarray}
 Analogous estimate holds for  $T'_{2, 3}(h)$. Combining these estimates we get the assertion of  Lemma \ref{L3.3}.
 
 The calculation of  asymptotic expansion of the function  $T_i(h)$, the period of the motion along the trajectory $E=h<0$ situated in the region $G_i, \ i=1, 2$,  is treated by the same method, but as independent variable near the saddle the variable $\eta$ is used. In this calculation the sum of integrals over separatrices in asymptotic expansions of $T_1$ and $T_2$ coincides with the sum of integrals over separatrices for asymptotic expansion of  the period of  motion in the the region $G_3$. The integral, which is calculated in this way, is reduced to the same form (\ref{3.22}); the only difference in that $h$ is replaced with $|h|$, and $\xi$ is replaced with $\eta$. Hence the assertion of the Remark to  Lemma \ref{L3.3} is valid.
 \medskip
 
 \subsubsection{Proof of  Lemma \ref{L3.4}}
  \label{sss3.4.5}
Let us denote $\psi_c$ the value of the function $\psi$ at the point $C$. Then 
  \begin{equation}  \label{3.23} 
 \frac{\partial}{\partial h}\frac{1}{T}\oint\limits_{E=h}\psi dt= \frac{\partial}{\partial h}\frac{1}{T}\oint\limits_{E=h}(\psi-\psi_c) dt\,.
  \end{equation}
Calculating the derivative and  making use of the estimates in  Lemma \ref{L3.2}, its corollary, and  corollary of Lemma \ref{L3.3}, we get the result of the Lemma \ref{L3.4}.
The transformation  (\ref{3.23}) can not be used in the problems where  the boundary of the domain contains several saddle points. But the result of  the lemma is valid in these cases too. Let us describe briefly the corresponding proof.
As in  Lemma \ref{L3.3} we can prove expansions
 \begin {eqnarray}
 T=\alpha\ln h+\chi, \ \chi=\chi(h, z)=O(1), \ \partial\chi/\partial h=O(\ln h) \, ,   \nonumber \\
 \oint\limits_{E=h}\psi dt=\beta\ln h+\mu, \ \mu=\mu(h, z)=O(1), \  \partial\mu/\partial h=O(\ln h) \, .  \nonumber \
 \end{eqnarray}
Now
 \begin {eqnarray}
  \frac{\partial}{\partial h}\,\frac{1}{T}\oint\limits_{E=h}\psi dt&=&\frac{1}{T^2}\left[\left( \frac{\beta}{h}+\frac{\partial\mu}{\partial h}\right)(\alpha\ln h+\chi)-(\beta\ln h+\mu)\left(\frac{\alpha}{h}+\frac{\partial\chi}{\partial h}\right)\right]= \   \nonumber \\
&=& \frac{1}{T^2} O( h^{-1})= O( h^{-1}\ln^{-2}h)  \,.   \nonumber \
 \end{eqnarray}
 \medskip
 
 \subsection{Proofs of  Lemmas on perturbed motion}
  \label{ss3.5}

 In this subsection   Lemma \ref{L3.5} is proved. The proof of  Lemma \ref{L3.6} is  completely analogous, and it is omitted.
 
In accordance with  Lemma \ref{L3.2} ($1 $)
  \begin {eqnarray}
 \oint\limits_{E=h}|\chi|dt<\nu,  \nonumber
  \end {eqnarray}
 where
 \begin {eqnarray}
    \chi=\chi(p, q, z, \eps)=\frac{\partial E}{\partial q}f_1+\frac{\partial E}{\partial p}f_2+\frac{\partial E}{\partial z}f_3 \, . \nonumber \
  \end {eqnarray}
Let us introduce $c_6=3\nu_1, \ c_7=2\max(c_2, c_4)$. Denote $\eta(t), \xi(t)$  the values of $\eta, \xi$ at the point $(p(t), q(t), z(t))$.
Let $c_6\eps<h(t')<c_7^{-1}$. Denote $t_{1*}$ the supremum of  moments of time $t_1$ such that for $t'\le t\le t_1$ the solution $(p(t), q(t), z(t))$ is defined and meets the conditions
 \begin {eqnarray}  \label{3.24} 
 |\xi(t)|<d_3^{-1}, \  |\eta(t)|<d_1^{-1}, \ z(t)\in B-\frac{3}{2}c_1^{-1}\, , \\
\frac12h(t')<h(t)<2h(t').\phantom{***********} \, \nonumber
 \end{eqnarray}
 Denote $\xi_{1*}=\xi(t_{1*})$. For  $t'\le t\le t_{1*}$ we have
 \begin{eqnarray}  
\dot\xi&=&\partial{\cal {H}}/\partial\eta +O(\eps)>d_2^{-1}\sqrt{\frac12h(t')+\xi^2}\, +O(\eps)>\frac12 d_2^{-1}\sqrt{\frac12h(t')+\xi^2} \nonumber\, ,\\
t-t'&=&\int\limits_{0}^{\xi(t)}\frac{d\xi}{\dot\xi}<2d_2\int\limits_{0}^{d_3^{-1}}\frac {d\eta}{\sqrt{\frac12h(t')+\xi^2}}=O(\ln h(t')) \, ,   \label{3.25} \\
 |z(t)-z(t')|&=&O(\eps\ln h(t')) \, ,   \nonumber \\
 |h(t)-h(t')|&\le&\eps\int_{t'}^{t}|\chi|dt<2d_2\eps\int\limits_{0}^{\xi(t)}\frac {|\chi| d\xi}{\sqrt{\frac12h(t')+\xi^2}}= \int\limits_{0}^{d_3^{-1}}O(1)d\eps= O(\eps) \, .   \nonumber\
 \end{eqnarray}
The obtained inequalities allow to get more accurate estimate for $h(t_{1*})$:
 \begin{eqnarray*}  
  h(t_{1*})&-&h(t')=\eps\int_{t'}^{t_{1*}}\chi dt=\eps\int\limits_{0}^{\xi_{1*}}\frac{\chi}{\partial{\cal {H}}/\partial\eta +O(\eps)}{d\xi}=  \   \nonumber \\
 &=&\eps\int\limits_{0}^{\xi_{1*}}\left(\frac{\chi}{\partial{\cal {H}}/\partial\eta}\right)_{\substack {E=h(t')\\z=z(t')\\ \eps=0\phantom{***}}}d\xi+ \eps\int\limits_{0}^{\xi_{1*}}\left[\left(\frac{\chi}{\partial{\cal {H}}/\partial\eta +O(\eps)}\right)_{\substack {E=h(t)\\z=z(t)}}-\left(\frac{\chi}{\partial{\cal {H}}/\partial\eta}\right) _{\substack {E=h(t')\\z=z(t')\\ \eps=0\phantom{***}}}\right]d\xi  \nonumber\, .
 \end{eqnarray*}

By means of (\ref{3.25}) we get that the integrand in the second integral is
 \begin{equation}  
 \frac{O(\eps)}{h(t')+\xi^2}+ \frac{O(\eps\ln h(t'))}{\sqrt{h(t')+\xi^2}}  \nonumber\, .
  \end{equation}
So
 \begin{equation}  
  h(t_{1*})-h(t')=\eps\int\limits_{0}^{\xi_{1*}}\left(\frac{\chi}{\partial{\cal {H}}/\partial\eta}\right)_{\substack {E=h(t')\\z=z(t')\\ \eps=0\phantom{***}}}d\xi+O(\eps^2/\sqrt{h(t')}\,) \nonumber\, .
  \end{equation}
In this expression the integral is calculated over the segment of the unperturbed trajectory in the unperturbed motion. Therefore
 \begin{equation}  
  |h(t_{1*})-h(t')|<\eps\nu_1\le\frac12\eps c_6, \ \frac12h(t')<h(t_{1*})<2h(t')\nonumber\, .
  \end{equation}
Therefore at $t=t_{1*}$ the conditions on $z, h, \eta$ in (\ref{3.24}) are satisfied as strict inequalities. By definition of  $t_{1*}$ there should be $\xi(t_{1*})=d_3$.

In the further motion the phase point evidently makes the curve that  is close to the unperturbed trajectory  $E=h(t'), z=z(t')$, and arrives at the segment  $\xi=d_3^{-1}, \ -d_1^{-1}<\eta<0$. Along this curve $h=h(t')+O(\eps), \ z=z(t')+O(\eps\ln h(t'))$. Along any part of this curve the change of $E$ with an accuracy $O(\eps^2\ln\eps)$ is equal to the integral of the function $\eps\chi(p, q, z, 0)$ over the segment of  the unperturbed trajectory and, therefore, this change does not exceed $\frac32\eps\nu_1$. Therefore in this motion $\frac12h(t')<h(t)<2h(t')$.

Further motion is considered in an analogous manner. The phase point arrives first  at the segment  $\xi=-d_3^{-1}, \ -d_1^{-1}<\eta<0$, then it makes the curve close to the unperturbed trajectory, arrives at  the segment  $\xi=-d_3^{-1}, \ 0<\eta<d_1^{-1}$ and, finally, at a moment of time $t''=t'+O(\ln h(t'))$ arrives at the ray $C\eta$ having $0<\eta<d_1^{-1}$. For $t'\le t\le t''$  estimates (\ref{3.6}) are satisfied and the estimate
 \begin{equation}  
  |j(t)-j(t')|\le\left(\frac{\partial I}{\partial h}\right)_{\substack {h=h_*\\z=z_*}}|h(t)-h(t')| +\left(\frac{\partial I}{\partial h}\right)_{\substack {h=h_*\\z=z_*}} |z(t)-z(t')|=O(\eps\ln h(t'))  \nonumber\
  \end{equation}
is valid (here $h_*\in(\frac12h(t'), 2h(t')), \ z_*=z(t')+O(\eps\ln h(t'))$).

For the function $\varphi(p, q, z, \eps)$, which vanishes identically at the point $C$, the integral along the motion is estimated in the same manner as above for the function $\chi$. This gives the second estimate (\ref {3.7}).

To estimate the integral  along the motion of the function $\psi$ in  (\ref {3.7}), we split it into integrals over the defined above segments of the trajectory, situated either far from or near the saddle. Integrals over the  segments, situated  far from the saddle, coincide with accuracy  $O(\eps\ln h(t'))$ with integrals over the  segments of the unperturbed trajectory $E=h(t')$. For the integral over the  segment of the perturbed trajectory with $0<\xi<d_3^{-1},\ \eta>0$, situated near the saddle, we have
 \begin{eqnarray*}  
  \int\limits_{0}^{t_{1*}}\psi dt=\int\limits_{0}^{d_3^{-1}}\left(\frac{\psi}{\partial{\cal {H}}/\partial\eta+O(\eps)}\right)_{\substack {E=h(t)\\z=z(t)}}d\xi=\int\limits_{0}^{d_3^{-1}}\left(\frac{\psi}{\partial{\cal {H}}/\partial\eta}\right)_{\substack {E=h(t')\\z=z(t')\\ \eps=0\phantom{***}}}d\xi+\  \nonumber\\
  + \int\limits_{0}^{d_3^{-1}}\left[\left(\frac{\psi}{\partial{\cal {H}}/\partial\eta +O(\eps)}\right)_{\substack {E=h(t)\\z=z(t)}}-\left(\frac{\psi}{\partial{\cal {H}}/\partial\eta}\right)_{\substack {E=h(t')\\z=z(t')\\ \eps=0\phantom{***}}}\right]d\xi \, .  \nonumber\
 \end{eqnarray*}
Taking into account already proved estimates  (\ref {3.6}) and Lemma  \ref {L3.7}, we estimate the integrand in the last integral as 
\begin{equation}  
  \frac{O(\eps)}{(h(t')+\xi^2)^{3/2}}+ \frac{O(\eps\ln h(t'))}{h(t')+\xi^2} \nonumber \,.
  \end{equation}
The integral of this function is $O(\eps/h(t'))$.
Therefore
\begin{equation}  
 \int\limits_{0}^{t_{1*}}\psi dt=\int\limits_{0}^{d_3^{-1}}\left(\frac{\psi}{\partial{\cal {H}}/\partial\eta}\right)_{\substack {E=h(t')\\z=z(t')\\ \eps=0\phantom{***}}}+O(\eps/h(t')) \nonumber \,.
  \end{equation}
The integral in the right hand side is just the integral of $\psi$  along the unperturbed trajectory $E=h(t'), z=z(t')$ when $\eps=0$. The integrals of $\psi$  along  other segments of the trajectory near the saddle are estimated in analogous manner. 
Combining   these estimates we get the first estimate  (\ref {3.7}).
 \medskip
 
 \subsection{Proofs of  Lemmas on averaged system}
  \label{ss3.6}

 In this subsection the proofs of   Lemmas \ref{L2.2} and \ref{L2.3} on the motion in averaged system are given. These proofs are based on Lemmas on unperturbed motion in Subsection \ref{ss3.1}.

 \subsubsection{Proof of   Lemma \ref{L2.2}}
  \label{sss3.6.1}
 a) Existence of the solution with initial condition $h=0, z=z_*$ {at}  $\tau=\tau_*$ follows from the standard existence theorem for ODEs (see, e.g., \cite {Hartman}, p. 21) as the right hand side of the averaged system is continuous.
 
To prove a uniqueness  let us suppose that there are two solutions, \\
 $(H^{(1)}(\tau), Z^{(1)}(\tau))$ and $(H^{(2)}(\tau), Z^{(2)}(\tau))$, say in the region $G_3$ (i.e. $\tau<\tau_*)$, crossing the separatrix at $\tau=\tau_*, z=z_*$. Denote $J^{(i)}(\tau)=I(H^{(i)}(\tau), Z^{(i)}(\tau))$,  $ i=1, 2$.  Denote $u(\tau)=|(J^{(2)}(\tau)- J^{(1)}(\tau)|+|(Z^{(2)}(\tau)- Z^{(1)}(\tau)|$. Suppose that $u(\tau_{**})\ne 0$ for  some $\tau_{**}<\tau_*$. Then $u(\tau)\ne 0$ for  $\tau_{**}<\tau<\tau_{*}$  (in the opposite case we have  a contradiction with the standard uniqueness theorem). We may assume that $H^{(i)}(\tau_{**})<1/2,\, i=1,2$.
 From the formulas for the  right hand side of the averaged equations for $z, J$,  (\ref {2.3}),  (\ref {2.7}) and from Lemmas  \ref{L3.1} - \ref{L3.4} we get that for $\tau_{**}<\tau<\tau_{*}$
\begin{equation} 
\frac{du(\tau)}{d\tau}=O\left(\frac{1}{H_*(\tau)\ln^3H_*(\tau)}\right)u(\tau)  \nonumber \, ,
 \end{equation}
where $H_*(\tau)=\min\{H^{(1)}(\tau), H^{(2)}(\tau)\}$. Therefore for $\tau_1\in(\tau_{**}, \tau_{*})$ we have
\begin{equation} 
\nonumber
u(\tau_{**})=u(\tau_1)\exp\left(\int\limits_{\tau_1}^{\tau_{**}}O\left(\frac{1}{H_*(\tau)\ln^3H_*(\tau)}\right)d\tau\right) \,. 
 \end{equation}

But  according to the averaged equation for $h$ (\ref {2.3}), estimate (\ref {3.2}), and  Lemma  \ref{L3.2} we have
$$
\frac{dH^{(i)}}{d\tau}<\nu_1^{-1}/\ln H^{(i)}\,.
$$
Therefore
\begin{equation}  
\int\limits_{\tau_1}^{\tau_{**}}O\left(\frac{1}{H_*(\tau)\ln^3H_*(\tau)}\right)d\tau=O\left(\int\limits_{0}^{1/2}\frac{dh}{h\ln^2h}\right)=O(1)  \ \nonumber \
 \end{equation}
and
\begin{equation}  \label{3.26} 
U(\tau_{**})=U(\tau_1)O(1)\,.
 \end{equation}
But $U(\tau_1)\to 0$ as $\tau_1\to\tau_*$. Therefore $U(\tau_{**})=0$ in contradiction with our hypothesis. The uniqueness is proved.

b) Assertion b) of the  Lemma  \ref{L2.2} is an evident corollary of formula  (\ref {2.7}) and condition  (\ref {2.2}), as it was discussed at the end of Subsection \ref{ss2.3}.

 \subsubsection{Proof of   Lemma \ref{L2.3}}
  \label{sss3.6.2}
We will omit index ``$\nu$" at solutions $(H_\nu, Z_\nu)$ and  $(H'_\nu, Z'_\nu)$. There exists a moment of the slow time $\bar\tau_0>\tau_0, \bar\tau_0=\tau_0+O(\delta\ln\delta)$ such that
$$
-\nu_1\delta <H(\bar\tau_0)<-\delta, \ -\nu_1\delta <H'(\bar\tau_0)<-\delta.
$$
Making use of the formula
\begin{equation}  \label{3.27} 
\frac{dZ}{d\tau}=f_{3C}(Z)+\frac{1}{T}\oint\limits_{E=h}(f_3^0-f_{3C})dt
 \end{equation}
we get that
$$
Z(\bar\tau_0)-Z'(\bar\tau_0)=O(\delta).
$$
For $\tau_0\le\tau\le\bar\tau_0$ we have
\begin{equation}  \label{3.28} 
|H(\tau)|+|H'(\tau)|+|Z(\tau)- Z'(\tau)|=O(\delta).
 \end{equation}
Therefore, the estimates of   Lemma \ref{L2.3} are valid for $\tau_0\le\tau\le\bar\tau_0$.
Denote
 \begin{eqnarray*}
 J(\tau)&=&I(H(\tau), Z(\tau)), \quad  J'(\tau)=I(H'(\tau), Z'(\tau)),\\ 
U(\tau)&=&|J(\tau)-J'(\tau)|+|Z(\tau)-Z'(\tau)|.
\end{eqnarray*}
Then $U(\bar\tau_0)=O(\delta\ln\delta)$.

Exactly as in (\ref {3.26}), we get  that $U(\tau)=O(U(\bar\tau_0))$ for $\bar\tau_0\le\tau\le K$.
So
\begin{equation}  \label{3.29} 
|J(\tau)-J'(\tau)|=O(\delta\ln\delta), \ |Z(\tau)-Z'(\tau)|=O(\delta\ln\delta).
 \end{equation}
From here we get that
\begin{equation}  \label{3.30} 
|H(\tau)-H'(\tau)|<\nu_2\frac{\delta|\ln\delta|}{1+|\ln|H_*(\tau)||}, \ 
 \end{equation}
where $H_*(\tau)$ lies between $H(\tau)$ and $H'(\tau)$.

Let us consider two cases: $H(\tau)\ge-2\nu_2\delta|\ln\delta|$ and $H(\tau)<-2\nu_2\delta|\ln\delta|$.

If $H(\tau)\ge-2\nu_2\delta|\ln\delta|$, then, from  (\ref {3.30}), $H'(\tau)>-\nu_2\delta|\ln\delta|, H_*(\tau)>-\nu_2\delta|\ln\delta|$ and
\begin{equation}  \label{3.31} 
|H(\tau)-H'(\tau)|=O(\delta). \
 \end{equation}
 
If $H(\tau)<-2\nu_2\delta|\ln\delta|$, then, from  (\ref {3.30}),  $\frac{1}{2}|H(\tau)|<|H'(\tau)|<2|H(\tau)|$, $\frac{1}{2}|H(\tau)|<|H_*(\tau)|<2|H(\tau)|$ and, therefore,
\begin{equation}  \label{3.32} 
|H(\tau)-H'(\tau)|=O\left(\frac{\delta\ln\delta}{1+|\ln|H(\tau)||}\right). \ 
 \end{equation}
In view of (\ref {3.28}), (\ref {3.31}), (\ref {3.32}) we get for $\tau_0\le\tau\le K$
\begin{equation} 
|H(\tau)-H'(\tau)|=O\left(\delta+\frac{\delta|\ln\delta|}{1+|\ln|H(\tau)||}\right), \ \nonumber \
 \end{equation}
i.e. the estimate of    Lemma \ref{L2.3} for $|H(\tau)-H'(\tau)|$.

\medskip

Now we should improve the estimate (\ref {3.29}) for $Z$. From  equality (\ref {3.27}), considered for $Z(\tau)$ {and}  $Z'(\tau)$,  Lemma \ref{L3.4}, and  (\ref{3.29}), (\ref{3.30}), making use of Lagrange formula, we get  for $\bar\tau_0<\tau< \bar\tau_0+\nu_3^{-1}$
$$
\frac{d}{dt}(Z-Z')=O(Z-Z')+O\left(\frac{\delta\ln\delta}{H_*(\tau)\ln^3|H_*(\tau)|}\right),
$$
where $H_*(\tau)$ lies between $H(\tau)$ and $H'(\tau)$.

From here  we get for $\bar\tau_0\le\tau\le \bar\tau_0+\nu_3^{-1}$
\begin{equation} 
Z(\tau)-Z'(\tau)= O(\delta)+
O\left(\int\limits_{\bar \tau_0}^{\tau}\frac{\delta\ln\delta \,d\tau}{H_*(\tau)\ln^3|H_*(\tau)|} \right)
= O(\delta)+ O\left(\int\limits_{\delta}^{|H_*(\tau)|}\frac{\delta\ln\delta \, dh}{h\ln^2h}\right) ,\ \nonumber 
 \end{equation}
 or
\begin{equation}  \label{3.33} 
|Z(\tau)-Z'(\tau)|=O\left(\delta+\frac{\delta|\ln\delta|}{|\ln|H_*(\tau)||}\right) . \  \
 \end{equation}
We know that  for $\bar\tau_0+\nu_3^{-1}\le\tau\le K$ the following estimate is valid:
\begin{equation}  \label{3.34} 
|Z(\tau)-Z'(\tau)|=O(\delta\ln\delta).
\end{equation}
 In view of  (\ref{3.28}), (\ref{3.33}), (\ref{3.34}) we have for $\tau_0\le\tau\le K$
\begin{equation} 
|Z(\tau)-Z'(\tau)|=O\left(\delta+\frac{\delta|\ln\delta|}{1+|\ln|H(\tau)||}\right).  \ \nonumber \
 \end{equation}
  Lemma \ref{L2.3} is proved.
  \medskip
   \section{Passage through a narrow vicinity of  separatrices}
  \label{S4}
  \medskip
In this section  Proposition \ref{Pr2.3} of  Subsection \ref{ss2.5} is proved. 
    \medskip
   \subsection{Preliminary transformation}
  \label{ss4.1}
\begin{figure}
 \begin{center}
            \includegraphics[scale=0.4]{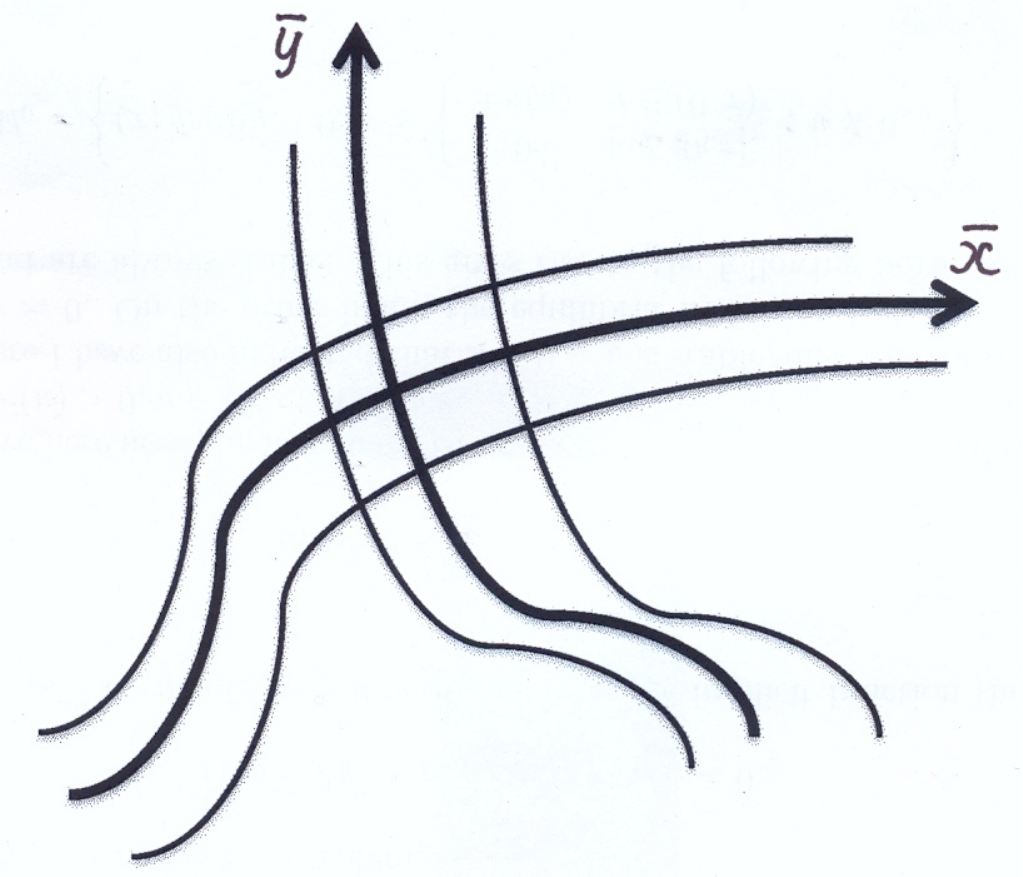}
            \end{center}
           \caption{New coordinates. }
            \label{xy_coordinates} 
\end{figure}

 Below in all estimates $z\in B-c_1^{-1}$, where the constant $c_1$ is chosen in such a way that  $3c_1^{-1}$-neighbourhood of the set $\{z\colon z=Z_\nu(\tau), \ \nu=1, 2; \ \tau \in[0, K]\}$ belongs to $B$, and  $3c_1^{-1}$-neighbourhood of the set $\{(p, q, z)\colon z=Z_\nu(\tau), \ E(p, q, z)=H_\nu(\tau), \ \nu=1, 2;\ \tau\in[0, K]\}$ belongs to $D$. 
In some neighbourhood of the  point $C$ we can introduce  new coordinates $\bar x, \bar y$  instead of $p, q$ in such a way  that the equations $\bar x=0$ and $\bar y=0$ define separatrices as it is shown in Fig. \ref{xy_coordinates},  and $\left.\partial(p, q)/\partial(\bar y,\bar x)\right|_C=1$ (we do not need the transformation $(p, q)\mapsto(\bar y,\bar x)$  be symplectic).
In the new coordinates the function $E$ has the form
 \begin{eqnarray*}
E=-\omega_0\bar y\bar x(1+O(|\bar y|+|\bar x|)), \quad 
\omega_0=\omega_0(z)>c_2^{-1}.\  \ \nonumber \
\end{eqnarray*}
 \medskip
 \begin{lemma} 
        \label{L4.1}
For   $z\in B-c_1^{-1}, \ |\bar y|<c_3^{-1}, |\bar x|<c_3^{-1}$  there  exists  a smooth transformation of variables $\bar F\colon y, x, z\to \bar y,\bar x,  z$ such that $\bar y=y+O(\eps), \ \bar x=x+O(\eps)$, and in the  variables $y, x$ the perturbed motion is described by equations  
    \begin{eqnarray*}    
    \dot x=-\omega_0 x(1+O(\eps+|x|+|y|)), \  
    \dot y=\omega_0 y(1+O(\eps+|x|+|y|)). \  \nonumber \
\end{eqnarray*}
\end{lemma} 

This Lemma is a direct corollary of a theorem by N. Fenichel (see \cite {Fenichel}, Theorem $11.1$). Coordinates of such type as 
$y, x$ here are often called   Fenichel's  coordinates.

 \medskip
Lemma \ref{L4.1} allows to define in the domain $ |\bar y|<c_3^{-1}, |\bar x|<c_3^{-1}$ the transformation of variables  $F\colon y, x\to p, q$ depending on $z$ as a parameter. The inverse transformation  $F^{-1}\colon p, q\to y, x$ is defined in the $c_4^{-1}$-neighbourhood of  the point $C$ in the plane $q, p$.
  \medskip
   \subsection{Lemmas on perturbed motion}
  \label{ss4.2}
\begin{figure}
 \begin{center}
            \includegraphics[scale=0.35]{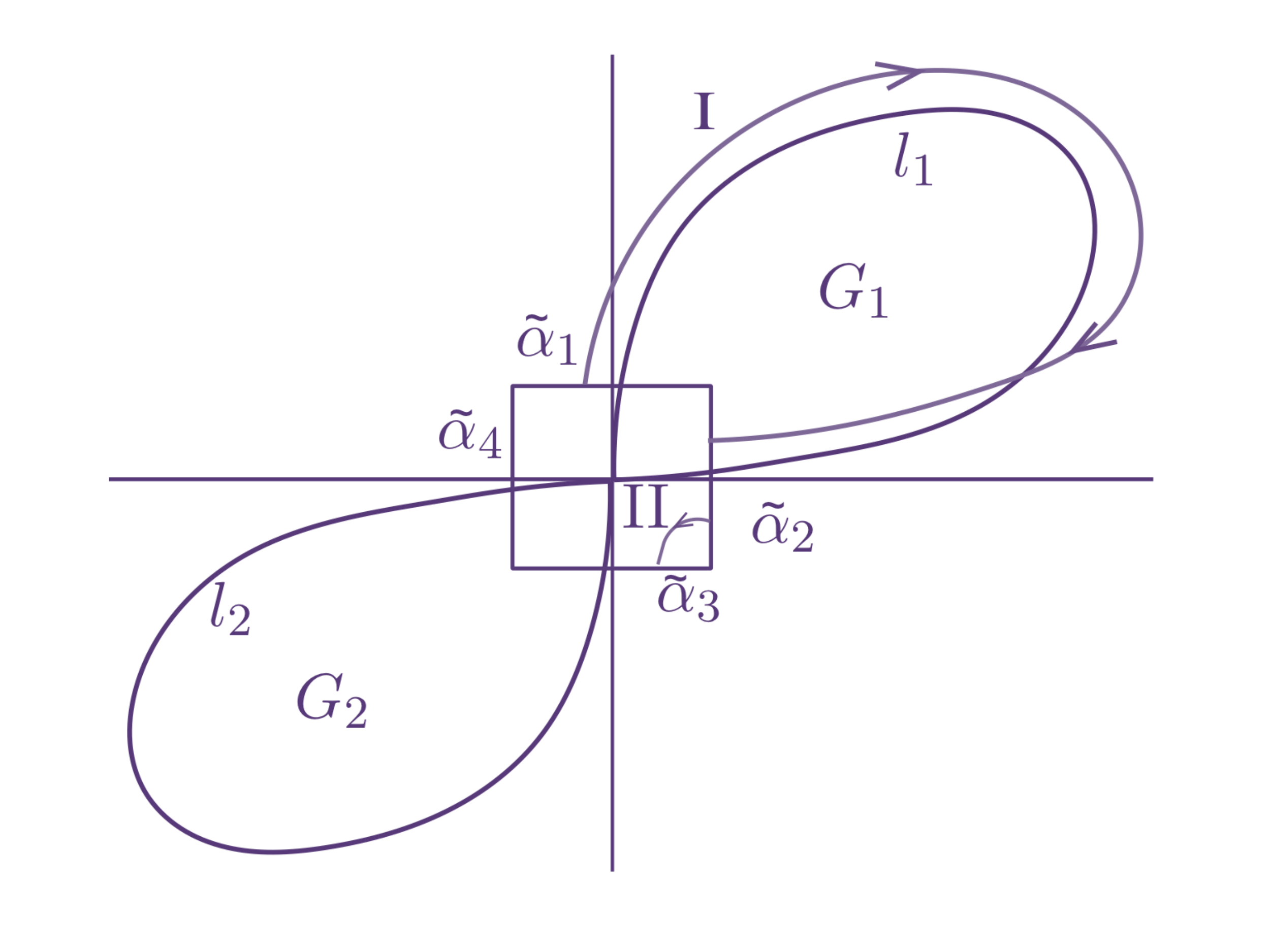}
            \end{center}
           \caption{Trajectories I and II. }
            \label{square+trajectories} 
\end{figure}

Consider in the plane $x, y$ a square ${\cal K}=\{x, y\colon |x|<2c_6\sqrt{\eps}, \ |y|<2c_6\sqrt{\eps},\}$. The constant $c_6$ is chosen below in Lemma \ref{L4.2} and Subsection \ref{ss4.3}.
\begin{figure}
 \begin{center}
            \includegraphics[angle=-0.5, scale=0.5]{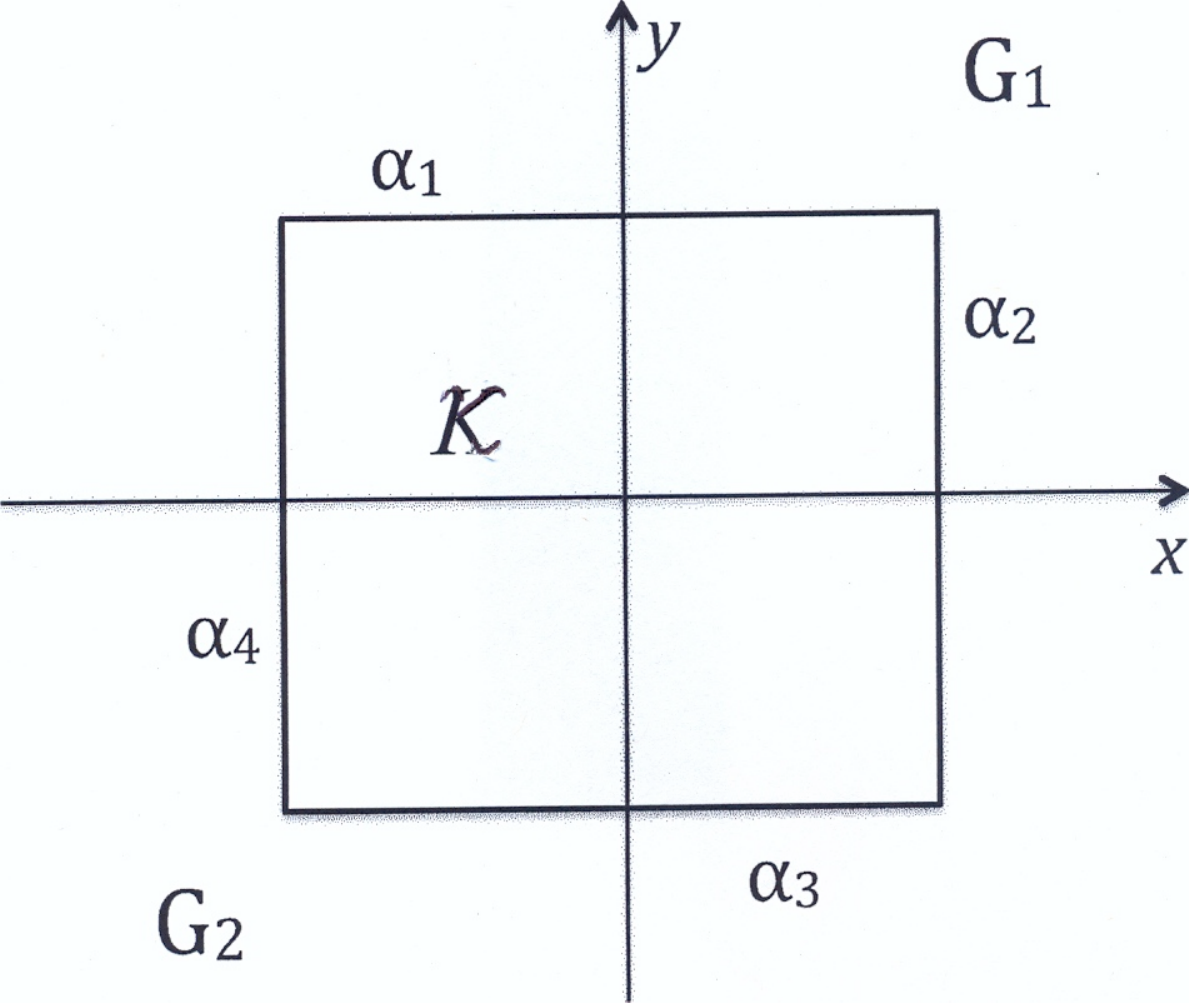}
            \end{center}
           
           \caption{The square ${\cal K}$. }
            \label{square} 
\end{figure}  
 Let $\alpha_1, ...\alpha_4$ be sides of this square, the numeration of the sides corresponds to Fig. \ref{square}. Denote $\tilde{\cal {K}}=\tilde{\cal {K}}(z)=F[{\cal {K}}], \ \tilde\alpha_i=\tilde\alpha_i(z)=F[\alpha_i], \ i=1, ..., 4$. For a point $(p(t), q(t))$, that belongs  to $c_4^{-1}$-neighbourhood of the point $C$, we denote $(y(t), x(t))=F^{-1}(p(t), q(t))$. 
The following two Lemmas show that points that started to move at  sides of the square $ \tilde{\cal {K}}$ far from its corners and not too close to the $x$-axis, will return to the sides of $ \tilde{\cal {K}}$.
 In  Fig. \ref{square+trajectories} the trajectories I and II correspond to Lemmas   \ref{L4.2} and  \ref{L4.3}, respectively.
 \medskip
 \begin{lemma} 
        \label{L4.2}
If  $z(t')\in B-2c_1^{-1}, \ (p(t'), q(t'))\in\tilde\alpha_i(z(t')), i=1, 3, \  |x(t')|<c_6\sqrt{\eps}$, then there exists a moment of time $t''$ such that
  \begin{eqnarray*}    
   (p(t''), q(t''))\in\tilde\alpha_{i+1}(z(t'')), \ t''-t'=O(\ln\eps),  \  \nonumber \
h(t')-h(t'')>c_7^{-1}\eps   \,. \nonumber \
\end{eqnarray*}
\end{lemma} 
{\bf{Remark}}. If we increase the value of the constant $c_6$, it would  result in increasing of  values of the constants $c_j, \, j>7$. The constant $c_7$ can be chosen independently of $c_6$.

 \medskip
 \begin{lemma} 
        \label{L4.3}
Let  $z(t')\in B-2c_1^{-1}, \ (p(t'), q(t'))\in\tilde\alpha_2(z(t'))\cup\tilde\alpha_4(z(t')),\  |y(t')|>\eps^{r_1}$, where $r_1$ is any given in advance positive number. If  $ (p(t'), q(t'))\in\tilde\alpha_2(z(t')), \ y(t')>\eps^{r_1}$, then there exists a moment of time $t''$ such that
 $$    
   (p(t''), q(t''))\in\tilde\alpha_1(z(t'')).
     $$
     In analogous way,
   \begin{eqnarray*}     
     \mbox {if} \  (p(t'), q(t'))\in\tilde\alpha_2(z(t')), \  y(t')<-\eps^{r_1}, \ \mbox{then} \ (p(t''), q(t''))\in\tilde\alpha_3(z(t'')),  \  \nonumber \\
      \mbox {if} \  (p(t'), q(t'))\in\tilde\alpha_4(z(t')), \  y(t')>\eps^{r_1}, \ \mbox{then} \ (p(t''), q(t''))\in\tilde\alpha_1(z(t'')),\phantom{**}  \  \nonumber \\
 \mbox {if} \  (p(t'), q(t'))\in\tilde\alpha_4(z(t')), \  y(t')<-\eps^{r_1}, \ \mbox{then} \ (p(t''), q(t''))\in\tilde\alpha_3(z(t'')).  \  \nonumber \
\end{eqnarray*}
In all cases $t''-t'=O(\ln\eps), \ h(t'')-h(t')=O(\eps^{3/2})$.
\end{lemma}
 
The following assertion shows that a phase point that  stared to move from a diagonal of the square  $ \tilde{\cal {K}}$ not too close to the origin of the coordinate system will arrive to a  side of $ \tilde{\cal {K}}$.
 \medskip
 \begin{lemma} 
        \label{L4.4}
If $z(t')\in B-2c_1^{-1}, \ (p(t'), q(t'))\in \tilde{\cal {K}}(z(t')), \  y(t')=-x(t')>\eps$, then there exists a moment of time $t''$ such that
 \begin{eqnarray*}    
   (p(t''), q(t''))\in\tilde\alpha_1(z(t'')), \  \nonumber \
   t''-t'=O(\ln\eps), \ h(t'')-h(t')=O(\eps^{3/2}).  \  \nonumber \
   \end{eqnarray*}
   \end{lemma} 

Solutions of the averaged system with initial conditions from $W^\delta$ at $\tau=0$ cross the separatrix at $\tau,  z$ such that $|\tau-\hat\tau_*|<c_8\delta, \ |z-\hat z_*|<c_9\delta$. Consider the surface $\Lambda$ in the space $x, y, z$:
$$
\Lambda=\{x, y, z\colon|x|=2c_6\sqrt{\eps}, \ |z-\hat z_*|<2c_9\delta, \ |y|<\eps^{r_1}\}.
$$
Let  $\tilde\Lambda$ be the image of this surface under  the map $x, y, z\mapsto p, q, z$.
Denote $\Xi$ the set in the space $p, q, z$ sweept  by  $\tilde\Lambda$ during the shift along the trajectories of the perturbed system (\ref{perturbed}) for time $t\in(-\hat\tau_*-2c_8\delta, -\hat\tau_*+2c_8\delta)$ before the first arrival to the boundary of the region $D-c_1^{-1}$.
 \medskip
 \begin{lemma} 
        \label{L4.5}
$$
 \mbox{\rm mes}\, \,\Xi=O(\eps^{r_1-1/2} \mbox{\rm mes}\,W^{\delta}).
 $$
  \end{lemma} 
These Lemmas are proved in Subsections \ref{sss4.4.2}, \ref{sss4.4.3}
 \medskip
   \subsection{Proof of  Proposition \ref{Pr2.3}} 
  \label{ss4.3}
Define $r_1=r+1/2$, where  $r$ is the integer number introduced before  the formulation of Theorem \ref{T2.1}. Introduce $w=W^\delta\cap\Xi, \ W_3^\delta=W^\delta\setminus w$, where $\Xi$ is defined at the end of Subsection \ref{ss4.2}. According to Lemma \ref{L4.5}  $\mbox{\rm mes}\,w=O(\eps^r\mbox{mes}\,W^\delta)$.
Consider motion of a phase point $(p(t), q(t), z(t))$ that starts  in $W_3^\delta$ at the moment  of time  $t=0$. In accordance with  Proposition \ref{Pr2.1}  the moment of time $t_-$ is defined   such that $h(t_-)=2k_1\eps$. Then there exists a moment of time $t'_-$ preceding  $t_-$ such that the point $(p(t'_-), q(t'_-))$ lies in the $c_4^{-1}$-neighbourhood of the saddle $C$ for $z=z(t'_-)$,  and $y(t'_-)=-x(t'_-)>0$, $c_5\eps<h(t'_-)<2c_5\eps$. The proof of the last assertion is completely analogous to the proofs of  Lemmas \ref{L3.5}, \ref{L3.6} and  is omitted here.

Now let us choose the constant $c_6$ that defines the  size of the square ${\cal {K}}$  in Subsection \ref{ss4.2} (here we use the Remark to   Lemma \ref{L4.2}, which allows to increase the value of $c_6$). Choose $c_6$  such that in the segments $\alpha_1$ and $\alpha_3$ the inequality $-2k_2\eps<E<4c_5\eps$ implies   the inequality $|x|<c_6\sqrt{\eps}$ for $z\in B-c_1^{-1}$.
In the segments $\alpha_1$ and $\alpha_3$  according  to  Lemma \ref{L4.1} we have $|E|=2\omega_0c_6\sqrt{\eps}\,|x|+O(\eps^{3/2}), \ \omega_0>c_2^{-1}$. 
Choose $c_6^2>4c_2\max(2c_5, k_2)$. Then for $|x|>c_6\sqrt{\eps}$ we have $|E|>\max(4c_5\eps, 2k_2)$. So, this choice of $c_6$ meets   our condition. It is easy to check that under the same  choice of $c_6$   the inequality  $|E|>2c_5\eps$ is satisfied on the line $x=-y$ outside of  ${\cal {K}}$.

 Let us return to the motion of the point  $(p(t), q(t), z(t))$. It is evident that $y(t')>\eps$ (otherwise  $h(t')=O(\eps^2)<c_5\eps$). In accordance with Lemma \ref{L4.4} there exists a moment of time $t''_-=t'_-+O(\ln\eps)$ such that  $ (p(t''_-), q(t''_-))\in\tilde\alpha_1(z(t''_-))$ and  $ h(t''_-)=h(t'_-)+O(\eps^{3/2})$. In particular $\frac12c_5\eps<h(t''_-)<3c_5\eps$ and, therefore, $|x(t''_-)|<c_6\sqrt{\eps}$ due to the choice of $c_6$.
  Now Lemmas \ref{L4.2}, \ref{L4.3} allow, while their hypotheses are satisfied, to definite inductively  the moments of time $t_s$ of consecutive arrivals of our phase point to the sides $\tilde\alpha_1$ {or} $\tilde\alpha_3$ of the square; \ $t_1=t''_-$.
 In accordance with Lemmas \ref{L4.2}, \ref{L4.3}, $h(t_s)-h(t_{s+1})>\frac12c_7^{-1}\eps$. As $|E|<c_{10}\sqrt{\eps}$ in the square  ${\cal {K}}$, so $s<4c_7c_{10}$. As $t_{s+1}-t_s=O(\ln\eps)$, so $z(t_s)-z(t_1)=O(\ln\eps)$ and, therefore, $z(t_s)\in B-\frac52c_1^{-1}$.
 If $|x(t_s)|<c_6\sqrt{\eps}$, then the moment of time $t_{s+1}$ do exists (here we use the following: because of the choice of initial conditions and the definition of the set $w$, the hypothesis $|y|>\eps^{r_1}$ of   Lemma \ref{L4.3} can not be violated). But the sequence $\{t_s\}$ should be finite: $s<4c_7c_{10}$. Therefore, there exists a number $s_*$ such that $|x(t_{s_*})|\ge c_6\sqrt{\eps}$. Then because of the choice of $c_6$ it  should be $h(t_{s_*})<-k_2\eps$. By continuity there exists a moment of time $t_+\in(t_-, t_{s_*})$ such that $h(t_+)=-k_2\eps$, as it was stated in Proposition  \ref{Pr2.3}.
 
 Estimate $h(t), H_\nu(\eps t), z(t),Z_{\nu}(\eps t)$ for $t_-\le t\le t_+$. For  $t\in (t_-, t_+)$ we have $h(t)=O(\eps)$ by definitions of  $ t_-, t_+$.
  As $H_\nu(\eps t_-)=O(\eps), \  t_+-t_- =O(\ln\eps)$ and, in correspondence with estimates of  Lemmas \ref{L3.2}, \ref{L3.3}, $\dot H_\nu= O(\eps/\ln |H_\nu|) \ \mbox{for} \ H_\nu\ne 0, \ \mbox{so} \ H_\nu(\eps t)=O(\eps)\  \mbox{for}\ t\in (t_-, t_+)$.
  As $\dot z=O(\eps), \dot Z_{\nu}=O(\eps)$, so for $  t\in (t_-, t_+)$  we have
  $$
 z(t)-Z_\nu(\eps t)=O(\eps\ln\eps).
 $$ 
Improve the last estimate. From identity (\ref{3.16}) we get
 \begin{equation}     
  \dot z-\dot Z_\nu=\eps O(|z-Z_\nu|)+\eps O(1/\ln|H_\nu|)+\eps f_3^0(p, q, z)- f_{3C}^0( z)+O(\eps^2), \  \nonumber \
   \end{equation} 
where $ f_{3C}^0$ is the value of the function $ f_3^0(p, q, z)$ at the point $C$. Calculate integrals of the left and right sides of this relation with respect to time from $t_s$ to some $t\in (t_s, t_{s+1})$. Making use of already established estimates  for $H_\nu, \ z-Z_\nu$, we get
 \begin{equation}     
   z(t)-Z_\nu(\eps t)= z(t_s)-Z_\nu(\eps t_s)+O(\eps)+ \eps\int_{t_s}^t( f_3^0(p, q, z)- f_{3C}^0( z))dt .\  \nonumber \
   \end{equation} 
The integrand in the last integral vanishes at the point $C$. Analogously to Lemma \ref{L3.5}, this integral is $O(1)$. Therefore, for  $t\in (t_s, t_{s+1})$ we have
 \begin{equation}     
   z(t)-Z_\nu(\eps t)= z(t_s)-Z_\nu(\eps t_s)+O(\eps) .\  \nonumber \
   \end{equation} 
Analogous estimate holds for $t_1\le t\le t_-\,$. As $s<4c_7c_9$ and,  in accordance with  Proposition  \ref{Pr2.1}, $z(t_-)-Z_\nu(\eps t_-)=O(\eps), \  \mbox{so} \  z(t)-Z_\nu(\eps t)= O(\eps) \ \mbox{for} \ t_-\le t\le t_+$. This  was the assertion of  Proposition  \ref{Pr2.3}.

\medskip
{\bf {Remark}}. In the proof of the last estimate the representation (\ref{3.16}) was used. We can not use this representation in  problems, where separatrices connect different saddle points, and at these points the function $ f_3^0(p, q, z)$ has different values. The accuracy of description  of $z$ in these cases is, in general, given by proved above estimate $z(t)-Z_\nu(\eps t)= O(\eps\ln\eps)$ (see also Remark at the end of Subsection \ref{ss3.3}, part III).
 \medskip
   \subsection{Proofs  of   Lemmas on motion in a narrow  vicinity of separatrices} 
  \label{ss4.4.}
  \subsubsection{Proof of   Lemma \ref{L4.2}} 
  \label{sss4.4.1}
 Let, for certainty, $(p(t'), q(t'))\in\tilde\alpha_1(z(t'))$. Then for any choice of  $c_6$ there exists a moment of time $t''>t', \ t''=t'+O(\ln\eps)$ such that  the point    $(p(t''), q(t''))$ lies  in the $c_4^{-1}$-neighbourhood of the saddle $C$, and $|x(t'')|=2c_6\sqrt{\eps}, \ h(t'')=h(t')-\eps\Theta_2(z(t'))+O(\eps^{3/2})$. The proof of this assertion is analogous to the proof of   Lemma \ref{3.5}, and we omit it. Last estimate shows, in particular, that $c_7^{-1}\eps\le h(t')-h(t'')\le d_1\eps$, where $c_7$ and $d_1$ do not depend on choice of $c_6$.
 As  $h(t')=O(\eps)$, so  $h(t'')=O(\eps)$ and, therefore,  $y(t'')=O(\sqrt{\eps})$. Now
 \begin{eqnarray*}    
  h(t')=-2\omega_0c_6\sqrt{\eps}x(t')+O(\eps^{3/2}) , \  \nonumber \\
    h(t'')=-2\omega_0c_6\sqrt{\eps}y(t'')+O(\eps^{3/2}).  \  \nonumber \
 \end{eqnarray*}
From here
 \begin{equation}     
   y(t'')-x(t')=\frac{h(t')-h(t'')}{2\omega_0c_6\sqrt{\eps}}+O(\eps). \  \nonumber \
   \end{equation} 
Choose $c_6^2>d_1c_2$. Then 
$$
|  y(t'')-x(t')|<\frac{d_1c_2}{2} \frac{\sqrt{\eps}}{c_6}+O(\eps) <\frac12c_6\sqrt{\eps}+O(\eps)<c_6\sqrt{\eps}.
$$
Therefore, if $|x(t')|<c_6\sqrt{\eps}$, then $|y(t'')|<2c_6\sqrt{\eps}$. This means that  $(p(t''), q(t''))\in\tilde\alpha_2(z(t''))$. This was the assertion of the Lemma \ref{L4.2}.

\medskip
 \subsubsection{Proofs of   Lemmas \ref{L4.3},  \ref{L4.4}} 
  \label{sss4.4.2}
 Let, for certainty, $(p(t'), q(t'))\in\tilde\alpha_2(z(t')), \ y(t')>\eps^{r_1}$. Denote $t''$ the supremum of the moments of time $\bar t$ such that for $t'\le t\le\bar t$ the following conditions hold: 
  \begin{eqnarray}   \label{4.1}
  z(t)\in B-c_1^{-1}, \ -2c_6\sqrt{\eps}\le x(t)\le 2c_6\sqrt{\eps}+\eps, \,
   \frac12\eps^{r_1}\le y(t)\le 2c_6\sqrt{\eps}\,.
 \end{eqnarray}
According to  Lemma \ref{L4.1},  for $t'\le t\le t''$ we have $\dot y>\frac{1}{2}c_2^{-1}y$.
Therefore $t''-t'=O(\ln\eps)$ and $ z(t'')\in B-\frac32c_2^{-1}$. Then, again in correspondence with Lemma     \ref{L4.1}, on the segment $ x= 2c_6\sqrt{\eps}, \ |y|<2c_6\sqrt{\eps}$ we have  $ \dot x<-c_6c_2^{-1}\sqrt{\eps}$. On the segment $ x= -2c_6\sqrt{\eps}, \ |y|<2c_6\sqrt{\eps}$ we have  $ \dot x> c_6c_2^{-1}\sqrt{\eps}$. Therefore,  for $t'\le t\le t''$ conditions for $x(t)$ in (\ref{4.1}) hold as strict inequalities, and $|x(t'')|<2c_6\sqrt{\eps}$. In correspondence with the estimate $\dot y>\frac{1}{2}c_2^{-1}y$ we have $y(t'')>\eps^{r_1}$. As at the moment $t=t''$ the  phase point  should arrive to the boundary of the domain defined in (\ref{4.1}), there exists the only possibility: $y(t'')=2c_6\sqrt{\eps}$, i.e.  $(p(t''), q(t''))\in\tilde\alpha_1(z(t''))$.

In accordance with Lemma \ref{L4.1} for $t'\le t\le t''$ we have
 \begin{eqnarray*}    
  h(t)=-\omega_0(z(t))y(t)x(t)+O(\eps^{3/2}),  \  \nonumber 
    \frac{d}{dt}(\omega_0(z(t))y(t)x(t))=O(\eps^2 + |x|^3 + |y|^3) \,. \nonumber \
 \end{eqnarray*}
Therefore, $ h(t'')-h(t')=O(\eps^{3/2})$, and this was the assertion of   Lemma \ref{L4.3}. 

The proof of  Lemma  \ref{L4.4} follows   the same lines.
\medskip
 \subsubsection{Proof of  Lemma \ref{L4.5}} 
  \label{sss4.4.3}

Denote $\tilde\Lambda\circ t$ the time-$t$ schift of the set $\tilde\Lambda$ along  trajectories of the perturbed system  (\ref{perturbed}). Denote $\tilde\Lambda\circ [a_1, a_2]=\cup_{\substack{t\in [a_1, a_2]}}\tilde\Lambda\circ t$. Assume for simplicity of explanation that $\tilde\Lambda\circ [-K/\eps, 0]\in D$ (it is easy to avoid this restriction by considering only the  part of the set $\tilde\Lambda\circ [-K/\eps, 0]$ which does not leave $D$). For the set $\Xi$ we have
$$
\Xi=\left(\tilde\Lambda\circ\left [\frac{-\hat\tau_*-2c_8\delta}{\eps}, \frac{-\hat\tau_*+2c_8\delta}{\eps}\right]\right)\cap (D-c_1^{-1}).
$$
Consider $\Xi_0=\tilde \Lambda\circ [-d_1^{-1}, 0]$. For small enough $d_1^{-1}$ the points of this set lie in $B-\frac32c_1^{-1}$ with respect to $z$ and in $\frac12c_4^{-1}$-neighbourhood of the saddle point $C$ with respect to $(p, q)$. Therefore the set $\Xi_0$ can be considered in the variables $x, y, z$ of   Lemma \ref{L4.1}. As in these variables the phase flux through $\tilde\Lambda$ is $O(\eps^{r_1+1/2}\delta^{l-2})$, so the phase volume of the set $\Xi_0$ is  $O(\eps^{r_1+1/2}\delta^{l-2})$, if the volume element is $dx dy dz$. As the Jacobian of the transformation $x, y, z\mapsto p, q, z$ is $O(1)$, so the same estimate of the phase volume is valid if the volume element is $dp dq dz$. 
Consider now sets
 \begin{eqnarray*}    
  \Xi_j&=&\Xi_0\circ(-jd_1^{-1}),  \  \nonumber \\
  j&=&\left [\frac{(\hat\tau_*-2c_8\delta)d_1}{\eps}\right]-1, ..., \left [\frac{(\hat\tau_*+2c_8\delta)d_1}{\eps}\right]+1.  \nonumber \
 \end{eqnarray*}
We have mes $\Xi_j=O(\mbox{mes}\,\, \Xi_0)$, as the divergence of the right hand side of  system  (\ref{perturbed}) is $O(\eps)$, and so the distortion of the phase volume in this system during the time  $O(1/\eps)$ does not exceed $O(1)$ times. As $\Xi\subset\cup_{\substack{j}}\,\Xi_j$, we get mes $\Xi=O(\eps^{r_1-1/2}\delta^{l-1})=O(\eps^{r_1-1/2}\mbox{mes}\, W^\delta)$.
Lemma \ref{L4.5} is proved.
\medskip
 \section{Calculation of  measures  captured into different regions at the separatrix} 
  \label{S5}
  In this section  Propositions  \ref{Pr2.4}  and \ref{Pr2.5} of  Subsection  \ref{ss2.5} are proved.
\medskip
 \subsection{Preliminary constructions} 
  \label{ss5.1}
Denote $I=I(\nu, z, h)$ the value of the action variable for the trajectory $E=h$ of the unperturbed system in the region $G_\nu, \ \nu=1, 2, 3$. Since  in each of the regions $G_\nu$ the dependence of $I$ on $h$ is monotonous,  we can rewrite in any  
$G_\nu$ the averaged system  (\ref{2.3}) as a system of differential equations for $I, z$. Making use of gluing of solutions of averaged system at the separatrix (see Subsection \ref{ss2.3}), we can consider the averaged system with respect to the variable $\gamma=(\nu, z, I)$. Separatrix crossing leads to a jump of  values $\nu$ and $I$.
Denote $\bar g_1^\tau(\gamma)$ the shift of the point $\gamma$ for slow time $\tau=\eps t$ along the trajectory of such averaged system,  $\bar g_1^0(\gamma)=\gamma$, and index ``1'' indicates  that the solutions with $\nu=3$ and $\nu=1$ are glued at the separatrix.
Denote $W=W^{\delta}, \ W_{\nu}=W_\nu^{\delta}, \ \alpha=(p, q, z)$. Let $\Pi(\cdot)$ be the standard projection from the $\alpha$-space to the $\gamma$-space; $\Pi(\alpha)=\gamma$.
Denote
$$
\Gamma=\Pi [W], \ \Gamma_1^{K}=\bar g_1^{K} [\Gamma], W_1^{K}=\Pi^{-1} [\Gamma_1^K].
$$
We will denote points from $\Gamma, \Gamma_1^K, W, W_1^{K}$ as $\gamma^0, \gamma^K, \alpha^0=(p^0, q^0, z^0), \alpha^K=(p^K, q^K, z^K)$ respectively.
\medskip
 \subsection{Proof of the Proposition  \ref{Pr2.4}} 
  \label{ss5.2}
Denote $g^t(\alpha)$ the time-$t$ shift of a point $\alpha=(p, q, z)$  along the trajectory of  system  (\ref{perturbed}). Consider $g^t(\alpha^K)$ for $\alpha^K\in W_1^K$. Estimates in Propositions  \ref{perturbed},  \ref{2.2},  \ref{2.3} were proved for initial conditions from $W$.  Analogous estimates are valid for initial conditions from $W_1^{K}$, if we consider the motion in reverse  direction in time. From this  we get the following assertion: there exists a set $w_1\subset W_1^K, \ \mbox{mes} \ w_1=O(\eps^r\delta^{l-1})$, such that for initial conditions $
\alpha^K\in W_1^K\setminus w_1$  the behaviour of $h, z$ along  $g^t(\alpha^K), \ -K\eps\le t\le 0$, is described with an accuracy $O(\eps\ln\eps)$ by the solution of averaged system glued of solutions in $G_1$ and $G_3$ (passage from  $G_1$ to $G_2$ is impossible for $\alpha^K\in W_1^K\setminus w_1$).  It follows from these estimates that there exists a set ${W'}_1^{K}\subseteq W_1^{K}$ such that $g^{-K\eps^{-1}}[{W'_1}^K]\subseteq W, \  \mbox{mes} \,(W_1^K\setminus {W'_1}^K)=O(\eps\ln\eps \, \, \delta^{l-2})$. Denote $W'_1=g^{-K\eps^{-1}}[{W'_1}^K]$. The   definition of $W_2$ implies that 
$$
W'_1\subseteq W\setminus W_2=W_1\cup w
$$
A standard calculation of change of a phase volume along a motion gives us
 \begin{equation}   \label{5.1}    
   \mbox{mes} \, W'_1= \int\limits_{{W'_1}^K}\exp\left(-\eps\int\limits_0^{K\eps^{-1}}\left(\frac{\partial f_1}{\partial q}+\frac{\partial f_2}{\partial p}+\frac{\partial f_3}{\partial z}\right)dt\right)dp^K dq^K dz^K. \
   \end{equation} 
Here the outer  integral  is calculated with respect to initial (or, better to say, ``final'') conditions from the set ${W'_1}^K$. The inner integral
is calculated with respect to time along a solution of system (\ref{perturbed}) with given  ``final'' condition. In the  integral with respect to time it is reasonable to replace exact solution with the averaged one $\bar g_1^{-\tau}(\gamma^K)$, and to estimate the accuracy of this approximation. The result is described by the following Lemma, which is proved in Subsection \ref{sss5.4.1}.  
 \medskip
 \begin{lemma} 
        \label{L5.1}
For  $\alpha^K \in {{W'_1}^K} $ the  following estimate holds: 
 \begin{eqnarray*}    
\eps\int\limits_0^{K\eps^{-1}}\left(\frac{\partial f_1}{\partial q}+\frac{\partial f_2}{\partial p}+\frac{\partial f_3}{\partial z}\right)dt=\int\limits_0^K {\cal F}(\bar\gamma)d\tau+O(\eps\ln\eps),  \\
\mbox{\rm where} \ \bar\gamma=\bar g_1^{-\tau}(\gamma^K), \  \gamma^K=\Pi(\alpha^K), \\
{\cal F}(\gamma)=\frac{1}{T}\oint\limits_\gamma\left(\frac{\partial f_1^0}{\partial q}+\frac{\partial f_2^0}{\partial p}+\frac{\partial f_3^0}{\partial z}\right)dt. \ 
 \end{eqnarray*}
 The integral with index $\gamma=(\nu, Z, J)$ is calculated along the level line of the Hamiltonian $E$ in the region $G_{\nu}(z)$, and the ``action'' for this level line is equal to $J$. The parameter along the level line is time $t$ of the unperturbed motion, $T$ is the period of this motion.
 \end{lemma}  
 Making use of Lemma \ref{L5.1} we get
 \begin{equation}   
   \mbox{mes} \  W_1'= \int\limits_{{W_1'}^{K}}\exp\left(-\int\limits_0^K {\cal F}(\bar\gamma)d\tau\right)dp^K dq^K dz^K+O(\eps\ln\eps \ \delta^{l-2}) .\   \nonumber \
  \end{equation} 
 Replace  domain of integration ${W_1'}^K$ with $W_1^K$ in the last expression. This gives an additional error $O(\eps\ln\eps \ \delta^{l-2})$. We get
 \begin{equation}   
   \mbox{mes}  \ W_1'= \int\limits_{W_1^{K}}\exp\left(-\int\limits_0^K F(\bar\gamma)d\tau\right)dp^K dq^K dz^K
   +O(\eps\ln\eps \ \delta^{l-2}). \   \nonumber \
  \end{equation} 
  In the last expression we can use action-angle variables $J^K, \varphi^K$ instead of $p^K, q^K$ as independent variables. As integrand does not depend on $\varphi^K$, we get
   \begin{equation}   
   \mbox{mes}  \ W_1'=2\pi \int\limits_{\Gamma_1^{K}}\exp\left(-\int\limits_0^K{\cal{F}}(\bar\gamma)d\tau\right)dJ^K dz^K
   +O(\eps\ln\eps \ \delta^{l-2}). \   \nonumber \
  \end{equation} 
Let us make in the outer integral the transformation of  variables $(J^K, z^K)\mapsto(J^0, z^0)$ by means of the formula $\bar g_1^K(\gamma^0)=\gamma^K$, where $\gamma^0=(3, z^0, J^0), \ \gamma^K=(1, z^K, J^K)$. We get
   \begin{equation}   
   \mbox{mes}  \ (W_1')=2\pi \int\limits_\Gamma \Omega_{0, K}(\gamma^0)dJ^0 dz^0+O(\eps\ln\eps \ \delta^{l-2}), \   \nonumber \
  \end{equation} 
where
$$ 
\Omega_{0, K}(\gamma^0)=\exp\left(-\int\limits_0^K{\cal{F}}(\bar\gamma)d\tau\right)\frac{\partial (J^K, z^K)}{\partial (J^0, z^0)}, \quad \bar \gamma=\bar g_1^{\tau}(\gamma^0).
$$
 \medskip
\begin{lemma}
 \label{L5.2}
  $$   
\Omega_{0, K}(\gamma^0)=P_1(z_*)=\Theta_1(z_*)/\Theta_3(z_*),
 $$
 where $z_*=z_*(\gamma^0)$ is the value of $z$ at the moment of the separatrix crossing for the solution of averaged system
  $\bar g_1^{\tau}(\gamma^0)$.
 \end{lemma}
 This Lemma is proved in Subsection  \ref{sss5.4.2}.
 
 \medskip
Because of Lemma \ref{L5.2} 
\begin{equation}   \label{5.2}    
   \mbox{mes}  \ {W}_1'=2\pi \int\limits_\Gamma P_1(z_*)dJ^0dz^0+O(\eps\ln\eps \ \delta^{l-2}). \ 
      \end{equation}
    As ${W}_1'\subseteq W_1\cup w$  and  $\mbox{mes}  \ w = O(\eps^r \ \delta^{l-1})$, so
    $$
    \mbox{mes}  \ W_1' <  \mbox{mes}  \ W_1 + O(\eps^r \ \delta^{l-1}).
    $$
    In completely analogous way, but for the index $\nu=2$, we get that there exists a set ${W}_2' \subseteq W_2\cup w$  
     such that  
  \begin{eqnarray}   \label{5.3}    
   \mbox{mes}  \ {W}_2'&=&2\pi \int\limits_\Gamma P_2(z_*)\, d J^0dZ^0+O(\eps\ln\eps \ \delta^{l-2}), \\ 
     \mbox{mes}  \  {W}_2'&<& \mbox{mes} \  W_2+O(\eps^r \ \delta^{l-1}).  \nonumber
      \end{eqnarray}
      Consider relations 
  \begin{eqnarray}   \label{5.4}    
    \mbox{mes} \ {W}_1'&< &\mbox{mes} \ W_1+O(\eps^r \ \delta^{l-1}) ,\\
      \mbox{mes} \  {W}_2'&< &\mbox{mes} \ W_2+O(\eps^r \ \delta^{l-1}),\nonumber
  \\
        \mbox{mes} \ W_1&+&\mbox{mes} \  W_2=\mbox{mes}  \ W+O(\eps^r \ \delta^{l-1}), \nonumber
 \\
          \mbox{mes} \ {W}_1'&+&\mbox{mes} \  {W}_2'=\mbox{mes} \ W+O(\eps\ln\eps \ \delta^{l-2}). \nonumber
 \
      \end{eqnarray}
      To get the last equality, it is enough to add (\ref{5.2}) and  (\ref{5.3}), and to take into account that $P_1+P_2=1$. From  (\ref{5.4}) we get
  \begin{eqnarray*}   
    \mbox{mes} \ W_1&=& \mbox{mes} \ W- \mbox{mes} \ W_2+O(\eps^r\delta^{l-1})= \mbox{mes} \ {W}_1'+\mbox{mes}  \ {W}_2' -\mbox{mes} \ W_2+O(\eps\ln\eps\delta^{l-2}) \   \nonumber \\
&< &\mbox{mes}  \ {W}_1' +O(\eps\ln\eps \ \delta^{l-2}), \ \   \nonumber \
 \end{eqnarray*}
 and, on the other hand,
 $$
\mbox{mes} \ W_1>\mbox{mes} \  {W}_1'+O(\eps^r \ \delta^{l-1}).
$$
Therefore we get
 \begin{eqnarray*} 
 \mbox{mes} \ W_1&=&\mbox{mes} \ {W}_1'+O(\eps\ln\eps \ \delta^{l-2})=2\pi \int\limits_\Gamma P_1(z_*)\, dJ^0dz^0+O(\eps\ln\eps \ \delta^{l-2})=   \   \nonumber \\
 &=&\int\limits_W P_1(z_*)dp^0dq^0dz^0+O(\eps\ln\eps \ \delta^{l-2}),  \   \nonumber \\
 \end{eqnarray*}
 and analogous expression for $W_2$.
 Proposition \ref{Pr2.4} is proved.
\medskip
 \subsection{Proof of Proposition  \ref{Pr2.5}} 
  \label{ss5.3}
Let $I, \varphi \ \mbox{mod} \ 2\pi$ be the action-angle variables of the system with Hamiltonian $E$  in a neighbourhood of the set $W^\delta$. Consider in system (\ref{perturbed}) variable $\varphi$ as a new time. For $I, z$ we get a nonautonomous system of $l-1$ equations. Consider for this system extended phase space with space variables $I, z, \varphi$ and a new time $\vte\colon d\vte/d\varphi=1$. Now for $\varphi, I, z$ we have
  \begin{equation*}
\varphi'=1, \ I'=\eps f_4, \ z'=\eps f_5\,.
  \end{equation*} 
Here ``prime" denotes derivative with respect to $\vte, \ f_i=f_i(\varphi, I, z, \eps), \ j=4,5$, are smooth functions. Denote $u^\vte$ the operator of the shift along the trajectories of this system  during the time $\vte$.
Consider the sequence of the sets
$$
u^{2s\delta}(U^\delta), \ s=0, 1, ..., N-1, N=[\pi/\delta]\,.
$$
Only adjoining sets in this sequence can intersect each other. The measure of any such intersection is $O(\eps\delta^{l-1})$. Making use of the fact that the shift along the trajectories  of  this system  during the time $O(1)$ distorts measure only with a coefficient  $1+O(\eps)$, we get
 \begin{equation}   \label{5.5} 
\mbox{mes} \ \bigcup\limits_{s=0}^{N-1}u^{2s\delta}(U^{\delta})=N\,\mbox{mes}\, U^{\delta}+O(\eps\delta^{l-2}).
 \end{equation} 
Analogous reasoning  for sets $U_\nu^\delta, \, \nu=1,2$, gives the estimate
 \begin{equation}   \label{5.6} 
\mbox{mes} \ \bigcup\limits_{s=0}^{N-1}u^{2s\delta}(U_\nu^{\delta})=N\,\mbox{mes}\,U_\nu^{\delta}+O(\eps\delta^{l-2}).
 \end{equation} 
Then
 \begin{eqnarray}   \label{5.7} 
\mbox{mes} \ (W^{\delta}\!\!\vartriangle\!\bigcup\limits_{s=0}^{N-1}u^{2s\delta}(U^{\delta}))=O(\eps\delta^{l-2}+\delta^l), \\
\mbox{mes} \ (W_\nu^{\delta}\!\!\vartriangle\!\bigcup\limits_{s=0}^{N-1}u^{2s\delta}(U_\nu^{\delta}))=O(\eps\delta^{l-2}+\delta^l). \   \   \nonumber 
 \end{eqnarray} 
Here it is taken into account that $u^{\vartheta}(U_\nu^{\delta})\cap W^\delta\subseteq W_\nu^\delta\cup w; \ \   \vartriangle$ is the symbol of symmetric difference of sets.
From  (\ref{5.5}) -  (\ref{5.7}) we get
 \begin{eqnarray*}   
\mbox{mes} \ U^{\delta}= \frac{1}{N} \ \mbox{mes} \ W^{\delta}+O(\eps\delta^{l-1}+\delta^{l+1}), \\
\mbox{mes} \ U_\nu^{\delta}= \frac{1}{N} \ \mbox{mes} \ W_\nu^{\delta}+O(\eps\delta^{l-1}+\delta^{l+1}), \   \nonumber \\
\frac{\mbox{mes} \ U_\nu^{\delta}}{\mbox{mes} \ U^{\delta}}=\frac{\mbox{mes} \ W_\nu^{\delta}}{\mbox{mes} \ W^{\delta}}+O(\delta+\frac{\eps}{\delta}). \   \   \nonumber \\
 \end{eqnarray*} 
From here and from the result of  Proposition \ref{Pr2.4} we get
 \begin{equation}
\frac{\mbox{mes} \ U_\nu^{\delta}}{\mbox{mes} \ U^{\delta}}=\frac{\Theta_\nu(\hat z_*)}{\Theta_3(\hat z_*)}+O\left(\delta+\frac{\eps|\ln\eps|}{\delta}\right). \   \   \nonumber
 \end{equation}
 This was the assertion of  Proposition \ref{Pr2.5}.
\medskip
 \subsection{Proofs of  Lemmas on measure estimates}   
  \label{ss5.4}
 \subsubsection{Proof of  Lemma  \ref{5.1}}   
  \label{sss5.4.1}
Denote 
$$
\chi=\chi(p, q, z, \eps)=\frac{\partial f_1}{\partial q}+\frac{\partial f_2}{\partial p}+\frac{\partial f_3}{\partial z}, \quad \chi^0=\chi(p, q, z, 0).
$$
We have
 
 \begin{eqnarray}   \label{5.8}   
\eps\int\limits_0^{K\eps^{-1}}\chi dt&=&\eps\int\limits_0^{K\eps^{-1}}\chi^0dt+O(\eps)\\
&=&\eps\int\limits_0^{{t'}_-}\chi^0 dt+\eps\int\limits_{{t'}_-}^{{t'}_+}\chi^0dt+\eps\int\limits_{{t'}_+}^{K\eps^{-1}}\chi^0 dt+O(\eps), \  \nonumber 
 \end{eqnarray}
 where ${t'_-}$ and  ${t'_+}$ are the moments of the time analogous to the moments  ${t_-}, \ {t_+}$ introduced in Subsection \ref{ss2.5}, but for  initial conditions from  $W_1^K$. As  ${t'_-}-{t'_+}=O(\ln\eps)$, so the second term in the right hand side of the last equality is $O(\eps\ln\eps)$.
 
 For $0\le t\le{t'}_-$ we will consider motion round by round, as it was done in  Section \ref{S3}. Suppose for simplicity of the exposition that during any round the phase point crosses the ray $C\eta$ just one time, and that the motion takes place in the region $|E|<1/2$. Denote $t_1< t_2< ...<t_N$ the successive moments of the crossing  of  $C\eta, \ t_i\in(0, {t'}_-)$. According to  Lemma  \ref{L3.5} and its Corollary, we have
  
  \begin{equation}   \label{5.9} 
\frac{\eps}{t_{i+1}-t_i}\int\limits_{t_i}^{t_{i+1}}\chi^0 dt=\frac{\eps}{T(h(t_i), z(t_i))}\oint\limits_{E=h(t_i)}\chi^0dt+\eps^2 O\left(h^{-1}(t_i)\ln^{-1}h(t_i) \right).
 \end{equation} 
According to Propositions \ref{Pr2.1} - \ref{Pr2.3} 
 \begin{eqnarray}   \label{5.10}   
|h(t_i)-H_1(t_i)|+|z(t_i)-Z_1(t_i)|=O\left(\frac{\eps\ln\eps}{\ln H_1(t_i)}\right),  \\
 \frac12H_1(t_i)<h(t_i)<2H_1(t_i).\phantom{************} \  \nonumber \
 \end{eqnarray}
From  (\ref{5.9}),  (\ref{5.10}), making use of Lemma  \ref{L3.4}  and estimate  ({4}) of  Lemma  \ref{L3.2} we get
 \begin{equation}   
\frac{\eps}{t_{i+1}-t_i}\int\limits_{t_i}^{t_{i+1}}\chi^0 dt=\frac{\eps}{T(H_1(t_i), Z_1(t_i))}\oint\limits_{\substack{E=H_1(t_i)\\z=Z_1(t_i)}}\chi^0dt+O\left(\frac{\eps^2\ln\eps}{H_1(t_i)\ln^3H_1(t_i)}\right).  \  \nonumber \
 \end{equation} 
Let us  multiplay the left and right hand sides of this equality by \\$t_{i+1}-t_i=O(\ln H(t_i))$ and sum  up the obtained estimates. Taking into account that $t_1=O(1), \ t_N={t'}_-+O(\ln \eps)$, we get
 \begin{equation}   
\eps\int\limits_0^{t'_-}\chi^0 dt=\sum\limits_{i=1}^N\frac{1}{T(H_1(t_i), Z_1(t_i))}\,\, (\!\oint\limits_{\substack{E=H_1(t_i)\\z=Z_1(t_i)}}\chi^0dt \ \ )\ \eps(t_{i+1}-t_i)+O(\eps\ln\eps).  \  \nonumber \
 \end{equation} 
The sum in the right hand side can be represented as an integral with an accuracy $O(\eps\ln\eps)$. Therefore, we have
 \begin{equation}   
\eps\int\limits_0^{t_-'}\chi^0 dt=\int\limits_0^{\tau_*}{\cal{F}}(\bar\gamma)d\tau+O(\eps\ln\eps),  \  \nonumber \
 \end{equation} 
 where $\tau_*$ is the moment of the separatrix crossing in the averaged system.
 In the analogous way
  \begin{equation}   
\eps\int\limits_{{t'}_+}^{K\eps^{-1}}\chi^0 dt=\int\limits_{\tau_*}^K{\cal{F}}(\bar\gamma)d\tau+O(\eps\ln\eps).  \  \nonumber \
 \end{equation} 
Results of  Lemma  \ref{L5.1} follow from these estimates and (\ref{5.8}).
\medskip
 \subsubsection{Proof of  Lemma  \ref{5.2}}   
  \label{sss5.4.2}
 Let $\tau_3$ and $\tau_1$ be any numbers such that $ 0\le\tau_3<\tau_*<\tau_1\le K$. Here $\tau_*=\tau_*(\gamma^0)$ is the moment of the separatrix crossing for $\bar g_1^{\tau}(\gamma^0)$. Denote
 $$
 \bar g_1^{\tau_3}(\gamma^0)=\gamma^{(3)}=(3, z^{(3)}, J^{(3)}), \  \bar g_1^{\tau_1}(\gamma^0)=\gamma^{(1)}=(1, z^{(1)}, J^{(1)}).
 $$
 Each of values $\gamma^0, \gamma^{(3)}, \gamma^{(1)}, \gamma^K$ defines all others. These values can be defined also through $\tau_*, z_*$.
  Denote
 \begin{equation}   \label{5.11} 
\Omega_{\tau_3, \tau_1}=\exp\left(-\int\limits_{\tau_3}^{\tau_1}{\cal{F}}(\bar\gamma)d\tau\right)\frac{\partial(J^{(1)} , z^{(1)})}{\partial(J^{(3)} , z^{(3)})}. \
 \end{equation} 
 Similarly define  $\Omega_{0, \tau_3}$ and  $\Omega_{\tau_1, K}$. Then
 $$
 \Omega_{0, K}=\Omega_{0, \tau_3}\Omega_{\tau_3, \tau_1}\Omega_{\tau_1, K}\,.
 $$
 \medskip
 \begin{lemma} 
        \label{L5.3}  
  $
  \Omega_{0, \tau_3}=1, \ \Omega_{\tau_1, K}=1\,.
  $
 \end {lemma} 
 This Lemma is proved in Subsection \ref{sss5.4.3}.
 
  \medskip
  \noindent
  
  \begin{corollary}
 
$ \Omega_{0, K}=\Omega_{\tau_3, \tau_1}$\,.
\end{corollary}

 Now let $\tau_3$ and $ \tau_1$ tend to $\tau_*$. The first multiplier in (\ref{5.11}) tends  to $1$. Therefore 
  
 \begin{equation} 
 \Omega_{0, K}=\lim_{\substack{\tau_3\to\tau_*-0 \\ \tau_1\to\tau_*+0}}\frac{\partial(J^{(1)} , z^{(1)})}{\partial(J^{(3)} , z^{(3)})} 
=\left.\lim_{\tau_1\to\tau_*+0 }\frac{\partial(J^{(1)} , z^{(1)})}{\partial(\tau_* , z_*)}\right / \lim_{\tau_1\to\tau_*-0}\frac{\partial(J^{(3)} , z^{(3)})}{\partial(\tau_* , z_*)}\,.
   \nonumber 
 \end{equation} 
Let us calculate these limits. Denote $F$ and $\Phi$ the right hand sides of the averaged equations for $I$ and  $z$ respectively. Then
 \begin{eqnarray*}    
  J^{(1)}=\frac{1}{2\pi}S_1(z_*)+\int\limits_{\tau_*}^{\tau_1}Fd\tau\, , \quad
z^{(1)}=z_*+\int\limits_{\tau_*}^{\tau_1}\Phi d\tau \, .  \nonumber \
 \end{eqnarray*}
Integrals here are calculated along the solution of the averaged system  $\bar g_1^{\tau}(\gamma^0)$, and $S_1=S_1(z)$ is the area of the region $G_1(z)$. Making use of the formulas  for the right hand sides of averaged system (\ref{2.3}), (\ref{2.7}), we get 
 \begin{eqnarray*} 
 \lim_{\tau_1\to\tau_*+0 }\frac{\partial J^{(1)}}{\partial \tau_*}&=&-\lim_{\tau_1\to\tau_*+0}F=\frac{1}{2\pi}\left(\Theta_1(z_*)+\left(\oint\limits_{l_1}\frac{\partial E}{\partial z}dt \right)\,  f_{3c}^0 \right), \\
 \lim_{\tau_1\to\tau_*+0 }\frac{\partial z^{(1)}}{\partial \tau_*}&=&-\lim_{\tau_1\to\tau_*+0}\Phi=-f_{3c}^0 , \\
   \lim_{\tau_1\to\tau_*+0 }\frac{\partial J^{(1)}}{\partial z_*}&=&\frac{1}{2\pi}\frac{\partial S_1}{\partial z_*}=-\frac{1}{2\pi}\oint\limits_{l_1}\frac{\partial E}{\partial z}dt,  \\
\lim_{\tau_1\to\tau_*+0 }\frac{\partial z^{(1)}}{\partial z_* }&=&{\bf 1}_{l-2} . \
 \end{eqnarray*}  
Here ${\bf 1}_{l-2}$ is the unit $(l-2)\times(l-2)$ matrix;  $f_{3c}^0$ is the value of the function  $f_3^0$ at the saddle point $C$ for $z=z_*$. Making use of these relations, we get
 \begin{equation}
\lim_{\tau_1\to\tau_*+0}\frac{\partial(J^{(1)} , z^{(1)})}{\partial(\tau_*, z_*)}= 
\begin{vmatrix}
\frac{1}{2\pi}\left(\Theta_1(z_*) +\left(\oint\limits_{l_1}\frac{\partial E}{\partial z}dt \right)\,  f_{3c}^0 \right),
  & -\frac{1}{2\pi}\oint\limits_{l_1}\frac{\partial E}{\partial z}dt \\
-f_{3c}^0 , & {\bf 1}_{l-2} \\
\end{vmatrix}
=\frac{1}{2\pi}\Theta_1(z_*) \, .  \nonumber \
 \end{equation}  
Similarly, 
 \begin{equation}  
 \lim_{\tau_3\to\tau_*-0 }\frac{\partial (J^{(3)}, z^{(3)})}{\partial  (\tau_*, z_*)}=\frac{1}{2\pi}\Theta_3(z_*)   \, . \nonumber 
 \end{equation} 
 Finally we have
 $$
 \Omega_{0, K}=\Theta_1(z_*)/\Theta_3(z_*)\,.
 $$
 This  was the assertion of  Lemma \ref{5.2}.
 \medskip
 \subsubsection{Proof of  Lemma  \ref{5.3}}   
  \label{sss5.4.3}
To avoid long calculations with derivatives we will use known results about the averaging method.
Denote $\Gamma_0$  a neighbourhood of the point $(3, z^0, J^0)$ such that  $\bar g_1^{\tau}(\gamma)$ does not cross the separatrix for $0\le\tau\le \tau_1$ and $\gamma\in\Gamma_0$. Let $W_0=\Pi^{-1}(\Gamma_0)$, where $\Pi$ is the standard projection from $p, q, z$-space to $\gamma=(\nu, z, J)$-space. The reasoning  of the Subsection \ref{5.1} shows  that there exists a set $W'_0\in W_0$ such that
 \begin{eqnarray*}   
    \mbox{mes} \ {W_0}'&= &\mbox{mes} \ W_0+O(\eps), \\
      \mbox{mes} \  {W_0}'&=&2\pi\int\limits_{\Gamma_0}\Omega_{0, \tau_3}dJ^0dz^0+O(\eps) .\
       \end{eqnarray*}
       From here
        $$
         \mbox{mes} \  {W_0}-2\pi\int\limits_{\Gamma_0}\Omega_{0, \tau_3}dJ^0dz^0=O(\eps).
         $$
         As the left hand side does not depend on $\eps$, so
         $$
        \mbox{mes} \  {W_0}=2\pi\int\limits_{\Gamma_0}\Omega_{0, \tau_3}dJ^0dz^0=\int\limits_{W_0}\Omega_{0, \tau_3}dp^0dq^0dz^0 .  
         $$
 Therefore $\Omega_{0, \tau_3}=1$. Similarly,   $\Omega_{\tau_1, K}=1$.      
 Lemma  \ref{5.3} is proved.        
  
 \medskip
 \subsection{A rule for calculation of probabilities}   
  \label{ss5.5}
 A  heuristic reasoning  of \cite {lif, goldr}, which leads to the formulas of  Section \ref{S2} for probabilities of capture into different regions, is exposed in this subsection. This reasoning can be used as, in some sense, a rule, as it allows to calculate probabilities in general case of systems   of  form  (\ref{perturbed}), for other than in Fig. \ref{unperturbed_plane} types of phase portraits. This reasoning is justified by Proposition \ref{Pr5.1} of this subsection.
 \medskip
 \subsubsection{ A scheme of calculation of probabilities}   
  \label{sss5.5.1}
The following reasoning does not pretend to be rigorous. A corresponding  rigorous assertion is formulated at the end of this subsection. 
  \medskip
  Let $C\xi\eta=C\xi\eta(z)$ be the system of principal axes for the saddle point $C$ oriented as in  Fig. \ref{unperturbed_plane}. Let a phase point $(p(t), q(t))$ start moving  at a  moment of time $t=t'$  with $z=z', \ E=h'$ from the ray $C\eta$. Denote $\Theta_\nu=\Theta_\nu(z')$. The point $(p(t), q(t))$ first makes a curve $l_2'$, which  is close to $l_2(z)$. At the end of this curve
  $$
  E=h''=h'+\int\limits_{l'_2}\frac{dE}{dt}dt\approx h'+\int\limits_{l_2}\frac{dE}{dt}dt=h'-\eps\Theta_2\,.
  $$
If $0<h'<\eps\Theta_2$, then $h''<0$, i.e. the phase point  is captured into the region $G_2$. If $h'>\eps\Theta_2$, then  in further motion the phase point  makes a curve $l'_1$, which  is close to $l_1(z')$. At the end of this curve $E=h'''\approx h'-\eps(\Theta_1+\Theta_2)$. If $\eps\Theta_2<h'<\eps(\Theta_1+\Theta_2)$,  then $h'''<0$, i.e. the phase point  is captured into the region $G_1$. If $h'>\eps(\Theta_1+\Theta_2)$,  then $h'''>0$, i.e. the phase point comes back to $C\eta$. Introduce intervals $\varkappa_1=(\eps\Theta_2, \ \eps(\Theta_1+\Theta_2)),  \ \varkappa_2=(0, \eps\Theta_2, ), \ \varkappa_3=( \eps(\Theta_1+\Theta_2), \infty), \ \varkappa_\nu=\varkappa_\nu(z')$. In accordance with the previous explanation, points with $h'\in\varkappa_3$ will come back to $C\eta$ and after several rounds will arrive to $\varkappa_1\cup\varkappa_2$.  Points with $h'\in\varkappa_\nu, \ \nu=1,2$, will be captured into $G_\nu$. The measure of the subset of $U^\delta$, which will be captured into $G_\nu$, is proportional to the length of the interval $\varkappa_\nu(\hat z_*)$ (because the majority of  points from  $U^\delta$ have $z\approx\hat z_*$ when cross $C\eta$ for the  last time, and the phase flux through $\varkappa_\nu$ is equal in the principal approximation to   the length of  $\varkappa_\nu$).
Thus
$$
Q_\nu(\hat M_0)=\frac{\mbox{length} \ \varkappa_\nu}{\mbox{length} \ \varkappa_1+\mbox{length} \ \varkappa_2}=\frac{\Theta_\nu(\hat z_*)}{\Theta_1(\hat z_*)+\Theta_2(\hat z_*)}\,.
$$
 
The following  assertion  corresponds to the previous reasoning.

\begin{figure}
 \begin{center}
            \includegraphics[scale=0.5, angle =0]{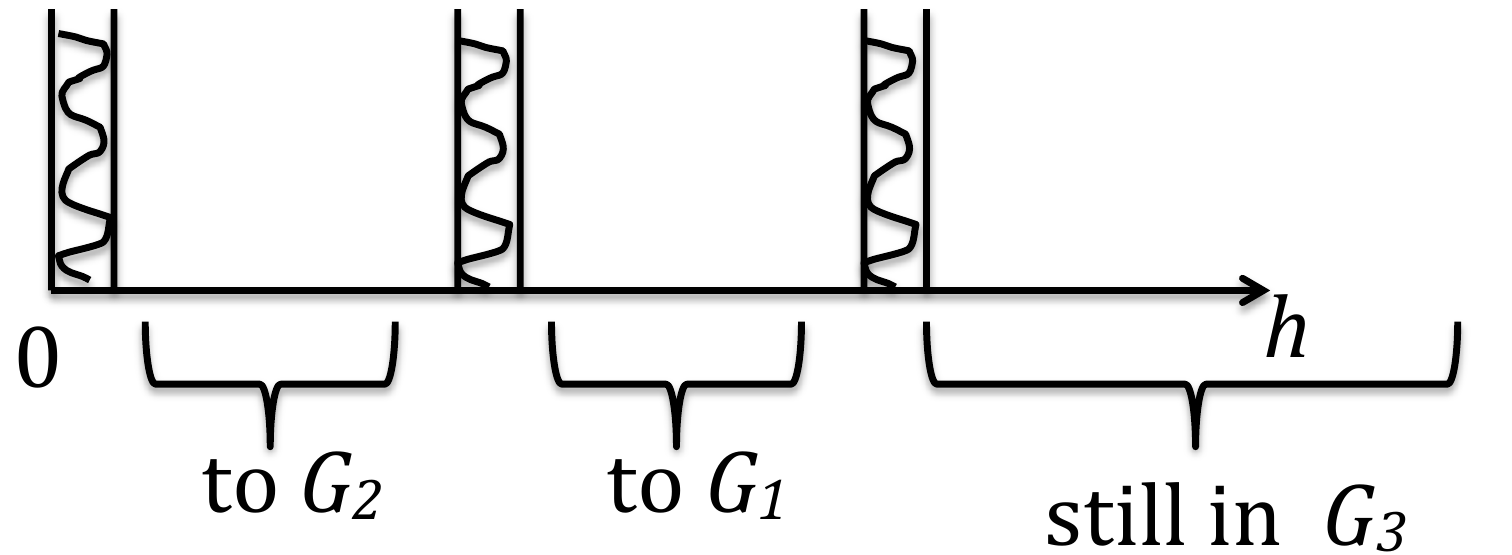}
            \end{center}
            
           \caption{For calculation of probabilities. }
            \label{intervals} 
\end{figure}

 \medskip
\begin{proposition}  
  \label{Pr5.1}
Let at a moment of time $t'$  a point $(p(t'), q(t'))$ lie on the axis  $C\eta(z')$ in $k_3^{-1}$-neighbourhood of the point $C$, and $z'\in B-k_4^{-1}, \ h'=E(p(t'), q(t'), z(t'))$.
Introduce intervals (Fig. \ref{intervals})
  \begin{eqnarray*}  
\varkappa'_1&=&(\eps\Theta_2+k_5\eps^{3/2}, \ \eps(\Theta_1+\Theta_2)-k_5\eps^{3/2}), \\
 \varkappa'_2&=&(k_5\eps^{3/2},  \eps\Theta_2-k_5\eps^{3/2}) , \\
  \varkappa'_3&=&(\eps(\Theta_1+\Theta_2)+k_5\eps^{3/2}, k_6^{-1})\,. \ 
 \end{eqnarray*}  
  Then the following holds.
 
 $1^0$. If $h' \in\varkappa'_\nu, \ \nu=1, 2$ then the phase point does not cross  $C\eta$ again for $t>t'$. There exists $t_\nu=t'+O(\ln\eps)$ such that $(p(t_\nu), q(t_\nu))\in G_\nu(z(t_\nu)), \ h(t_\nu)=-k_2\eps$.
 
  $2^0$. If $h' \in\varkappa'_3$, then there exists $t_3>t'$ such that $(p(t_3), q(t_3))\in C\eta(z(t_3)), \ h'-h(t_3)>k_7^{-1}\eps$. 
 \end{proposition}
 
    \medskip 
 {\bf {Remarks}}.
 
 {\bf 1}. Constant $k_2$ was introduced in Proposition \ref{perturbed}.

 {\bf 2}. Making use of Proposition \ref{Pr5.1} it is possible to prove formula for the probability (\ref{eq:cor}).  It is possible to prove also the result analogous to Proposition \ref{Pr2.4} but with more rough estimate: $O(\sqrt\eps)$ instead of  $O(\eps\ln\eps)$. 
 
{\bf 3}. The proof of the formula for the probability in previous sections uses essentially  that separatrices divide phase space into three regions. But sometimes, because of an additional symmetry of the problem, a system has several saddle points connected by separatrices, and in these cases phase space can be divided into four or more regions (see, for example, \cite{nei_mtt}). Formulas for probabilities for these cases can be obtained by means of a reasoning, analogous to that used at the beginning of this section. These formulas can be justified by means of an assertion analogous to  Proposition \ref{Pr5.1}.
 
  \medskip
 \subsubsection{Proof of  Proposition \ref{Pr5.1}}   
  \label{sss5.5.2}
Let us restrict ourselves by proving  of assertion $1^0$ for $\nu=2$. The proofs of other assertions are  completely analogous.
Denote
$$
\chi=\chi(p, q, z, \eps)=\frac{\partial E}{\partial q}f_1+\frac{\partial E}{\partial p}f_2+\frac{\partial E}{\partial z}f_3\, .
$$
Choose any $k_4$ such that $B-k_4^{-1}$ is not empty. For $z\in B-\frac12k_4^{-1}$ we have
$$
\oint\limits_{l_2}|\chi|dt<c_1\,.
$$
We may assume that  the quadratic part of the Hamiltonian $E$ near the saddle point $C$  in  the variables $\eta, \xi$   has the form $\frac12\omega_0(z)(\eta^2-\xi^2), \ \omega_0>0$. Denote $x=(\xi-\eta)/\sqrt2, \ y= (\xi+\eta)/\sqrt2$ (cf. Fig. \ref{square+trajectories}). Let $\Psi=\Psi(y, x, z)$ be  the Hamiltonian $E$, expressed through $y, x, z$:
$$
\Psi=-\omega_0yx+O(|y|^3+|x|^3)\,.
$$
 \medskip
\begin{lemma}  
  \label{L5.4}
Let $z\in B-\frac12k_4^{-1}, \ |x|<d_1^{-1}, |y|<d_1^{-1}$.
If $|y|\ge |x|$, then
$$
 \frac{\partial \Psi}{\partial x}=-\omega_0y+O(y^2)\,.
 $$
If $|y|< |x|$, then
$$
 \frac{\partial \Psi}{\partial y}=-\omega_0x+O(x^2)\,.
 $$
\end{lemma}  
The proof is evident.

   \medskip 
 \begin {corollary} For $0<|y|<d_2^{-1}<d_1^{-1}, \ |h|<d_3^{-1}, \   |x|\le y$  the equation $\Psi(y, x, z)=h$ defines a unique  $x=\tilde x(y, h, z)$ such that $ |x|<d_1^{-1}$. Function $\tilde x$ is smooth and
 $$
 \frac{\partial \tilde x}{\partial h}= \frac{1}{\partial \Psi/\partial x}, \   \frac{\partial \tilde x}{\partial z}= -\frac{\partial \Psi/\partial z}{\partial \Psi/\partial x}\,.
 $$
 For $0< |x|<d_2^{-1}, \  |h|<d_3^{-1}, \  |y|\le |x|$ the equation $\Psi=h$ defines  in the analogous manner   $y=\tilde y(x, h, z)$.
 \end{corollary}
 
 Denote $y(t), \ x(t)$ the values of $y, x$ at the point $p(t), q(t), z(t)$. We have $y(t')=-x(t')>0$. Let $\eps^{3/2}<h(t')<\eps\Theta_2(z')$. Denote $t_{1*}$ the supremum of  moments of time $t_1>t'$ such that for $t'\le t\le t_1$ the solution  $p(t), q(t), z(t)$ is defined and meets the conditions
 \begin{eqnarray}  \label{5.12}
|x(t)|&<&d_1^{-1}, \ 0< y(t)<d_2^{-1}, \ z(t)\in B-\frac34k_4^{-1}\,, \\
 |y(t)|&>&|x(t)|, \ y(t)>\eps^{7/8} \  \ \nonumber 
 \end{eqnarray} 
(the set of such $t_1$ is not empty, as $\dot y(t')>-\dot x(t')$\,). Denote $y_{1*}=y(t_{1*})$. 

For  $t'\le t\le t_{1*}$ we have
  \begin{eqnarray}  \label{5.13}
 \dot y&=&\omega_0y+O(y^2)+O(\eps)>c_2^{-1}y\,, \\
t_{1*}-t'&=&\int\limits_{y(t')}^{y_{1*}}\frac{dy}{\dot y}<c_2\int\limits_{y(t')}^{y_{1*}}\frac{dy}{y}<c_2\ln\frac{y_{1*}}{y(t')}=O(\ln\eps) \,,   \ \nonumber \\
 z(t)-z(t')&=&O(\eps\ln\eps)\,,   \ \nonumber  \\
  |h(t)-h(t')|&\le&\eps\int\limits_{t'}^{t}|\chi|dt<c_2\eps\int\limits_{y(t')}^{y(t)}\frac{|\chi|dy}{y}=\eps\int\limits_{y(t')}^{y(t)}O(1)dy=O(\eps)\,.   \ \nonumber  \
   \end{eqnarray}
The obtained inequalities show that at the moment of time  $t_{1*}$ all conditions in (\ref{5.12}) but the inequality  $y(t)<d_2^{-1}$ are satisfied with some margins. Therefore $y(t_{1*})=d_2^{-1}$. 

The obtained inequalities allow to give more accurate estimate of $h(t_{1*})$:
   \begin{eqnarray*}  
  h(t_{1*})&-&h(t')=\eps\int\limits_{t'}^{t_{1*}}\chi dt=
  \eps\int\limits_{y(t')}^{d_2^{-1}}\frac{\chi  dy}{-\partial\Psi/\partial x+O(\eps)}\\
  &=& \eps\int\limits_{y(t')}^{d_2^{-1}}
  \left(-\frac{\chi}{\partial\Psi/\partial x}\right)
  _{\substack{E=0 \\ z=z' \\ \eps=0}}dy
  -\eps\int\limits_{y(t')}^{d_2^{-1}}\left[\left(\frac{\chi}{\partial\Psi/\partial x+O(\eps)}\right)_{\substack{E=h(t) \\z=z(t)}}-\left(\frac{\chi}{\partial\Psi/\partial x}\right)_{\substack{E=0 \\ z=z' \\ \eps=0}}\right]dy \,.  \ \nonumber  \
   \end{eqnarray*}
 By means of   (\ref{5.12}), Lemma \ref{L5.4} and its Corollary, the integrand in the second integral for $y>\sqrt\eps$ is estimated as
 $$
 O\left(\frac{\eps}{y^2}\right)+O\left(\frac{\eps\ln\eps}{y}\right)\,.
$$
For $\eps^{7/8}<y< \sqrt\eps$ the integrands in both integrals are $O(1)$.  Therefore
$$
 h(t_{1*})-h(t')=-\eps\int\limits_0^{d_2^{-1}}
  \left(\frac{\chi}{\partial\Psi/\partial x}\right)
  _{\substack{E=0 \\ z=z(t') \\ \eps=0}}dy+O(\eps^{3/2})\,.
$$
Here the integral is calculated over a segment of the unperturbed separatrix.

In the further motion the phase point makes a curve situated in $O(\eps\ln\eps)$-neighbourhood of the unperturbed separatrix for $z=z'$ and arrives at the segment $x=d_2^{-1}, |y|<d_1^{-1}$ at a moment of time $t_{2*}=t_{1*}+O(1)$ having $E=O(\eps)$. The change of $E$ along this curve with an accuracy $O(\eps^2\ln\eps)$ is equal to the integral of function $\eps\chi(p, q, z', 0)$ along the corresponding part of the unperturbed separatrix.

Then, through the time  $O(\ln\eps)$, at some moment of time $t_{3*}$, the phase point arrives either at the ray $x=y+c_3\eps>0$ or at the ray $x=-y+c_3\eps>0$. The motion for $t\in(t_{2*}, t_{3*})$ is considered in completely analogous way to that for $t'\le t\le t_{1*}$. The change of $E$ with an accuracy $O(\eps^{3/2})$ is equal to the integral of  $\eps\chi(p, q, z', 0)$ along the segment of the unperturbed separatrix with $0\le x\le d_2^{-1}, \  |y|<d_1^{-1}$. Therefore $h( t_{3*})=h(t')-\eps\Theta_2(z')+O(\eps^{3/2})$.

If
$$
k_5\eps^{3/2}<h(t')<\eps\Theta_2(z')-k_5\eps^{3/2}\,,
$$
then
$$
-\eps\Theta_2(z')+\frac12k_5\eps^{3/2}<h( t_{3*})<-\frac12k_5\eps^{3/2}
$$
Therefore at the moment of time $t_{3*}$ the phase point lies in the region $G_2(z(t_{3*}))$. Condition $h( t_{3*})<-\frac12k_5\eps^{3/2}$ and Lemma \ref{5.4} allow to estimate $\dot \eta$ from below by a value of order $\eps^{3/4}$. Making use of this estimate we can show that  
at a moment of time $t''=t_{3*}+O(\eps^{1/4})$ the phase point arrives at the ray $C\xi$ having $h(t'')=O(\eps), \ h( t'')<-\frac14k_5\eps^{3/2}$.

 Further motion is considered in an analogous manner. While $-k_2\eps<h(t)<2c_1\eps$ the phase point moves round by round making curves near the unperturbed separatrix $l_2$. One round takes time $O(\ln\eps)$, the value of $E$ during one round decays by $\eps\Theta_2(z')+O(\eps^{3/2})>c_3^{-1}\eps$.
Therefore, there exists a moment of time $t_\nu=t'+O(\ln\eps)$ such that $h(t_\nu)=-k_2\eps, \ (p(t_\nu), q(t_\nu))\in G_\nu(z(t_\nu))$. This is the assertion of  Proposition \ref{Pr5.1}.

\bigskip {\bf Acknowledgment.} The author is thankful  to N.R.Lebovitz for  comments, discussions, and  help, to A.Bolsinov for advices  on integrable systems.

\newpage
\appendix
\section{Appendix. Perturbations of polyintegrable  systems and separatrix crosings}
The goal of this Appendix is to give a general description  of the problem of separatrix crossing
in single-frequency  systems and to demonstrate,  that  under rather general assumptions the study of this problem can be reduced to study of separatrix crossing in  system (\ref{perturbed}).  The exposition here follows mainly \cite{akn}, Subsection 6.1.10. 

A natural framework for studying one-frequency averaging is the framework of perturbations of polyintegrable superintegrable (also called  Nambu) systems\,\footnote{See \cite{superintegrable} for description of properties of polyintegrable systems.}. In this problem the equations of  motion have the form
\begin{equation*}
\label{A1}
\dot x= v(x, \eps),\quad x\in D \subseteq  \R^l, \quad 0<\eps\ll 1,\quad v(x, \eps)=v_0(x)+v_1(x, \eps)\,. \hskip 0.7cm (A.1)
\end{equation*}
Here  $D$ is a bounded domain in $\R^l$. We assume that the unperturbed ($\eps = 0$) system is polyintegrable, i.e. it  has $l-1$ smooth first  integrals $H_1,\ldots, H_{l-1}$ which are independent almost everywhere in $D$. 
We assume that the domain $D$ contains, together with each point, also the entire connected  component of the  common level set (a level line) of the integrals passing through this point. Then a   level line on which the first integrals are independent is  a smooth closed curve. In any  domain filled by  such  level lines system ({A.1}) can be reduced  to the standard form of a  system with one rotating phase.  

To introduce a framework for separatrix crossing we assume that: 

\noindent
a) the rank of the Jacobi matrix of the map ${\cal H}\, : \, D \to  \R ^{l -1}$ given by ${\cal H} (x) =
(H_1(x),\ldots,  H_{l-1}(x))$ is equal to $ l-1$  everywhere  but on a smooth  $l-2$ dimensional surface, where it equals to $ l-2$
;

\noindent
b) at each point,  where the rank equals $l-2$, the restriction of one of integrals
onto the  joint level of other integrals has a non-degenerate critical point;

\noindent
c) at  equilibrium positions of the unperturbed system ({A.1})  two eigenvalues
are non-zero real numbers (the other eigenvalues are equal to 0 because of the existence
of the integrals).

Then points, where the rank of the map  ${\cal H}$ equals $l-2$,  coincide with equilibria of system  ({A.1}) for $\eps=0$, the sum of non-zero eigenvalues equals 0 for such an equilibrium. We  call separatrices the  common level lines that pass through these points 
 as well as a union of such level lines.  Under the action of the perturbation    phase points can cross  separatrices.
 
 Assume that functions  $H_1,\ldots, H_{l-2}$ are independent on separatrices.   The values  $z_1,\ldots, z_{l-2}$ of these functions from some ball in $\R^{l-2}$ can be taken as new variables.  Joint levels of these functions form $l-2$-parametric family of 2-dimensional  surfaces ${\cal{S}}_z$, $z=( z_1,\ldots, z_{l-2})$. Unperturbed dynamics on each of these surfaces is described by a Hamiltonian system with one degree of freedom for which the restriction $E$ of the  function $H_{l-1}$ onto this surface  is  Hamilton's function, but the symplectic structure may be non-canonical.  The phase portrait of each of these systems contains  a saddle point and passing through it separatrices. In a neighbourhood of separatrices  the phase portrait has the same form as in Fig. \ref{unperturbed_plane} and can be considered as a portrait in $\R^2$. (Notice that this does not depend on topology of  ${\cal{S}}_z$.    For example, the phase portrait of the pendulum, Fig. \ref{pendulum}, should be considered on a cylinder, but a neighbourhood of separatrices can be put in $\R^2$ as a neighbourhood of separatrices of the form shown  in Fig. \ref{unperturbed_plane}.)  
 
 Let $\tilde p, q$ be Cartesian coordinates in   $\R^2$. In these coordinates, in a neighbourhood of separatrices, the symplectic structure has a form
 $\mu( \tilde p, q,z) d\tilde p\wedge dq$,  \, $\mu( \tilde p, q,z)\ne 0$. Define in a neighbourhood of separatrices a function $p=p( \tilde p, q,z)$ such that $\partial p/\partial \tilde p= \mu( \tilde p, q,z)$. In the variables $p,q$ the symplectic structure takes  the canonical  form  $d p\wedge dq$, and equation ({A.1}) takes the form (\ref{perturbed}). Thus, the results in Subsection \ref {s2.est} for system (\ref{perturbed}) describe also separatrix crossing for ({A.1}).
 
 One can also consider separatrix crossings directly for perturbations of a polyintegrable system. The phase space of the averaged system is the set of common level lines of the integrals of the unperturbed system, which has
the natural structure of a manifold with singularities \cite{bf} (singularities correspond to a separatrix). The averaged system approximately describes the evolution of the slow variables -   values of the integrals of the unperturbed system. The probabilities of falling into different domains after a separatrix crossing are expressed in terms of ratios of the quantities
 $$
 \tilde \Theta_i(z)=-\oint\limits_{l_i(z)}\left(\beta_1(z)\frac{\partial H_1}{\partial x}+\ldots+ \beta_{l-1}(z)\frac{\partial H_{l-1}}{\partial x}\right)v_1(x,0)dt \,,
 $$
 where $z$ parametrises the surface of singular points (``saddles'')  of the unperturbed system,  $\beta_j$ are coefficients  such that the expression inside the parentheses in the integrand vanishes at  singular points, and $l_i=l_i(z)$ is a separatrix.

   \end{document}